\pdfoutput=1
\RequirePackage{fix-cm}
\documentclass[smallextended]{svjour3}
\smartqed
\usepackage[top=2.5cm, bottom=2.5cm, left=2.5cm, right=2.5cm]{geometry}
\usepackage{amssymb,amsmath,mathtools,etoolbox,bm,booktabs,graphicx,float,enumitem}
\usepackage{tikz,overpic}
\usepackage{subfig}
\usepackage[numbers,sort&compress]{natbib}
\definecolor{lightblue}{rgb}{0.22,0.45,0.70}
\usepackage[colorlinks=true,breaklinks=true,linkcolor=lightblue,citecolor=lightblue,urlcolor=black]{hyperref}
\makeatletter\let\cl@chapter\relax\makeatother 
\usepackage[capitalise,nameinlink]{cleveref}

\usepackage[acronym]{glossaries}

\newacronym{bsk}{BSK}{Bazant--Storey--Kornyshev}
\newacronym{pb}{PB}{Poisson--Boltzmann}
\newacronym{pf}{PF}{Poisson--Fermi}
\newacronym{pnp}{PNP}{Poisson--Nernst--Planck}
\newacronym{fem}{FEM}{finite element methods}
\newacronym{ipam}{iPAM}{improved polygonal area mapping}

\numberwithin{equation}{section}
\numberwithin{theorem}{section}
\numberwithin{figure}{section}
\numberwithin{table}{section}
\numberwithin{remark}{section}
\numberwithin{lemma}{section}
\numberwithin{corollary}{section}
\numberwithin{proposition}{section}

\renewenvironment{proof}{\noindent{\it Proof.}}{\hfill$\square$}
\newenvironment{myproof}{\par\noindent}{\hfill$\square$\par}
\allowdisplaybreaks
\newcommand\vdiv{\mathop{\mathrm{div}}\nolimits}
\newcommand{\bx}{\boldsymbol{x}}
\newcommand{\by}{\boldsymbol{y}}
\newcommand{\bp}{\boldsymbol{p}}
\newcommand{\bq}{\boldsymbol{q}}
\newcommand{\bxi}{\boldsymbol{\xi}}
\newcommand{\bzeta}{\boldsymbol{\zeta}}
\newcommand{\btau}{\boldsymbol{\tau}}
\newcommand{\bY}{\mathbf{Y}}
\newcommand{\bH}{\mathbf{H}}
\newcommand{\bL}{\mathbf{L}}
\newcommand{\bQ}{\mathbf{Q}}
\newcommand{\bPi}{\mathbf{\Pi}}
\newcommand{\ds}{\ensuremath{\, \mathrm{d}s}}

\newcommand*{\normiii}[1]{\left|\!\left|\!\left|\ifblank{#1}{\:\cdot\:}{#1}\right|\!\right|\!\right|}

\begin{document}
\titlerunning{Mixed-primal FEM for a Poisson--Fermi interface problem}
\authorrunning{Mandottil, Nataraj \& Ruiz-Baier}

\title{A mixed-primal FEM for  Poisson--Fermi interfacial problems\thanks{\textbf{Funding:} This work has been partially supported by the IITB--Monash Research Academy; J.C. Bose grant ANRF/JBG/2025/000209/HAA; and by the Australian Research Council through the \textsc{Future Fellowship} grant FT220100496.}}
\author{
Siyad Rahman Mandottil
\and
Neela Nataraj
\and
Ricardo Ruiz-Baier
}
\institute{S. R. Mandottil \at IITB--Monash Research Academy, Indian Institute of Technology Bombay, Mumbai, Maharashtra 400076, India. \email{siyadrahmanm@monash.edu}.
\and N. Nataraj \at Department of Mathematics, Indian Institute of Technology Bombay, Mumbai, Maharashtra 400076, India. \email{neela@math.iitb.ac.in}.
\and R. Ruiz-Baier (corresponding author) \at School of Mathematics, Monash University, 9 Rainforest Walk, Melbourne 3800 VIC, Australia. \email{ricardo.ruizbaier@monash.edu}.
}
\date{\today}
\maketitle

\begin{abstract}
This paper presents a mixed-primal finite element method for a coupled, linearised Poisson--Fermi interface problem arising from the Bazant--Storey--Kornyshev theory, which describes electrostatic interactions in ionic solutions. It leads to a transmission problem coupling a fourth-order equation in the solvent region with a second-order equation in the solute region. To precisely compute the electrostatic stress and the interactive force in the electrolyte, we reformulate the fourth-order problem as a mixed system by introducing the gradient of the electrostatic potential, the divergence of its gradient, and the gradient of its divergence as additional independent variables. The second-order problem is studied in its primal form as this approach eliminates the requirement for Lagrange multipliers to enforce the transmission conditions. The resulting formulation is discretised using Raviart--Thomas finite elements for the vector variables and piecewise polynomial finite elements for the scalar variables. Well-posedness of the mixed-primal formulation and optimal a priori error estimates for the field variables as well as for the postprocessed Maxwell tensor  are established under suitable regularity assumptions, in parameter-weighted norms in which the inf-sup constant is independent of the correlation length and of the permittivity ratio. The constants in the error estimates are independent of both parameters as well, up to a loss of at most half an order when the mesh does not resolve the length scale (square root of the relative permittivity); the Sobolev norms of the solution may still depend on the parameters through boundary layers. The postprocessed Maxwell stress converges with optimal order on quasi-uniform meshes, except for the lowest-order scheme in 3D. Numerical experiments validate the theoretical convergence rates and demonstrate the properties of the proposed method for problems in 2D and 3D.
\end{abstract}



\section{Introduction}

Charged biomolecules interacting with their environment are subject to electrostatic forces that are fundamental to biological processes such as the function of proteins, biomolecular self-assembly, ion transport through channels and membrane interactions \cite{Tresset2007, RocchiaAlexovHonig2001}. Understanding ion behaviour in electrolytes is critical to describing cellular functions and designing biomedical systems \cite{Draper2005, Lipfert2014}. This is important because nucleic acids, proteins and biological membranes are all inherently charged.
Traditionally, ion distributions and transport phenomena have been described by the \gls{pb} and the \gls{pnp} theories through a mean-field approximation \cite{Bazant2007}. Although these models have been successful in explaining a wide variety of electrokinetic and transport phenomena \cite{SchochHanRenaud2008}, they neglect steric and electrostatic correlation effects. Consequently, these models do not capture charge inversion by multivalent ions and its sensitive dependence on surface charge density and dielectric parameters \cite{Besteman2005, He2022FEMPF, HeSun2021}.

These considerations have motivated the development of alternative continuum formulations \cite{Vlachy1999}. Among them, the \gls{bsk} theory \cite{BSK2011, BSK2012} incorporates electrostatic correlations via a higher-order diffusion-reaction equation. This model features a Landau--Ginzburg-like free-energy functional, whose extremisation yields a fourth-order modified \gls{pf} equation. Further theoretical, experimental, and numerical investigations (including simulations of three-dimensional biological ion channels \cite{Liu2013}) have corroborated the \gls{bsk} framework in accurately describing dense ionic systems \cite{BSK2012, HeSun2021, JungelMassimini2024}. 
The \gls{bsk} theory can be cast as an interface problem that couples a fourth-order \gls{pf} equation in the solvent region with a second-order Poisson equation in the solute domain.  We recall that \gls{fem} (fitted, unfitted, immersed, etc) for second-order interface problems are well-established  \cite{Babuska1970,Cai2011,Ewing1999,He2010,Chen1998, Li1998, GuoLin2019, LinSheenZhang2019,Hansbo2002, BCHLM2015, LiuZhang2020, HuangWuXiao, BHLZ2025}.  For fourth-order interface problems, such as biharmonic interface problems, constructing globally $C^1$-conforming elements is particularly difficult, motivating mixed and immersed methodologies to relax continuity constraints. See, e.g., \cite{lin2011immersed,zhang2014fourth,ma2025nitsche,chen2025immersedc0interiorpenalty} for unfitted,   
\gls{ipam} schemes (based on fourth-order accurate interface reconstruction), Nitsche nonconforming virtual element methods, and immersed $C^0$ interior penalty methods.  

Mixed \gls{fem} approaches for the modified \gls{pf} equation have received recent attention. While the work in \cite{He2022FEMPF} focuses on analysing and numerically implementing these methods for a single-domain problem, subsequent studies have extended this framework to the more complex nonlinear fourth-/second-order \gls{pf} interface problem \cite{Liu2025InterfacePF, Liu2025MixedPFInterface}. For instance,   \cite{Liu2025InterfacePF}  
recasts the original fourth-order \gls{pf} problem into two coupled second-order Poisson-like equations. That work then presents an interface-fitted \gls{fem} and derives optimal error estimates. In contrast, \cite{Liu2025MixedPFInterface} casts the fourth-order modified \gls{pf} equation into a mixed formulation, choosing the electrostatic potential, its gradient, its Laplacian, and the gradient of its Laplacian as the primary variables in the solvent. For the solute region, governed by a second-order Poisson equation, a standard $H^1$-conforming approximation is applied to the electrostatic potential.   This approach requires artificial stabilisation terms to accommodate standard Lagrange finite element approximations. A unifying objective across these studies is the high-fidelity computation of electrostatic stress to rigorously investigate the charge reversal phenomenon \cite{He2022FEMPF, Liu2025InterfacePF, Liu2025MixedPFInterface}.

In the present work, we restrict our attention to the linearised \gls{bsk} model, as it sufficiently captures the dominant phenomenological interfacial mechanisms, 
including charge reversal, 
within relevant parametric regimes. This model retains the mechanism behind overscreening in the \gls{bsk} theory: the characteristic roots $\lambda_\pm=(1\pm\sqrt{1-4\delta_c^2})/(2\delta_c^2)$ of the operator $\delta_c^2\Delta^2-\Delta+1$, with $\delta_c$ the dimensionless correlation length, are complex for $\delta_c>1/2$, so that the potential and the ionic charge density decay in an oscillatory manner and the first layer of counter-ions can overcompensate the surface charge (see   \Cref{sec:numerical}). However, the linearisation is quantitatively accurate only when the potential is small compared with the thermal voltage $k_{B}T/e$ divided by the ionic valencies and, being symmetric with respect to the sign of the potential, the linearised model does not distinguish multivalent counter-ions from co-ions. For the surface charges of interest the scaled surface potentials are of the order of a few thermal voltages, so we regard the linearised model as a qualitative description of correlation effects, and we compare it with the nonlinear model  in the last test of  \Cref{sec:numerical}. The linearised problem is also the building block for the nonlinear one, since each Newton step leads to a problem with the structure of \eqref{IP}, in which the constant reaction coefficient is replaced by a positive function (minus the derivative of the ionic charge density with respect to the potential). An important objective is the precise evaluation of the Maxwell electrostatic stress tensor and the resulting interaction forces. Because this tensor depends explicitly on the electric field, highly accurate field approximations are essential. To circumvent the order reduction inherent to the post-processing of standard \gls{fem} \cite{He2022FEMPF}, we advocate for a mixed \gls{fem} framework for the fourth-order sub-problem.

While the coupled solvent-solute configuration, transmission conditions, and four-variable split for the mixed formulation are analogous to  \cite{Liu2025MixedPFInterface} (noting additionally the recent three-field formulations in \cite{Farrell2022, KLR2026}), our novel continuous mixed formulation, discrete formulation, and numerical realisation for the coupled system avoids the need for the artificial stabilisation terms. Furthermore, we resolve the fully coupled interface problem via a monolithic approach, which provides a natural alternative to the Dirichlet--Neumann alternating iterations employed in previous studies. 

The main contributions of this work are the following. First, we show that the mixed--primal formulation is well-posed in parameter-weighted norms, with an inf-sup constant that is independent of the dimensionless correlation length $\delta_c$ and of the permittivity ratio $\varepsilon_r$, both at the continuous and at the discrete level. {We also} 
derive quasi-optimal a priori error estimates whose constants do not depend on $\delta_c$, as well as error estimates whose constants are independent of both $\delta_c$ and $\varepsilon_r$; the latter are of optimal order as soon as the mesh size is below $\varepsilon_r^{1/2}$, and lose at most half an order otherwise. This is relevant because $\delta_c$ varies over orders of magnitude with the ionic concentration, and $\varepsilon_r=\varepsilon_2/\varepsilon_1$ is typically small. Second, the discretisation uses $\vdiv$-conforming Raviart--Thomas elements for the fluxes, which renders it stable without any stabilisation or Lagrange multipliers, and it yields an $L^2$-optimal approximation of the Maxwell stress tensor. We establish the error estimate of stress tensor under weaker regularity assumptions in comparison to  \cite{Liu2025MixedPFInterface}. Third, we provide a set of numerical examples corroborating optimal convergence of primal, dual, and postprocessed quantities, and we also include tests of the linear and nonlinear models.

The paper is organised as follows. The model   and the variational formulation are presented in  \Cref{sec:model}.  Next, \Cref{sec:wellposedness} gives a detailed well-posedness analysis of the continuous formulation.  \Cref{sec:discrete} is devoted to the finite element formulation and well-posedness. The error analysis of the presented \gls{fem}  is carried out in \Cref{sec:error}. \Cref{sec:numerical} provides the results of numerical experiments for validating the proposed method. Finally, we draw some conclusions and discuss possible extensions in \Cref{sec:conclusion}.

Standard notations on Sobolev spaces  are followed.
For a domain $D$ and $s\ge0$, we denote by $\|\bullet\|_{s,D}$ and $|\bullet|_{s,D}$ the norm and seminorm of $H^s(D)$, with $H^0(D)=L^2(D)$, so that $\|\bullet\|_{0,D}$ is the $L^2(D)$-norm; the same notation is used for vector- and tensor-valued functions, whose spaces are written in bold (e.g., $\bH^s(D)$, $\bL^2(D)$). For $p\in[1,\infty]$, $\|\bullet\|_{0,p,D}$ is the norm of $L^p(D)$ (also for vector- and tensor-valued functions, with the Euclidean and Frobenius norms, respectively). Integrals over $D$ are written without the measure, e.g., $\int_D u\,v$ for the $L^2(D)$ inner product of $u$ and $v$. The notation $A\lesssim B$ abbreviates $A\leq CB$ for some positive generic constant $C$, which depends only on the domain and the shape regularity of the triangulation.

\section{Model problem and weak form}\label{sec:model} 
\begin{figure}[t!]
  \centering
  \subfloat[Immersed case.\label{fig:domainA}]{\begin{minipage}[b]{0.4\textwidth}
    \centering
    \includegraphics{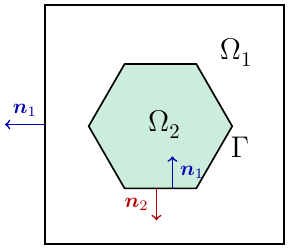}

  \end{minipage}}
  \quad 
  \subfloat[Back-to-back case.\label{fig:domainB}]{\begin{minipage}[b]{0.55\textwidth}
    \centering
    \includegraphics{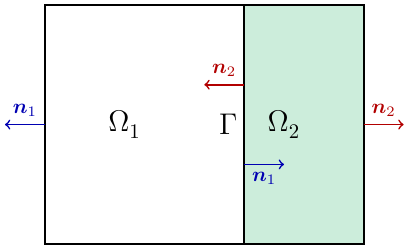}

  \end{minipage}}
  \caption{Graphical representation of sub-domains $\Omega_{1}^{},\Omega_{2}^{}$ with interface $\Gamma$. Other options are possible for analysis and computation, as mentioned in  \Cref{rem:domains}.}
  \label{fig:domain}
\end{figure}

Let $\Omega\subset\mathbb{R}^{d},\ d\in \{2,3\}$, be an open bounded domain with Lipschitz boundary $\partial\Omega$ and $\bm{n}$ be the unit outward normal on $\partial\Omega$. Subdivide $\Omega$ into two non-overlapping Lipschitz domains $\Omega_{1}^{}$ and $\Omega_{2}^{}$ with an interface $\Gamma$ between them. That is, $\overline\Omega=\overline\Omega_{1}\cup\overline\Omega_{2}$, $\Gamma=\partial\Omega_{1}^{}\cap\partial\Omega_{2}^{}$ and $\Omega_{1}^{}\cap\Omega_{2}^{}=\emptyset$. By $\bm{n}_i$ we denote the outward unit normal on $\partial\Omega_i$; on $\Gamma$ it points from $\Omega_i$ towards its neighbouring sub-domain, so that $\bm n_2=-\bm n_1$ there. The two main cases of domain decomposition consist of  $\Gamma=\partial\Omega_{2}^{}$ in the \textit{immersed case} (see \Cref{fig:domainA}) and $\Gamma\subsetneq\partial\Omega_2$ with $\overline\Gamma\cap\partial\Omega\neq\emptyset$ in the \textit{back-to-back case} (see \Cref{fig:domainB}).

For sake of the application herein, $\Omega_{1}$ is considered the solvent or electrolyte region and $\Omega_{2}$ the solute or molecular region. We consider the following linear fourth-order/second-order interface problem comprising a linear version of the fourth-order modified \gls{pf} equation proposed by \gls{bsk} theory \cite{BSK2011} in the solvent $\Omega_{1}$ and a second-order Poisson-type equation in the solute region $\Omega_{2}$:
\begin{align}\label{IP}
\delta_{c}^{2}\Delta^{2}\varphi_{1}^{}-\Delta\varphi_{1}^{} + \varphi_{1}^{} &= f_{1}^{} && \text{ in }\Omega_{1}^{},\nonumber\\
    -\varepsilon_{r}\Delta\varphi_{2}^{}+\varphi_{2}^{}&=f_{2}^{} && \text{ in }\Omega_{2}^{},\nonumber\\
    \varphi_{1}^{}-\varphi_{2}^{}&=j && \text{ on }\Gamma,\nonumber\\
    \nabla\varphi_{1}^{}\cdot \bm{n}_1+\varepsilon_{r}\nabla\varphi_2\cdot\bm{n}_2&=g && \text{ on }\Gamma,\\
    \varphi_{1}^{}&=g_{1}^{} && \text{ on }\partial\Omega_{1}^{}\setminus\Gamma ,\nonumber\\
    \varphi_{2}^{}&=g_{2}^{} && \text{ on }\partial\Omega_{2}^{}\setminus\Gamma,\nonumber\\
    -\Delta\varphi_{1}^{}&=g_{3} && \text{ on }\partial\Omega_{1}^{},\nonumber
\end{align}
where $f_{i}\in L^2(\Omega_{i})$ are given source/sink terms, and $j$, $g$, $g_{1}^{}$, $g_{2}^{}$, $g_3$ are  boundary and interface data (their precise regularity and compatibility are specified in \eqref{eq:compat} below). The constant $\delta_{c}^{}$ in (\ref{IP})$_{1}$ is the dimensionless correlation length  given by $\delta_{c}=\frac{l_{c}}{\lambda_{D}}$, where $l_{c}$ is the electrostatic correlation length and $\lambda_{D}$ is the Debye length \cite{BSK2011}. The parameter $\varepsilon_{r}$ in (\ref{IP})$_{2}$ and (\ref{IP})$_{4}$ is the relative permittivity. Equations (\ref{IP})$_{3}-$(\ref{IP})$_{4}$ are the  transmission conditions on $\Gamma$ and (\ref{IP})$_{5}-$(\ref{IP})$_{7}$ are the boundary conditions.
The relations (\ref{IP})$_3$--(\ref{IP})$_4$ and   (\ref{IP})$_7$   which prescribes $-\Delta\varphi_1$ (the scaled charge density in the absence of correlations, $\delta_c=0$) on the solvent side of the interface, are  used in \cite{Liu2025InterfacePF,Liu2025MixedPFInterface}. The zeroth-order term in (\ref{IP})$_2$ stems from the linearisation of a nonlinear response $f_2(\varphi_2)$ in the solute, as in \cite{Liu2025MixedPFInterface}, and it can also be regarded as a regularisation term, but for a dielectric solute that contains no mobile ions, the physically appropriate equation is the Poisson equation $-\varepsilon_r\Delta\varphi_2=f_2$. While the analysis uses the zeroth-order term as in \cite{Liu2025MixedPFInterface}, some of the simulations in \Cref{sec:numerical} are performed with the Poisson case.

\begin{remark}\label{rem:domains}
Although our analysis assumes a general back-to-back domain configuration, the results naturally extend to other standard geometries. First, if $\Omega_{2}$ is completely immersed within $\Omega_{1}$ (see \Cref{fig:domainA}), we simply have $\Gamma = \partial\Omega_{2}$. In this case, the boundary condition (\ref{IP})$_{6}$ is void and $W_0=W$, and the arguments hold without further modification. 

Second, for clarity of presentation, we treat $\Omega_{2}$ as a single domain. However, because the transmission conditions are strictly local to the interface, all constructions and results remain valid if $\Omega_{2}$ comprises a finite union of pairwise disjoint Lipschitz subdomains, $\Omega_{2}^{} = \bigcup_{i=1}^{m} \Omega_{2}^{i}$, with $\Gamma = \bigcup_{i=1}^{m} \partial\Omega_{2}^{i}$ (as illustrated for $m=2$ in \Cref{fig:domain2} and \Cref{fig:domain3}), provided that no part of $\Omega_1$ is enclosed by $\Gamma$, that is, every connected component of $\Omega_1$ touches $\partial\Omega$ on a set of positive measure (this is needed in \Cref{infsupb}). In such cases, the analysis proceeds component-wise with  minor  adjustments.
\end{remark}

In order to compute the interactive force, it is important to precisely approximate the electric and displacement fields. We reformulate the model problem (\ref{IP}) using the mixed variables $\bp\equiv \nabla \varphi_1$ (potential gradient, that is, minus the electric field), $u\equiv -\vdiv\bp$ (which reduces to the scaled charge density for $\delta_c=0$) and $\bzeta\equiv \delta_c^2\nabla u$ (correlation contribution to the displacement field). 
Then, (\ref{IP}) takes the strong mixed-primal form
\begin{align}\label{IPM}
		\bp&=\nabla \varphi_1 && \text{ in }\Omega_{1}^{},\nonumber\\
        -\vdiv\bp &= u&& \text{ in }\Omega_{1}^{},\nonumber\\
        \bzeta&=\delta_c^2\nabla u && \text{ in }\Omega_{1}^{},\nonumber\\
        -\vdiv\bzeta + u + \varphi_{1}^{} &= f_{1}^{} && \text{ in }\Omega_{1}^{},\nonumber\\
        -\varepsilon_{r}\Delta\varphi_{2}^{}+\varphi_{2}^{}&=f_{2}^{} && \text{ in }\Omega_{2}^{},\\
        \varphi_{1}^{}-\varphi_{2}^{}&=j && \text{ on }\Gamma,\nonumber\\
        \bp\cdot\bm{n}_{1}+\varepsilon_{r}\nabla\varphi_2\cdot\bm{n}_{2}&=g && \text{ on }\Gamma,\nonumber\\
        \varphi_{1}^{}&=g_{1}^{} && \text{ on }\partial\Omega_{1}^{}\setminus\Gamma ,\nonumber\\
        \varphi_{2}^{}&=g_{2}^{} && \text{ on }\partial\Omega_{2}^{}\setminus\Gamma,\nonumber\\
        u&=g_{3} && \text{ on }\partial\Omega_{1}.\nonumber
\end{align}   
Owing to the subsequent imposition of boundary and interfacial conditions, we consider the spaces \begin{gather*}
\bQ = \bH(\mathrm{div};\Omega_{1}^{}) := \{\bm{v} \in \bL^2(\Omega_{1}^{}) : \vdiv \bm{v} \in {L}^2(\Omega_{1}^{})\},\qquad V:=L^2(\Omega_{1}^{}),
\\
W := H^1(\Omega_{2}^{}), \qquad W_{0} := \{v \in H^1(\Omega_{2}^{}) : v = 0 \text{ on }\partial\Omega_{2}^{}\setminus\Gamma\}.
\end{gather*}
\textit{Data and traces.} We assume that the Dirichlet data are compatible in the following sense: there exist $\tilde g_1\in H^1(\Omega_1)$ and $\tilde g_2\in H^1(\Omega_2)$ such that
\begin{equation}\label{eq:compat}
\tilde g_1=g_1 \ \text{ on }\partial\Omega_1\setminus\Gamma,\qquad \tilde g_2=g_2\ \text{ on }\partial\Omega_2\setminus\Gamma,\qquad \tilde g_1-\tilde g_2=j\ \text{ on }\Gamma,
\end{equation}
that is, $(g_1,g_2,j)$ are the traces of a piecewise $H^1$ function with the prescribed jump across $\Gamma$. In the immersed case ($\Gamma\cap\partial\Omega=\emptyset$ and $\partial\Omega_2\setminus\Gamma=\emptyset$), \eqref{eq:compat} only requires $g_1\in H^{1/2}(\partial\Omega)$ and $j\in H^{1/2}(\Gamma)$. In the back-to-back case it is a genuine compatibility condition on $\partial\Gamma=\overline{\Gamma}\cap\partial\Omega$; for instance, for piecewise smooth  continuous data it amounts to $g_1=g_2+j$ on $\partial\Gamma$. We denote by $g_D:=\tilde g_1|_{\partial\Omega_1}\in H^{1/2}(\partial\Omega_1)$ the Dirichlet datum of $\varphi_1$ on the whole of $\partial\Omega_1$, which coincides with $g_1$ on $\partial\Omega_1\setminus\Gamma$ and with $\tilde g_2+j$ on $\Gamma$.
For $\psi_2\in W_0$, the trace $\psi_2|_\Gamma$ belongs to $H^{1/2}_{00}(\Gamma)$, so that its extension by zero to $\partial\Omega_1$, denoted by $E_\Gamma\psi_2$, belongs to $H^{1/2}(\partial\Omega_1)$ (in the immersed case, $H^{1/2}_{00}(\Gamma)=H^{1/2}(\Gamma)$). Here we use that $H^{1/2}_{00}(\Gamma)$ is intrinsic to $\Gamma$ (its norm is equivalent to $\big(\|\mu\|_{1/2,\Gamma}^2+\int_\Gamma\mu^2\,\mathrm{dist}(\bx,\partial\Gamma)^{-1}\ds\big)^{1/2}$), so that, {since $\Gamma$ is a Lipschitz $(d-1)$-manifold with Lipschitz
boundary $\partial\Gamma$,} 
the extension by zero is bounded both into $H^{1/2}(\partial\Omega_2)$ and into $H^{1/2}(\partial\Omega_1)$.  Accordingly, we assume $g\in H^{1/2}_{00}(\Gamma)'$ (which is $H^{-1/2}(\Gamma)$ in the immersed case) and $g_3\in H^{1/2}(\partial\Omega_1)$, and we write $\langle\bullet,\bullet\rangle_{\partial\Omega_1}$ for the duality pairing between $H^{-1/2}(\partial\Omega_1)$ and $H^{1/2}(\partial\Omega_1)$, and $\langle\bullet,\bullet\rangle_{\Gamma}$ for the one between $H^{1/2}_{00}(\Gamma)'$ and $H^{1/2}_{00}(\Gamma)$.

Multiplying by suitable test functions, integrating by parts, and applying the boundary and transmission conditions (on $\Gamma$, $\varphi_1=\varphi_2+j=(\varphi_2-\tilde g_2)+g_D$, while $\varphi_1=g_D$ on $\partial\Omega_1\setminus\Gamma$; that is, $\varphi_1|_{\partial\Omega_1}=g_D+E_\Gamma(\varphi_2-\tilde g_2)$), we obtain a mixed variational formulation for (\ref{IPM}): find  $(\bp,\bzeta,\varphi_{1}^{},u,\varphi_{2}^{})\in \bQ\times \bQ\times V\times V\times W$ with $\varphi_{2}^{}-\tilde g_{2}^{}\in W_0$ such that
\begin{align}\label{IPw}
    a_{1}^{}(\bp,\bq) + b(\bq,\varphi_{1}^{}) - b_{\Gamma}(\bq,\varphi_2-\tilde g_2)&= l_{1}^{}(\bq) &&\text{ for all } \bq\in \bQ,\nonumber\\
    a_{2}^{}(\bzeta,\bxi) + b(\bxi,u) &= l_{2}^{}(\bxi)&&\text{ for all } \bxi\in \bQ,\nonumber\\
    b(\bp,\psi_{1}^{}) + a(u,\psi_{1}^{}) &= 0&&\text{ for all } \psi_{1}^{}\in V,\\
    b(\bzeta,v) - a(u,v) - a(\varphi_{1}^{},v) &= l_{4}^{}(v)&&\text{ for all } v\in V,\nonumber\\
    a_{3}^{}(\varphi_{2}^{},\psi_{2}^{}) + b_{\Gamma}(\bp,\psi_{2}^{})& = l_{5}^{}(\psi_{2}^{})&&\text{ for all } \psi_{2}^{}\in W_{0},\nonumber
\end{align}
where the bilinear and linear forms are defined by
{\begin{align*}
       & a_{1}^{}(\bp,\bq):=\int_{\Omega_1}\bp\cdot\bq,\quad b(\bq,\varphi_{1}^{}):=\int_{\Omega_1}\varphi_{1}^{}\,\vdiv \bq,\quad b_{\Gamma}(\bq,\varphi_{2}^{}):=\langle \bq\cdot\bm{n}_{1}^{}, E_\Gamma\varphi_{2}^{} \rangle_{\partial\Omega_1}, \\
       & a_{2}^{}(\bzeta,\bxi):=\delta_{c}^{-2}\int_{\Omega_1}\bzeta\cdot\bxi,\quad a(u,v):=\int_{\Omega_1}u\,v,\quad a_{3}^{}(\varphi_{2}^{},\psi_{2}^{}):=\varepsilon_{r}^{} \int_{\Omega_2}\nabla\varphi_{2}^{}\cdot\nabla\psi_{2}^{} + \int_{\Omega_2}\varphi_{2}^{}\,\psi_{2}^{},  \\
       &l_{1}^{}(\bq):=\langle \bq\cdot\bm{n}_{1}^{},g_D \rangle_{\partial\Omega_{1}},\ l_{2}^{}(\bxi):=\langle \bxi\cdot\bm{n}_{1}^{},g_{3}^{} \rangle_{\partial\Omega_{1}^{}},\    l_{4}^{}(v):= - \int_{\Omega_1}f_{1}^{}\,v,\ l_{5}^{}(\psi_{2}^{}):=\int_{\Omega_2}f_{2}^{}\,\psi_{2}^{} + \langle g,\psi_{2}^{} \rangle_{\Gamma}^{}.
\end{align*}}
Since $E_\Gamma\psi_2$ vanishes on $\partial\Omega_1\setminus\Gamma$, the form $b_\Gamma$ only involves the interface: it is the duality pairing between $\bq\cdot\bm n_1|_\Gamma\in H^{1/2}_{00}(\Gamma)'$ and $\psi_2|_\Gamma\in H^{1/2}_{00}(\Gamma)$. In particular, if $\bq\cdot\bm n_1\in L^2(\partial\Omega_1)$, which is the case for the discrete fluxes below, then $b_\Gamma(\bq,\psi_{2})=\langle \bq\cdot\bm{n}_{1}^{}, \psi_{2}^{} \rangle_{\Gamma}$. 
The pairing on $\partial\Omega_1$ is only needed to give a meaning to this integral for a general $\bq\in\bQ$, whose normal trace belongs to $H^{-1/2}(\partial\Omega_1)$. In the same situation, since $g_D=\tilde g_2+j$ on $\Gamma$, the first equation in \eqref{IPw} reads
\[
a_1(\bp,\bq)+b(\bq,\varphi_1)-{\langle \bq\cdot\bm{n}_{1}^{}, \varphi_{2}^{} \rangle_{\Gamma}}={\langle \bq\cdot\bm{n}_{1}^{},g_1 \rangle_{\partial\Omega_{1}\setminus\Gamma} + \langle \bq\cdot\bm{n}_{1}^{},j \rangle_{\Gamma}} 
\]
which does not involve the lifting $\tilde g_2$ and is the form used in the implementation.

\textit{Homogenisation.} Define the product space $\vec{\bY} := \bQ\times \bQ\times V\times V\times W_{0}$ with elements $\vec{\bx}=(\bp,\bzeta,\varphi_1,u,\varphi_2)$ and $\vec{\by}=(\bq,\bxi,\psi_{1}^{},v,\psi_{2}^{})$, and the bilinear form $\mathcal{M}:\vec{\bY}\times\vec{\bY}\to\mathbb{R}$ and the linear form $l:\vec{\bY}\to\mathbb{R}$ given by
{\begin{align*}
   \mathcal{M}(\vec{\bx};\vec{\by}) &:= a_1(\bp,\bq)+b(\bq,\varphi_1)-b_\Gamma(\bq,\varphi_2)+a_2(\bzeta,\bxi)+b(\bxi,u)+ b(\bp,\psi_1)\\
   &\qquad +a(u,\psi_1)+b(\bzeta,v)-a(u,v)-a(\varphi_1,v)+a_3(\varphi_2,\psi_2)+b_\Gamma(\bp,\psi_2) ,\\
    l(\vec{\by}) &:= l_1(\bq)+l_{2}^{}(\bxi)+l_4(v)+l_5(\psi_2)-a_3(\tilde g_2,\psi_2).
\end{align*}}
Writing the solution of \eqref{IPw} as $(\bp,\bzeta,\varphi_1,u,\varphi_2^0+\tilde g_2)$ with $\varphi_2^0:=\varphi_2-\tilde g_2\in W_0$, problem \eqref{IPw} is equivalent to: find $\vec{\bx}=(\bp,\bzeta,\varphi_1,u,\varphi_2^0)\in\vec{\bY}$ such that
\begin{equation}\label{Pw}
   \mathcal{M}(\vec{\bx};\vec{\by})=l(\vec{\by})\quad\text{ for all } \vec{\by}\in \vec{\bY}.
\end{equation}
Trial and test spaces now coincide, which is needed for the inf-sup arguments below.
To simplify notation, in the analysis we write $\varphi_2$ instead of $\varphi_2^0$.

\section{Well-posedness of the weak formulation}\label{sec:wellposedness}
Since neither the $\bQ\times\bQ\times W_0$-block nor the $V\times V$-block of $\mathcal M$ is symmetric, the abstract theory of perturbed saddle-point problems \cite[Chapter 4]{Boffi}, \cite{Braess} does not apply directly. Instead, we establish a global inf-sup condition for $\mathcal M$ by an explicit construction of test functions. A main point is that this is done in parameter-weighted norms, so that the stability constant is independent of $\delta_c$ and $\varepsilon_r$.

\textit{Parameter-weighted norms.} We define
\begin{gather*}
\|\bq\|_{\bQ}^2:=\|\bq\|_{0,\Omega_1}^2+\|\vdiv\bq\|_{0,\Omega_1}^2,\quad
\|\bxi\|_{\bQ,\delta_c}^2:=\delta_c^{-2}\|\bxi\|_{0,\Omega_1}^2+\|\vdiv\bxi\|_{0,\Omega_1}^2, \\
\|\psi_2\|_{W,\varepsilon_r}^2:=\varepsilon_r\|\nabla\psi_2\|_{0,\Omega_2}^2+\|\psi_2\|_{0,\Omega_2}^2,\quad 
\normiii{\vec{\by}}^2:=\|\bq\|_{\bQ}^2+\|\bxi\|_{\bQ,\delta_c}^2+\|\psi_1\|_{0,\Omega_1}^2+\|v\|_{0,\Omega_1}^2+\|\psi_2\|_{W,\varepsilon_r}^2,
\end{gather*}
for $\vec{\by}=(\bq,\bxi,\psi_1,v,\psi_2)\in\vec{\bY}$. These weights are dictated by  \eqref{IPM}:   $(\bzeta,u)$ solves the singularly perturbed reaction--diffusion problem $\delta_c^{-2}\bzeta=\nabla u$, $-\vdiv\bzeta+u=f_1-\varphi_1$, whose natural energy is $\delta_c^{-2}\|\bzeta\|_{0,\Omega_1}^2+\|\vdiv\bzeta\|_{0,\Omega_1}^2+\|u\|_{0,\Omega_1}^2$, the pair $(\bp,\varphi_1)$ solves a Darcy-type problem, and $\|\bullet\|_{W,\varepsilon_r}$ is the energy norm induced by $a_3$.

\medskip \noindent
{With respect to this norm setting, we can show that $b_\Gamma(\bullet,\bullet)$ is bounded by utilising the normal trace theorem in $\bQ$ \cite[Lemma 2.1.1]{Boffi}, the trace theorem in $H^1(\Omega_2)$ \cite[Theorem 3.10]{Ern1}, and the boundedness of $E_\Gamma$. To this end, we consider $H^{1/2}(\partial\Omega_1)$ with the quotient norm $\|\mu\|_{1/2,\partial\Omega_1}:=\inf\{\|G\|_{1,\Omega_1}:\ G\in H^1(\Omega_1),\ G|_{\partial\Omega_1}=\mu\}$ and $H^{1/2}_{00}(\Gamma)$ with the norm $\|\mu\|_{H^{1/2}_{00}(\Gamma)}:=\|E_\Gamma\mu\|_{1/2,\partial\Omega_1}$. Then, for $\bq\in\bQ$ and any extension $G$ of $\mu$, Green's formula gives $\langle\bq\cdot\bm n_1,\mu\rangle_{\partial\Omega_1}=\int_{\Omega_1}\bq\cdot\nabla G+\int_{\Omega_1}\vdiv\bq\,G\le\|\bq\|_{\bQ}\|G\|_{1,\Omega_1}$, and taking the infimum over $G$,
\begin{equation}\label{eq:normaltrace}
|\langle\bq\cdot\bm n_1,\mu\rangle_{\partial\Omega_1}|\le\|\bq\|_{\bQ}\|\mu\|_{1/2,\partial\Omega_1}\qquad\text{for all }\bq\in\bQ,\ \mu\in H^{1/2}(\partial\Omega_1).
\end{equation}
This yields\begin{equation}\label{eq:tr}
    |b_\Gamma(\bq,\psi_2)| \le C_{\rm tr}\|\bq\|_{\bQ}\|\psi_2\|_{1,\Omega_2},
\end{equation}
where $C_{\rm tr}$ in \eqref{eq:tr} is the norm of the map $W_0\ni\psi_2\mapsto E_\Gamma\psi_2\in H^{1/2}(\partial\Omega_1)$ (with $W_0$ normed by $\|\bullet\|_{1,\Omega_2}$), that is, $\|\psi_2\|_{H^{1/2}_{00}(\Gamma)}\le C_{\rm tr}\|\psi_2\|_{1,\Omega_2}$. Finally, it follows by utilising $\|\psi_2\|_{1,\Omega_2}^2=\varepsilon_r^{-1}\,\varepsilon_r\|\nabla\psi_2\|_{0,\Omega_2}^2+\|\psi_2\|_{0,\Omega_2}^2\le\max\{1,\varepsilon_r^{-1}\}\|\psi_2\|_{W,\varepsilon_r}^2$ in \eqref{eq:tr} that
\begin{equation}\label{eq:trace}
|b_\Gamma(\bq,\psi_2)| \le C_\Gamma \|\bq\|_{\bQ}\|\psi_2\|_{W,\varepsilon_r}
\qquad\text{for all }\bq\in\bQ,\ \psi_2\in W_0,
\end{equation}
where $C_\Gamma:=C_{\rm tr}\max\{1,\varepsilon_r^{-1/2}\}$.}
\begin{lemma}[Continuity]\label{continuity}
For all $\vec{\bx},\vec{\by}\in\vec{\bY}$ there holds $|\mathcal{M}(\vec{\bx};\vec{\by})|\le C_{\mathcal M}\normiii{\vec{\bx}}\,\normiii{\vec{\by}}$ with $C_{\mathcal M}:=3+2C_\Gamma$. Moreover, under the assumptions of \Cref{sec:model},
\begin{gather*}
|l(\vec{\by})|\le C_l\normiii{\vec{\by}}, \quad \text{with} \\
 C_l:= \|g_D\|_{1/2,\partial\Omega_1}+\max\{1,\delta_c\}\|g_3\|_{1/2,\partial\Omega_1}+\|f_1\|_{0,\Omega_1}+\|f_2\|_{0,\Omega_2}   +C_\Gamma\|g\|_{H^{1/2}_{00}(\Gamma)'}+\|\tilde g_2\|_{W,\varepsilon_r}.
\end{gather*}
\end{lemma}
\begin{proof}
All the terms of $\mathcal M$ except the two interface terms are $L^2$ inner products of one component of $\vec{\bx}$ and one component of $\vec{\by}$ (or of their divergences), weighted consistently with $\normiii{\bullet}$; for instance $a_2(\bzeta,\bxi)=\int_{\Omega_1}(\delta_c^{-1}\bzeta)\cdot(\delta_c^{-1}\bxi)$. More precisely, let $\vec{\bx}=(\bp,\bzeta,\varphi_1,u,\varphi_2)$, $\vec{\by}=(\bq,\bxi,\psi_1,v,\psi_2)\in\vec{\bY}$, and collect the seven quantities entering $\normiii{\bullet}$ in the vectors
\[
X=(X_1,\ldots,X_7):=\big(\|\bp\|_{0,\Omega_1},\|\vdiv\bp\|_{0,\Omega_1},\delta_c^{-1}\|\bzeta\|_{0,\Omega_1},\|\vdiv\bzeta\|_{0,\Omega_1},\|\varphi_1\|_{0,\Omega_1},\|u\|_{0,\Omega_1},\|\varphi_2\|_{W,\varepsilon_r}\big),
\]
and $Y{=(Y_1,\ldots,Y_7)}:=\big(\|\bq\|_{0,\Omega_1},\|\vdiv\bq\|_{0,\Omega_1},\delta_c^{-1}\|\bxi\|_{0,\Omega_1},\|\vdiv\bxi\|_{0,\Omega_1},\|\psi_1\|_{0,\Omega_1},\|v\|_{0,\Omega_1},\|\psi_2\|_{W,\varepsilon_r}\big)$.
The Cauchy--Schwarz inequality (in $L^2$, and for $a_3$ in the inner product that induces $\|\bullet\|_{W,\varepsilon_r}$) gives
{\begin{align*}
&|a_1(\bp,\bq)|\le X_1Y_1,\quad |b(\bq,\varphi_1)|\le X_5Y_2,\quad |a_2(\bzeta,\bxi)|\le X_3Y_3,\quad |b(\bxi,u)|\le X_6Y_4,\quad |b(\bp,\psi_1)|\le X_2Y_5,\\
&|a(u,\psi_1)|\le X_6Y_5,\quad |b(\bzeta,v)|\le X_4Y_6,\quad |a(u,v)|\le X_6Y_6,\quad |a(\varphi_1,v)|\le X_5Y_6,\quad |a_3(\varphi_2,\psi_2)|\le X_7Y_7.
\end{align*}}
Denote by $S$ the set of the ten index pairs $(i,j)$ appearing on the right-hand sides. The index $i=6$ (that is, $\|u\|_{0,\Omega_1}$) appears in three pairs, $i=5$ in two, and every other $i$ in one; likewise $j=6$ (that is, $\|v\|_{0,\Omega_1}$) appears in three pairs, $j=5$ in two, and every other $j$ in one. Assume first that 
{$\normiii{\vec{\bx}}=\normiii{\vec{\by}}=1$.} Then, by $X_iY_j\le\frac12(X_i^2+Y_j^2)$,
\[
\sum_{(i,j)\in S}X_iY_j\le\frac12\sum_{i=1}^7 r_iX_i^2+\frac12\sum_{j=1}^7c_jY_j^2\le \frac32 (\normiii{\vec{\bx}}^2+
{\normiii{\vec{\by}}}^2)=3,
\]
where $r_i\le3$ and $c_j\le3$ are the numbers of pairs in $S$ containing $i$ and $j$, respectively. For general $\vec{\bx},\vec{\by}\ne\bm0$, we apply this to $X/\normiii{\vec{\bx}} $ and $Y/\normiii{\vec{\by}}$ and obtain, by homogeneity, that the ten non-interface terms of $\mathcal M(\vec{\bx};\vec{\by})$ are bounded in absolute value by $3\normiii{\vec{\bx}}\normiii{\vec{\by}}$ (the bound is trivial if $\vec{\bx}=\bm0$ or $\vec{\by}=\bm0$). For the two interface terms, \eqref{eq:trace} gives
\[
|b_\Gamma(\bq,\varphi_2)|\le C_\Gamma\|\bq\|_{\bQ}\|\varphi_2\|_{W,\varepsilon_r}\le C_\Gamma\normiii{\vec{\by}}\normiii{\vec{\bx}},\qquad
|b_\Gamma(\bp,\psi_2)|\le C_\Gamma\|\bp\|_{\bQ}\|\psi_2\|_{W,\varepsilon_r}\le C_\Gamma\normiii{\vec{\bx}}\normiii{\vec{\by}}.
\]
Adding up, $|\mathcal M(\vec{\bx};\vec{\by})|\le(3+2C_\Gamma)\normiii{\vec{\bx}}\normiii{\vec{\by}}$. 


For $l$, we bound each of its terms. By \eqref{eq:normaltrace}, we have $|l_1(\bq)|=|\langle\bq\cdot\bm n_1,g_D\rangle_{\partial\Omega_1}|\le\|\bq\|_{\bQ}\|g_D\|_{1/2,\partial\Omega_1}$. Again \eqref{eq:normaltrace} and the definition of $\|\cdot\|_{\bQ,\delta_c}$ lead to
\[
|l_{2}^{}(\bxi)|\le(\|\bxi\|_{0,\Omega_1}^2+\|\vdiv\bxi\|_{0,\Omega_1}^2)^{\frac12}\|g_3\|_{1/2,\partial\Omega_1}\le\max\{1,\delta_c\}\|\bxi\|_{\bQ,\delta_c}\|g_3\|_{1/2,\partial\Omega_1}.
\]
Moreover, $|l_4(v)|\le\|f_1\|_{0,\Omega_1}\|v\|_{0,\Omega_1}$, $|\int_{\Omega_2}f_2\,\psi_2|\le\|f_2\|_{0,\Omega_2}\|\psi_2\|_{0,\Omega_2}\le\|f_2\|_{0,\Omega_2}\|\psi_2\|_{W,\varepsilon_r}$, $|a_3(\tilde g_2,\psi_2)|\le\|\tilde g_2\|_{W,\varepsilon_r}\|\psi_2\|_{W,\varepsilon_r}$, and, by duality and the trace estimates following \eqref{eq:trace},
\[
|\langle g,\psi_2\rangle_\Gamma|\le\|g\|_{H^{1/2}_{00}(\Gamma)'}\|\psi_2\|_{H^{1/2}_{00}(\Gamma)}\le C_{\rm tr}\|g\|_{H^{1/2}_{00}(\Gamma)'}\|\psi_2\|_{1,\Omega_2}\le C_\Gamma\|g\|_{H^{1/2}_{00}(\Gamma)'}\|\psi_2\|_{W,\varepsilon_r}.
\]
Since each of $\|\bq\|_{\bQ}$, $\|\bxi\|_{\bQ,\delta_c}$, $\|v\|_{0,\Omega_1}$, and $\|\psi_2\|_{W,\varepsilon_r}$ is bounded by $\normiii{\vec{\by}}$, summing these six bounds {yields} $|l(\vec{\by})|\le C_l\normiii{\vec{\by}}$.
\end{proof}

The next result is the surjectivity of the divergence on fluxes with vanishing normal trace on $\Gamma$. Let $\bQ_\Gamma:=\{\bq\in\bQ:\ b_\Gamma(\bq,\psi_2)=0\ \text{for all }\psi_2\in W_0\}$.
\begin{lemma}[{Surjectivity of the divergence}]\label{infsupb}
There exists {$0<\beta\leq 1$}, depending only on $\Omega_1$ and $\Gamma$, such that for every $v\in V$ there is $\widehat{\bq}\in\bQ_\Gamma$ with $\vdiv\widehat{\bq}=v$ and $\|\widehat{\bq}\|_{\bQ}\le\beta^{-1}\|v\|_{0,\Omega_1}$. 
\end{lemma}
\begin{proof}
{This is standard (see, e.g., \cite[Section 51.1]{Ern2}): take $\widehat{\bq}:=-\nabla w$, where $w$ solves $-\Delta w=v$ in $\Omega_1$, $w=0$ on $\partial\Omega_1\setminus\Gamma$, $\nabla w\cdot\bm n_1=0$ on $\Gamma$, and utilise the Poincar\'e inequality $\|w\|_{0,\Omega_1}\le C_P\|\nabla w\|_{0,\Omega_1}$, which gives $\beta=(1+C_P^2)^{-1/2}\le1$. Since $\bQ$ carries the unweighted norm, $\beta$ is independent of $\delta_c$ and $\varepsilon_r$.}
\end{proof}

\begin{theorem}[Parameter-robust stability]\label{stability}
Let $\gamma:=\beta^2/(16+8\beta+2\beta^2)$. Then
\begin{gather}\label{iinfsup}
\sup_{\vec{\by}\in\vec{\bY}\setminus\{\bm0\}}\frac{\mathcal{M}(\vec{\bx};\vec{\by})}{\normiii{\vec{\by}}}\ge\gamma\normiii{\vec{\bx}}\quad\text{for all }\vec{\bx}\in\vec{\bY},\qquad 
\sup_{\vec{\bx}\in\vec{\bY}\setminus\{\bm0\}}\frac{\mathcal{M}(\vec{\bx};\vec{\by})}{\normiii{\vec{\bx}}}\ge\gamma\normiii{\vec{\by}}\quad\text{for all }\vec{\by}\in\vec{\bY}.
\end{gather}
In particular, $\gamma$ is independent of $\delta_c$ and $\varepsilon_r$.
\end{theorem}
\begin{proof}
Let $\vec{\bx}=(\bp,\bzeta,\varphi_1,u,\varphi_2)\in\vec{\bY}$, and let $\widehat{\bq}\in\bQ_\Gamma$ be given by \Cref{infsupb} with $\vdiv\widehat{\bq}=\varphi_1$. We test with three functions.

\emph{(i)} For $\vec{\by}_1:=(\bp,\bzeta,-\varphi_1,-u,\varphi_2)$, the terms involving $b$ and $b_\Gamma$ cancel pairwise 
and also $a(u,-\varphi_1)-a(\varphi_1,-u)=0$, while $-a(u,-u)=\|u\|_{0,\Omega_1}^2$. Therefore,
\[
\mathcal{M}(\vec{\bx};\vec{\by}_1)=\|\bp\|_{0,\Omega_1}^2+\delta_c^{-2}\|\bzeta\|_{0,\Omega_1}^2+\|\varphi_2\|_{W,\varepsilon_r}^2+\|u\|_{0,\Omega_1}^2 .
\]
\emph{(ii)} For $\vec{\by}_2:=(\widehat{\bq},\bm 0,0,0,0)$, only the terms $a_1(\bp,\widehat{\bq})+b(\widehat{\bq},\varphi_1)-b_\Gamma(\widehat{\bq},\varphi_2)$ remain. Since $b_\Gamma(\widehat{\bq},\varphi_2)=0$, $\vdiv\widehat{\bq}=\varphi_1$. 
\Cref{infsupb}, Cauchy--Schwarz and Young's inequalities yield $|\int_{\Omega_1}\bp\cdot\widehat{\bq}|\le  \beta^{-1} \|\bp\|_{0,\Omega_1}\|\varphi_1\|_{0,\Omega_1}\le\frac12\beta^{-2}\|\bp\|_{0,\Omega_1}^2+\frac12\|\varphi_1\|_{0,\Omega_1}^2$.
Hence 
\[
\mathcal{M}(\vec{\bx};\vec{\by}_2)=\int_{\Omega_1}\bp\cdot\widehat{\bq}+\|\varphi_1\|_{0,\Omega_1}^2\ge\tfrac12\|\varphi_1\|_{0,\Omega_1}^2-\tfrac12\beta^{-2}\|\bp\|_{0,\Omega_1}^2 .
\]
\emph{(iii)} For $\vec{\by}_3:=(\bm 0,\bm 0,\vdiv\bp,\vdiv\bzeta,0)$, which belongs to $\vec{\bY}$ since $\vdiv\bp,\vdiv\bzeta\in L^2(\Omega_1)$, the non-zero terms of $\mathcal{M}(\vec{\bx};\vec{\by}_3)$ are $b(\bp,\vdiv\bp)+a(u,\vdiv\bp)+b(\bzeta,\vdiv\bzeta)-a(u,\vdiv\bzeta)-a(\varphi_1,\vdiv\bzeta)$. {Employ} Young's inequality to arrive at $\int_{\Omega_1}u\,\vdiv\bp\ge-\frac12\|u\|_{0,\Omega_1}^2-\frac12\|\vdiv\bp\|_{0,\Omega_1}^2$, $-\int_{\Omega_1}u\,\vdiv\bzeta\ge-\|u\|_{0,\Omega_1}^2-\frac14\|\vdiv\bzeta\|_{0,\Omega_1}^2$ and $-\int_{\Omega_1}\varphi_1\,\vdiv\bzeta\ge-\|\varphi_1\|_{0,\Omega_1}^2-\frac14\|\vdiv\bzeta\|_{0,\Omega_1}^2$. A combination of all this yields
\begin{align*}
\mathcal{M}(\vec{\bx};\vec{\by}_3)&=\|\vdiv\bp\|_{0,\Omega_1}^2+\|\vdiv\bzeta\|_{0,\Omega_1}^2+\int_{\Omega_1}u\,\vdiv\bp-\int_{\Omega_1}\varphi_1\,\vdiv\bzeta-\int_{\Omega_1}u\,\vdiv\bzeta\\
&\ge\tfrac12\|\vdiv\bp\|_{0,\Omega_1}^2+\tfrac12\|\vdiv\bzeta\|_{0,\Omega_1}^2-\tfrac32\|u\|_{0,\Omega_1}^2-\|\varphi_1\|_{0,\Omega_1}^2 .
\end{align*}
Set $\vec{\by}:=\vec{\by}_1+\theta_1\vec{\by}_2+\theta_2\vec{\by}_3$ with $\theta_1:=\beta^2/2$ and $\theta_2:=\beta^2/8$. Adding (i)--(iii), the coefficients of the squared norms are: $1-\frac{\theta_1}{2\beta^2}=\frac34$ for $\|\bp\|_{0,\Omega_1}^2$; $1$ for $\delta_c^{-2}\|\bzeta\|_{0,\Omega_1}^2$ and $\|\varphi_2\|_{W,\varepsilon_r}^2$; $1-\frac32\theta_2=1-\frac{3\beta^2}{16}$ for $\|u\|_{0,\Omega_1}^2$; $\frac{\theta_1}{2}-\theta_2=\frac{\beta^2}{8}$ for $\|\varphi_1\|_{0,\Omega_1}^2$; and $\frac{\theta_2}{2}=\frac{\beta^2}{16}$ for $\|\vdiv\bp\|_{0,\Omega_1}^2$ and $\|\vdiv\bzeta\|_{0,\Omega_1}^2$. Hence
\begin{align*}
\mathcal{M}(\vec{\bx};\vec{\by})&\ge\tfrac34\|\bp\|_{0,\Omega_1}^2+\delta_c^{-2}\|\bzeta\|_{0,\Omega_1}^2+\|\varphi_2\|_{W,\varepsilon_r}^2+\Big(1-\tfrac{3\beta^2}{16}\Big)\|u\|_{0,\Omega_1}^2\\
&\quad +\tfrac{\beta^2}{8}\|\varphi_1\|_{0,\Omega_1}^2+\tfrac{\beta^2}{16}\big(\|\vdiv\bp\|_{0,\Omega_1}^2+\|\vdiv\bzeta\|_{0,\Omega_1}^2\big)\\
&\ge\tfrac{\beta^2}{16}\normiii{\vec{\bx}}^2 ,
\end{align*}
where in the last step we used $\beta\le1$, so that each coefficient is at least $\beta^2/16$ (in particular $1-\frac{3\beta^2}{16}\ge\frac{13\beta^2}{16}\ge\frac{\beta^2}{16}$). On the other hand, $\normiii{\vec{\by}_1}=\normiii{\vec{\bx}}$, $\normiii{\vec{\by}_2}=\|\widehat{\bq}\|_{\bQ}\le\beta^{-1}\|\varphi_1\|_{0,\Omega_1}\le\beta^{-1}\normiii{\vec{\bx}}$ and $\normiii{\vec{\by}_3}=(\|\vdiv\bp\|_{0,\Omega_1}^2+\|\vdiv\bzeta\|_{0,\Omega_1}^2)^{1/2}\le\normiii{\vec{\bx}}$, so that, by the triangle inequality, $\normiii{\vec{\by}}\le(1+\theta_1\beta^{-1}+\theta_2)\normiii{\vec{\bx}}=(1+\beta/2+\beta^2/8)\normiii{\vec{\bx}}$. Consequently,
\[
\sup_{\vec{\by}'\in\vec{\bY}\setminus\{\bm0\}}\frac{\mathcal{M}(\vec{\bx};\vec{\by}')}{\normiii{\vec{\by}'}}\ge\frac{\mathcal{M}(\vec{\bx};\vec{\by})}{\normiii{\vec{\by}}}\ge\frac{\beta^2/16}{1+\beta/2+\beta^2/8}\normiii{\vec{\bx}}=\gamma\normiii{\vec{\bx}},
\]
which proves \eqref{iinfsup}$_1$ (if $\vec{\bx}=\bm0$ there is nothing to prove, and otherwise $\vec{\by}\ne\bm0$ because $\mathcal M(\vec{\bx};\vec{\by})>0$).

For the second one, given $\vec{\by}=(\bq,\bxi,\psi_1,v,\psi_2)$ we take $\vec{\bx}:=\vec{\bx}_1+\theta_1\vec{\bx}_2+\theta_2\vec{\bx}_3$ with $\vec{\bx}_1:=(\bq,\bxi,-\psi_1,-v,\psi_2)$, $\vec{\bx}_2:=(\tilde{\bq},\bm0,0,0,0)$, and $\vec{\bx}_3:=(\bm0,\bm0,\vdiv\bq,\vdiv\bxi,0)$ with $\tilde{\bq}\in\bQ_\Gamma$ such that $\vdiv\tilde{\bq}=\psi_1$ and $\|\tilde{\bq}\|_{\bQ}\le\beta^{-1}\|\psi_1\|_{0,\Omega_1}$ (\Cref{infsupb}). The same cancellations as in (i) yield
\[
\mathcal M(\vec{\bx}_1;\vec{\by})=\|\bq\|_{0,\Omega_1}^2+\delta_c^{-2}\|\bxi\|_{0,\Omega_1}^2+\|\psi_2\|_{W,\varepsilon_r}^2+\|v\|_{0,\Omega_1}^2,
\]
and, since $b_\Gamma(\tilde{\bq},\psi_2)=0$,
\begin{align*}
\mathcal M(\vec{\bx}_2;\vec{\by})&=\int_{\Omega_1}\tilde{\bq}\cdot\bq+\|\psi_1\|_{0,\Omega_1}^2\ge\tfrac12\|\psi_1\|_{0,\Omega_1}^2-\tfrac12\beta^{-2}\|\bq\|_{0,\Omega_1}^2,\\
\mathcal M(\vec{\bx}_3;\vec{\by})&=\|\vdiv\bq\|_{0,\Omega_1}^2+\|\vdiv\bxi\|_{0,\Omega_1}^2+\int_{\Omega_1}\vdiv\bxi\,\psi_1-\int_{\Omega_1}\vdiv\bxi\,v-\int_{\Omega_1}\vdiv\bq\,v\\
&\ge\tfrac12\|\vdiv\bq\|_{0,\Omega_1}^2+\tfrac12\|\vdiv\bxi\|_{0,\Omega_1}^2-\|\psi_1\|_{0,\Omega_1}^2-\tfrac32\|v\|_{0,\Omega_1}^2,
\end{align*}
where we used Young's inequality as in (iii) (with weights $\frac14$ on $\|\vdiv\bxi\|_{0,\Omega_1}^2$ in the first two cross terms and $\frac12$ in the third one). These are the bounds (i)--(iii) with the roles of $(\varphi_1,u)$ and $(\psi_1,v)$ exchanged, so the same choice of $\theta_1,\theta_2$ yields $\mathcal M(\vec{\bx};\vec{\by})\ge\frac{\beta^2}{16}\normiii{\vec{\by}}^2$ and $\normiii{\vec{\bx}}\le(1+\beta/2+\beta^2/8)\normiii{\vec{\by}}$.
\end{proof}

\begin{theorem}[Well-posedness]\label{th:wellposed}
Under the assumptions of \Cref{sec:model}, problem \eqref{Pw} has a unique solution $\vec{\bx}\in\vec{\bY}$, and $\normiii{\vec{\bx}}\le\gamma^{-1}C_l$, with $\gamma$ and $C_l$ as in \Cref{stability} and \Cref{continuity}. Consequently, \eqref{IPw} has a unique solution, which does not depend on the choice of the liftings in \eqref{eq:compat}.
\end{theorem}
\begin{proof}
It is a direct consequence of \Cref{continuity}, \Cref{stability}, and the Banach--Ne\v{c}as--Babu\v{s}ka theorem (see, e.g., \cite{Ern2}). In more detail, \Cref{continuity} {yields} the boundedness of $\mathcal M$ and $l$, the first inequality in \eqref{iinfsup} gives injectivity and closed range of the operator induced by $\mathcal M$, and the second one shows that for every $\vec{\by}\ne\bm0$ there is $\vec{\bx}$ with $\mathcal M(\vec{\bx};\vec{\by})\ne0$, which yields surjectivity. The a priori bound follows from
\[
\gamma\normiii{\vec{\bx}}\le\sup_{\vec{\by}\in\vec{\bY}\setminus\{\bm0\}}\frac{\mathcal M(\vec{\bx};\vec{\by})}{\normiii{\vec{\by}}}=\sup_{\vec{\by}\in\vec{\bY}\setminus\{\bm0\}}\frac{l(\vec{\by})}{\normiii{\vec{\by}}}\le C_l.
\]
Finally, let $(\tilde g_1,\tilde g_2)$ and $(\tilde g_1',\tilde g_2')$ be two liftings satisfying \eqref{eq:compat}, with associated $g_D$ and $g_D'$. Then $\tilde g_2-\tilde g_2'\in W_0$, and $g_D-g_D'=E_\Gamma(\tilde g_2-\tilde g_2')$, since both sides vanish on $\partial\Omega_1\setminus\Gamma$ and coincide with $\tilde g_2-\tilde g_2'$ on $\Gamma$. Hence $g_D+E_\Gamma(\varphi_2-\tilde g_2)=g_D'+E_\Gamma(\varphi_2-\tilde g_2')$, so that the first equation in \eqref{IPw}, which reads $a_1(\bp,\bq)+b(\bq,\varphi_1)=\langle\bq\cdot\bm n_1,g_D+E_\Gamma(\varphi_2-\tilde g_2)\rangle_{\partial\Omega_1}$ is the same for both liftings. The other equations in \eqref{IPw} do not involve the liftings, and the constraint $\varphi_2-\tilde g_2\in W_0$ is equivalent to $\varphi_2-\tilde g_2'\in W_0$. Therefore both liftings lead to the same problem \eqref{IPw}, with a unique solution.
\end{proof}

\begin{remark}[Dependence on the parameters]\label{rem:robust}
The inf-sup constant $\gamma$ depends only on $\Omega_1$ and $\Gamma$, and the continuity constant $C_{\mathcal M}=3+2C_\Gamma$ depends on $\varepsilon_r$ only through the trace inequality \eqref{eq:trace}, as $\varepsilon_r^{-1/2}$ for $\varepsilon_r<1$; neither depends on $\delta_c$. The data are in $C_l$, which is not uniform in the parameters since it contains $C_\Gamma\|g\|_{H^{1/2}_{00}(\Gamma)'}$, which grows as $\varepsilon_r^{-1/2}$ for $\varepsilon_r<1$ when $g\neq0$, the term $\max\{1,\delta_c\}\|g_3\|_{1/2,\partial\Omega_1}$, which grows with $\delta_c>1$ when $g_3\neq0$, and the lifting-dependent term $\|\tilde g_2\|_{W,\varepsilon_r}$. Thus, parameter-robustness refers to the inf-sup constant $\gamma$. Note also that the bound in \Cref{th:wellposed} is for $\varphi_2^0=\varphi_2-\tilde g_2$; the bound for $\varphi_2$ itself follows by adding $\|\tilde g_2\|_{W,\varepsilon_r}$. This is in contrast with a stability analysis in the unweighted $\bH(\mathrm{div})\times L^2\times H^1$-norms, whose constants behave like $\min\{1,\delta_c^{-2},\varepsilon_r\}$ and $\max\{1,\delta_c^{-2}\}$. Since $\varepsilon_r=\varepsilon_2/\varepsilon_1$ is typically small and $\delta_c$ varies over orders of magnitude with the ionic concentration (see \Cref{sec:numerical}), this robustness is relevant in practice. The dependence of $C_{\mathcal M}$ on $\varepsilon_r$ does not propagate to the error estimates of \Cref{sec:error}, up to a loss of at most half an order in the pre-asymptotic regime (see \Cref{errthm}(b), \Cref{cor:rates} and \Cref{rem:errrobust}).
\end{remark}

\begin{remark}[On the choice of the weights]\label{rem:weights}
The weights in $\normiii{\bullet}$ cannot be modified without losing the uniformity of either \Cref{continuity} or \Cref{stability}. To see this, consider first the norm $\normiii{\bullet}_w$ obtained from $\normiii{\bullet}$ by replacing, both for trial and test functions, the weight of $\|\vdiv\bxi\|_{0,\Omega_1}^2$ by a constant $w>0$ (for instance, $w=\delta_c^{-2}$ corresponds to measuring $\bzeta$ in $\delta_c^{-1}\|\bullet\|_{\bQ}$), and let $\bzeta_m:=(m^{-1}\sin(mx_1),0,\ldots,0)^{\tt t}\in\bQ$, $m\in\mathbb N$, for which $\vdiv\bzeta_m=\cos(mx_1)$, so that $\|\bzeta_m\|_{0,\Omega_1}\to0$ and $\|\vdiv\bzeta_m\|_{0,\Omega_1}^2\to|\Omega_1|/2$ as $m\to\infty$. For $\vec{\bx}=(\bm0,\bm0,0,\vdiv\bzeta_m,0)$ and $\vec{\by}=(\bm0,\bzeta_m,0,0,0)$ we have $\mathcal M(\vec{\bx};\vec{\by})=b(\bzeta_m,\vdiv\bzeta_m)=\|\vdiv\bzeta_m\|_{0,\Omega_1}^2$, $\normiii{\vec{\bx}}_w=\|\vdiv\bzeta_m\|_{0,\Omega_1}$ and $\normiii{\vec{\by}}_w=(\delta_c^{-2}\|\bzeta_m\|_{0,\Omega_1}^2+w\|\vdiv\bzeta_m\|_{0,\Omega_1}^2)^{1/2}$, which shows, letting $m\to\infty$, that the continuity constant of $\mathcal M$ with respect to $\normiii{\bullet}_w$ is at least $w^{-1/2}$. On the other hand, for $\vec{\bx}=(\bm0,\bzeta,0,0,0)$ and any $\vec{\by}\in\vec{\bY}$, the definition of $\mathcal M$ and the Cauchy--Schwarz inequality give
\[
\mathcal M(\vec{\bx};\vec{\by})=\delta_c^{-2}\int_{\Omega_1}\bzeta\cdot\bxi+\int_{\Omega_1}\vdiv\bzeta\,v\le\big(\delta_c^{-1}\|\bzeta\|_{0,\Omega_1}+\|\vdiv\bzeta\|_{0,\Omega_1}\big)\normiii{\vec{\by}}_w,
\]
whereas $\normiii{\vec{\bx}}_w\ge w^{1/2}\|\vdiv\bzeta\|_{0,\Omega_1}$. Taking $\bzeta=\bzeta_m$ and letting $m\to\infty$ gives an inf-sup constant (in $\normiii{\bullet}_w$)  at most $w^{-1/2}$. Hence, uniform continuity and inf-sup constants require $w\simeq1$; in particular, $w=\delta_c^{-2}$ leads to an inf-sup constant of order up to $\delta_c$. Similarly, if the weight of $\|u\|_{0,\Omega_1}^2$ and of $\|v\|_{0,\Omega_1}^2$ is replaced by $w_u>0$, and $\normiii{\bullet}_{w_u}$ denotes the resulting norm, then taking $\vec{\bx}=(\bm0,\bm0,0,u,0)$ and $\vec{\by}=(\bm0,\bm0,0,-u,0)$, for which $\mathcal M(\vec{\bx};\vec{\by})=-a(u,-u)=\|u\|_{0,\Omega_1}^2$, shows that the continuity constant is at least $w_u^{-1}$, while for $\vec{\bx}=(\bm0,\bm0,0,u,0)$ one has $\mathcal M(\vec{\bx};\vec{\by})=\int_{\Omega_1}u\,(\vdiv\bxi+\psi_1-v)\le(1+1+w_u^{-1/2})\|u\|_{0,\Omega_1}\normiii{\vec{\by}}_{w_u}$ and $\normiii{\vec{\bx}}_{w_u}=w_u^{1/2}\|u\|_{0,\Omega_1}$, so that the inf-sup constant is at most $(2+w_u^{-1/2})w_u^{-1/2}$. Therefore $w_u\simeq1$ as well; in particular, a weight $\delta_c^{2}$ on $\|u\|_{0,\Omega_1}^2$ cannot compensate a weight $\delta_c^{-2}$ on $\|\vdiv\bxi\|_{0,\Omega_1}^2$, since it yields a continuity constant of order $\delta_c^{-2}$. A similar  argument applies to $\varphi_1$ and $\psi_1$, which are coupled with $u$ and $v$ through the parameter-free terms $a(u,\psi_1)$ and $a(\varphi_1,v)$. Finally, the weight $\delta_c^{-2}$ of $\|\bxi\|_{0,\Omega_1}^2$ is dictated by $a_2$. Note also that, by \eqref{IPM}, $\vdiv\bzeta=u+\varphi_1-f_1$, so that $\|\vdiv\bzeta\|_{0,\Omega_1}$ is controlled by the data uniformly in $\delta_c$, whereas $\delta_c^{-1}\|\vdiv\bzeta\|_{0,\Omega_1}$ is not (for instance, if $u$ has a boundary layer of width $\delta_c$).
\end{remark}

\section{Discrete formulation and well-posedness}\label{sec:discrete}
Consider a \textit{shape-regular simplicial triangulation} $\mathcal{T}_{h}^{}$ of $\overline{\Omega}$, fitted to the interface: every $T\in\mathcal T_h$ is contained in either $\overline{\Omega}_1$ or $\overline{\Omega}_2$, so that $\Gamma$ is a union of facets of $\mathcal T_h$. For curved interfaces this means that $\Gamma$ (and $\Omega_1$, $\Omega_2$) are replaced by polygonal/polyhedral approximations; see Example 4 in \Cref{sec:numerical}. Accordingly, throughout this section and \Cref{sec:error} we assume that $\Omega_1$ and $\Omega_2$ are polygonal/polyhedral, so that $\Gamma$ is a finite union of flat faces; for a curved $\Gamma$ the analysis applies to the problem posed on the polygonal/polyhedral approximation of the subdomains, and the geometric error is not analysed. Let $T\in \mathcal{T}_{h}^{}$ denote an element and $\partial T$ be its boundary.

For integers $k\geq 0$, let $\mathbb{P}_{k}^{}(T)$ denote the space of polynomials of total degree at most $k$ in $\mathbb{R}^d$ and $\mathbb{P}_{k}^{H}(T)$ be the space of homogeneous polynomials of degree $k$ in $\mathbb{R}^d$. Then the local Raviart--Thomas finite element space $\mathbb{RT}_{k}(T)$ is defined by
\begin{equation*}
    \mathbb{RT}_{k}^{}(T):=[\mathbb{P}_{k}^{}(T)]^d\oplus\bx\,\mathbb{P}_{k}^{H}(T).
\end{equation*}
Note that $\mathbb{RT}_k(T)\subset[\mathbb{P}_{k+1}(T)]^d$. 
Given integers $k_{1}^{}\geq 0$ and $k_{2}^{}\geq 1$, {we define} the following finite-dimensional subspaces of $\bQ$, $V$, $W$, and $W_{0}^{}$, respectively:
\begin{align*}
        \bQ_{h}^{}&:=\big\{\bq_{h}^{}\in\bQ:\bq_{h}^{}\big|_{T}\in\mathbb{RT}_{k_{1}^{}}^{}(T)\quad\text{for all } T\in\mathcal{T}_{h}^{},\ T\subset\overline\Omega_1\big\},\\
        V_{h}^{}&:=\big\{v_{h}^{}\in V: v_{h}^{}\big|_{T}\in\mathbb{P}_{k_{1}^{}}^{}(T)\quad\text{for all } T\in\mathcal{T}_{h}^{},\ T\subset\overline\Omega_1\big\},\\
        W_{h}^{}&:=\big\{w_{h}^{}\in W: w_{h}^{}\big|_{T}\in\mathbb{P}_{k_{2}^{}}^{}(T)\quad\text{for all } T\in\mathcal{T}_{h}^{},\ T\subset\overline\Omega_2\big\},\qquad
        W_{h}^{0}:=W_h\cap W_0,
\end{align*}
and $\vec{\bY}_h:=\bQ_h\times\bQ_h\times V_h\times V_h\times W_h^0\subset\vec{\bY}$. For simplicity of the analysis, we assume that the lifting $\tilde g_2$ in \eqref{eq:compat} can be taken in $W_h$, i.e., that $g_2$ is the trace of a function of $W_h$ (e.g., $g_2=0$). In general, $g_2$ is replaced by its Lagrange interpolant $\mathcal I_hg_2$ on $\partial\Omega_2\setminus\Gamma$, which adds a standard consistency term to the error estimates below.
The discrete counterpart of \eqref{IPw} is: find $(\bp_{h}^{},\bzeta_{h}^{}, \varphi_{1,h}^{}, u_{h}^{}, \varphi_{2,h}^{}) \in \bQ_{h}^{}\times \bQ_{h}^{}\times V_{h}^{}\times V_{h}^{}\times W_{h}^{}$ with $\varphi_{2,h}^{}-\tilde g_2\in W_h^0$ such that
\begin{align}\label{IPh}
    a_{1}^{}(\bp_{h}^{},\bq_{h}^{}) + b(\bq_h,\varphi_{1,h}^{}) - b_{\Gamma}(\bq_h,\varphi_{2,h}^{}-\tilde g_2) &= l_{1}^{}(\bq_{h}^{}) &&\text{ for all } \bq_{h}^{}\in \bQ_{h}^{},\nonumber\\
    a_{2}^{}(\bzeta_{h}^{},\bxi_{h}^{}) + b(\bxi_h,u_{h}^{}) &= l_{2}^{}(\bxi_{h}^{})&&\text{ for all } \bxi_{h}^{}\in \bQ_{h}^{},\nonumber\\
    b(\bp_h,\psi_{1,h}^{}) + a(u_{h}^{},\psi_{1,h}^{}) &= 0&&\text{ for all } \psi_{1,h}^{}\in V_{h}^{},\\
    b(\bzeta_h,v_{h}^{}) - a(u_{h}^{},v_{h}^{}) - a(\varphi_{1,h}^{},v_{h}^{})& = l_{4}^{}(v_{h}^{})&&\text{ for all } v_{h}^{}\in V_{h}^{},\nonumber\\
    a_{3}^{}(\varphi_{2,h}^{},\psi_{2,h}^{}) + b_{\Gamma}(\bp_h,\psi_{2,h}^{})& = l_{5}^{}(\psi_{2,h}^{})&&\text{ for all } \psi_{2,h}^{}\in W_{h}^{0}.\nonumber
\end{align}
As in the continuous case, writing $\varphi_{2,h}=\varphi_{2,h}^0+\tilde g_2$, \eqref{IPh} is equivalent to: find $\vec{\bx}_h=(\bp_h,\bzeta_h,\varphi_{1,h},u_h,\varphi_{2,h}^0)\in\vec{\bY}_h$ such that
\begin{equation}\label{Ph}
\mathcal{M}(\vec{\bx}_{h};\vec{\by}_{h}^{})=l(\vec{\by}_{h}^{})\qquad\text{ for all }\vec{\by}_{h}^{}\in \vec{\bY}_{h}^{}.
\end{equation}

Now we address the discrete analogue of \Cref{infsupb}. Let $\bQ_{h,\Gamma}:=\bQ_h\cap\bQ_\Gamma$, the RT fluxes with vanishing normal component on $\Gamma$.  
\begin{lemma}[{Discrete surjectivity of the divergence}]\label{dinfsupb}
There exists $\beta_\sharp>0$, independent of $h$, such that for every $v_h\in V_h$ there is $\widehat{\bq}_h\in\bQ_{h,\Gamma}$ with $\vdiv\widehat{\bq}_h=v_h$ and $\|\widehat{\bq}_h\|_{\bQ}\le\beta_\sharp^{-1}\|v_h\|_{0,\Omega_1}$.
\end{lemma}
\begin{proof}
{This is also standard (see, e.g., \cite[Lemma 51.10 and Remark 51.15]{Ern2}): take $\widehat{\bq}_h$ as the canonical Raviart--Thomas interpolant of $\widehat{\bq}$ from \Cref{infsupb} with $v=v_h$, which is   uniformly bounded thanks to the $H^{1+s}$-regularity of $w$ for some $s>0$ \cite{Dauge1992,HallerDintelmannMeinlschmidtWollner2019} and \cite[Theorem 16.4]{Ern1}, and satisfies $\vdiv\widehat{\bq}_h=v_h$ and $\widehat{\bq}_h\cdot\bm n_1=0$ on $\Gamma$. Hence {there exists such $\beta_\sharp\ \le1$ and it} depends only on $\Omega_1$, $\Gamma$, and the mesh shape regularity.}
\end{proof}

\begin{theorem}[Discrete well-posedness]\label{th:dwellposed}
Let $\gamma_\sharp:=\beta_\sharp^2/(16+8\beta_\sharp+2\beta_\sharp^2)$. Then
\begin{equation}\label{idinfsup}
\sup_{\vec{\by}_h\in\vec{\bY}_h\setminus\{\bm0\}}\frac{\mathcal{M}(\vec{\bx}_h;\vec{\by}_h)}{\normiii{\vec{\by}_h}}\ge\gamma_\sharp\normiii{\vec{\bx}_h}\qquad\text{for all }\vec{\bx}_h\in\vec{\bY}_h,
\end{equation}
and \eqref{Ph} (equivalently, \eqref{IPh}) has a unique solution, which satisfies $\normiii{\vec{\bx}_h}\le\gamma_\sharp^{-1}C_l$. The constant $\gamma_\sharp$ is independent of $h$, $\delta_c$ and $\varepsilon_r$.
\end{theorem}
\begin{proof}
The proof of \Cref{stability} carries over verbatim: $\vec{\by}_1$ is in $\vec{\bY}_h$; $\widehat{\bq}$ is replaced by $\widehat{\bq}_h$ from \Cref{dinfsupb}; and $\vec{\by}_3\in\vec{\bY}_h$ because $\vdiv\bQ_h\subseteq V_h$. The computations (i)--(iii) only utilise the identities $b_\Gamma(\widehat{\bq}_h,\varphi_{2,h})=0$ and $\vdiv\widehat{\bq}_h=\varphi_{1,h}$, the bound $\|\widehat{\bq}_h\|_{\bQ}\le\beta_\sharp^{-1}\|\varphi_{1,h}\|_{0,\Omega_1}$ and $\beta_\sharp\le1$, so that they hold with $\beta$ replaced by $\beta_\sharp$. Since $\vec{\bY}_h$ is finite dimensional, \eqref{idinfsup} implies unique solvability (injectivity of a linear map between spaces of equal finite dimension implies surjectivity). The bound follows as in the proof of \Cref{th:wellposed}, since $|l(\vec{\by}_h)|\le C_l\normiii{\vec{\by}_h}$ for all $\vec{\by}_h\in\vec{\bY}_h\subset\vec{\bY}$.
\end{proof}

\section{Error analysis of the mixed-primal finite element scheme}\label{sec:error}
We utilise the $L^2$-projection $P_{h}^{}:V\to V_{h}^{}$, the canonical Raviart--Thomas interpolant $\bPi_{h}^{}:\bQ\cap\bH^s(\Omega_1)\to\bQ_{h}^{}$, $s>0$, which satisfies the commuting property \cite[Section 2.5]{Boffi}\begin{equation}\label{eq:divcommute}
        \vdiv\bPi_h\bq=P_h\vdiv\bq,
    \end{equation} and the Lagrange interpolant $\mathcal{I}_{h}^{}:H^{k_2+1}(\Omega_{2}^{})\to W_{h}^{}$, which preserves the boundary values in $W_h$.
These canonical operators satisfy the approximation properties
\begin{subequations}\begin{align}
    &\|v-P_hv\|_{0,T}\lesssim h_T^{k_1+1}|v|_{k_1+1,T} \quad\forall v\in H^{k_1+1}(T), T\in\mathcal T_h\text{ with }T\subset\overline\Omega_1,\label{eq:L2approx}\\
    &\|\bq-\bPi_h\bq\|_{0,T}\lesssim h_T^{k_1+1}|\bq|_{k_1+1,T} \quad\forall \bq\in \bH^{k_1+1}(T),T\in\mathcal T_h\text{ with }T\subset\overline\Omega_1,\label{eq:Hdivapprox}\\
    &\|\bq-\bPi_h\bq\|_{0,4,T}\lesssim h_T^{k_1+1-d/4}|\bq|_{k_1+1,T} \quad\forall \bq\in \bH^{k_1+1}(T),T\in\mathcal T_h\text{ with }T\subset\overline\Omega_1,\label{eq:L4approx}\\
    &\|\psi_2-\mathcal I_h\psi_2\|_{0,T}+h_T|\psi_2-\mathcal I_h\psi_2|_{1,T}\lesssim h_T^{k_2+1}|\psi_2|_{k_2+1,T} \ \forall \psi_2\in H^{k_2+1}(T),T\in\mathcal T_h\text{ with }T\subset\overline\Omega_2.\label{eq:H1approx}
\end{align}\end{subequations}
These follow from \cite[Chapter 11]{Ern1}, \cite[Theorem 16.4]{Ern1}, \cite[Chapter 11]{Ern1} and \cite[Chapter 11]{Ern1}, respectively.

For $\vec{\bx}=(\bp,\bzeta,\varphi_1,u,\varphi_2)$ we set $\mathcal P_h\vec{\bx}:=(\bPi_h\bp,\bPi_h\bzeta,P_h\varphi_1,P_hu,\mathcal I_h\varphi_2)$. Here $\varphi_2$ denotes the full potential: since $\tilde g_2\in W_h$ and $\mathcal I_h\tilde g_2=\tilde g_2$, for the homogenised unknown $\varphi_2^0=\varphi_2-\tilde g_2$ the fifth component is $\mathcal I_h\varphi_2-\tilde g_2$, and the interpolation error is the same. 

We also need the action of $\bPi_h$ on the interface. Let $\mathcal F_h^\Gamma$ be the set of facets of $\mathcal T_h$ contained in $\Gamma$ (recall that $\mathcal T_h$ is fitted to $\Gamma$, so that $\overline\Gamma=\bigcup_{F\in\mathcal F_h^\Gamma}\overline F$). For $F\in\mathcal F_h^\Gamma$, we denote by $T_F$ the unique element of $\mathcal T_h$ such that $T_F\subset\overline\Omega_2$ and $F$ is a facet of $T_F$, by $h_F$ and $h_{T_F}$ the diameters of $F$ and $T_F$ (so that $h_F\le h_{T_F}$), and by $Q_F:L^2(F)\to\mathbb P_{k_1}(F)$ the $L^2(F)$-orthogonal projection. Since the degrees of freedom of $\mathbb{RT}_{k_1}(T)$ attached to a facet $F$ of $T$ are the moments $\int_F\bq\cdot\bm n_F\,\mu\ds$, $\mu\in\mathbb P_{k_1}(F)$, and since the normal component on $F$ of any element of $\mathbb{RT}_{k_1}(T)$ belongs to $\mathbb P_{k_1}(F)$ \cite[Section 2.5]{Boffi}, there holds
\begin{equation}\label{eq:RTtrace}
(\bPi_h\bq)\cdot\bm n_1|_F=Q_F\big(\bq\cdot\bm n_1|_F\big)\qquad\text{for all }F\in\mathcal F_h^\Gamma,
\end{equation}
for every $\bq$ in the domain of $\bPi_h$ such that $\bq\cdot\bm n_1|_F\in L^2(F)$ for all $F\in\mathcal F_h^\Gamma$. We also write $h_\Gamma:=\max_{F\in\mathcal F_h^\Gamma}h_{T_F}\le h$ and, for $w$ such that $w|_F\in H^{m}(F)$ for all $F\in\mathcal F_h^\Gamma$, $|w|_{m,\Gamma_h}:=\big(\sum_{F\in\mathcal F_h^\Gamma}|w|_{m,F}^2\big)^{1/2}$.

\begin{lemma}[Approximation properties]\label{iestimatelemma}
Let $k_1\ge0$, $k_2\ge1$, and assume $\bp,\bzeta\in\bH^{k_1+1}(\Omega_1)$, $\vdiv\bzeta,\varphi_1,u\in H^{k_1+1}(\Omega_1)$ and $\varphi_2\in H^{k_2+1}(\Omega_2)$. Then
\begin{align}
\normiii{\vec{\bx}-\mathcal P_h\vec{\bx}}&\lesssim h^{k_1+1}\Big(|\bp|_{k_1+1,\Omega_1}+|\vdiv\bp|_{k_1+1,\Omega_1}+\delta_c^{-1}|\bzeta|_{k_1+1,\Omega_1}+|\vdiv\bzeta|_{k_1+1,\Omega_1}\nonumber\\
&\qquad +|\varphi_1|_{k_1+1,\Omega_1}+|u|_{k_1+1,\Omega_1}\Big) +h^{k_2}\big(\varepsilon_r^{1/2}+h\big)|\varphi_2|_{k_2+1,\Omega_2},\label{ierrestimate}
\end{align}
where the hidden constant depends only on the shape regularity of $\mathcal T_h$ and on $k_1,k_2$. Finally, if $(\bp\cdot\bm n_1)|_F\in H^{k_1+1}(F)$ for all $F\in\mathcal F_h^\Gamma$, then
\begin{equation}\label{eq:facetapprox}
\|(\bp-\bPi_h\bp)\cdot\bm n_1\|_{0,F}\le C_A\,h_F^{k_1+1}|\bp\cdot\bm n_1|_{k_1+1,F}\qquad\text{for all }F\in\mathcal F_h^\Gamma,
\end{equation}
with $C_A>0$ depending only on $d$, $k_1$ and the shape regularity of $\mathcal T_h$.
\end{lemma}
\begin{proof}
We argue element by element and then add the local contributions. Let $T\in\mathcal T_h$ with $T\subset\overline\Omega_1$. 
{We apply \eqref{eq:L2approx}} 
to $v\in\{\varphi_1,u,\vdiv\bp,\vdiv\bzeta\}$; for the last two choices, the commuting property \eqref{eq:divcommute} 
shows that $\vdiv(\bq-\bPi_h\bq)=(I-P_h)\vdiv\bq$ for $\bq\in\{\bp,\bzeta\}$. Similarly, {we employ\eqref{eq:Hdivapprox}} 
for $\bq\in\{\bp,\bzeta\}$, and the {estimate for $\bzeta$} 
is multiplied by $\delta_c^{-1}$ in $\normiii{\bullet}$. Next, let $T\in\mathcal T_h$ with $T\subset\overline\Omega_2$. Since $k_2+1\ge2>d/2$, $H^{k_2+1}(T)\subset C^0(\overline T)$ and {by utilising \eqref{eq:H1approx} yields} 
\[
\|\varphi_2-\mathcal I_h\varphi_2\|_{W,\varepsilon_r}\le\varepsilon_r^{1/2}|\varphi_2-\mathcal I_h\varphi_2|_{1,\Omega_2}+\|\varphi_2-\mathcal I_h\varphi_2\|_{0,\Omega_2}\lesssim(\varepsilon_r^{1/2}h^{k_2}+h^{k_2+1})|\varphi_2|_{k_2+1,\Omega_2}.
\]
Squaring the local estimates, adding them over the elements, and recalling the definition of $\normiii{\bullet}$ gives \eqref{ierrestimate}. Finally, given $F\in\mathcal F_h^\Gamma$, \eqref{eq:RTtrace} yields $(\bp-\bPi_h\bp)\cdot\bm n_1|_F=(I-Q_F)(\bp\cdot\bm n_1|_F)$, and \eqref{eq:facetapprox} is the Bramble--Hilbert estimate for the $L^2(F)$-orthogonal projection onto $\mathbb P_{k_1}(F)$ on the $(d-1)$-dimensional simplex $F$, whose shape regularity is inherited from that of the element containing it \cite[Chapter 11]{Ern1}.
\end{proof}

In the error analysis, the interface consistency term will be measured by
\begin{equation}\label{eq:Thetadef}
\Theta_h(\bq):=\sup_{\psi_h\in W_h^0\setminus\{0\}}\frac{|b_\Gamma(\bq-\bPi_h\bq,\psi_h)|}{\|\psi_h\|_{W,\varepsilon_r}},
\end{equation}
defined for $\bq$ in the domain of $\bPi_h$.

\begin{theorem}[Error bounds]\label{errthm}
Let $\vec{\bx}$ and $\vec{\bx}_h$ be the solutions of \eqref{Pw} and \eqref{Ph}, respectively, and assume that $\bp,\bzeta\in\bH^s(\Omega_1)$ for some $s>0$ and $\varphi_2\in C^0(\overline\Omega_2)$, so that $\mathcal P_h\vec{\bx}$ is well defined. Then the following estimates hold:
\begin{equation}\label{error}
\mathrm{(a)} \qquad \text{(Quasi-{best approximation})} 
\qquad \normiii{\vec{\bx}- \vec{\bx}_h}\le\Big(1+{\gamma_\sharp^{-1}C_{\mathcal M}}\Big)\normiii{\vec{\bx}-\mathcal P_h\vec{\bx}},
\end{equation}
and, under the regularity assumptions of \Cref{iestimatelemma}, $\normiii{\vec{\bx}-\vec{\bx}_h}\lesssim(1+C_\Gamma)\big(h^{k_1+1}+h^{k_2}(\varepsilon_r^{1/2}+h)\big)$, where the hidden constant depends on the norms of the solution in \eqref{ierrestimate} but not otherwise on $\delta_c$ or $\varepsilon_r$. Without the above regularity, the same proof with an arbitrary $\vec{\by}_h\in\vec{\bY}_h$ in place of $\mathcal P_h\vec{\bx}$ gives $\normiii{\vec{\bx}-\vec{\bx}_h}\le(1+{\gamma_\sharp^{-1}C_{\mathcal M}})\inf_{\vec{\by}_h\in\vec{\bY}_h}\normiii{\vec{\bx}-\vec{\by}_h}$. 
\begin{equation}\label{error2}
\mathrm{(b)} \quad \text{(Parameter-robust bound)}
\quad\normiii{\vec{\bx}-\vec{\bx}_h}\le\Big(1+{\gamma_\sharp^{-1}}\Big)\normiii{\vec{\bx}-\mathcal P_h\vec{\bx}}+{\gamma_\sharp^{-1}}\Big(C_{\rm tr}\|\varphi_2-\mathcal I_h\varphi_2\|_{1,\Omega_2}+\Theta_h(\bp)\Big),
\end{equation}
where $\gamma_\sharp$ and $C_{\rm tr}$ (cf. \Cref{th:dwellposed} and \eqref{eq:trace}) are independent of $h$, $\delta_c$, and $\varepsilon_r$.
\end{theorem}
\begin{myproof}
\emph{Proof of (a).} Since $\vec{\bY}_h\subset\vec{\bY}$, we have the Galerkin orthogonality $\mathcal M(\vec{\bx}-\vec{\bx}_h;\vec{\by}_h)=0$ for all $\vec{\by}_h\in\vec{\bY}_h$. Since $\tilde g_2\in W_h$, the last component of $\mathcal P_h\vec{\bx}$ (that is, $\mathcal I_h\varphi_2^0=\mathcal I_h\varphi_2-\tilde g_2$) vanishes on $\partial\Omega_2\setminus\Gamma$, so $\mathcal P_h\vec{\bx}-\vec{\bx}_h\in\vec{\bY}_h$. {It follows from Galerkin orthogonality that $\mathcal M(\mathcal P_h\vec{\bx}-\vec{\bx}_h;\vec{\by}_h)=\mathcal M(\mathcal P_h\vec{\bx}-\vec{\bx};\vec{\by}_h)$}. {Utilising} this along with \eqref{idinfsup} and \Cref{continuity}, we obtain
\[
\gamma_\sharp\normiii{\mathcal P_h\vec{\bx}-\vec{\bx}_h}\le\sup_{\vec{\by}_h\in\vec{\bY}_h\setminus\{\bm0\}}\frac{\mathcal M(\mathcal P_h\vec{\bx}-\vec{\bx};\vec{\by}_h)}{\normiii{\vec{\by}_h}}\le C_{\mathcal M}\normiii{\vec{\bx}-\mathcal P_h\vec{\bx}}.
\]
Then
\[
\normiii{\vec{\bx}-\vec{\bx}_h}\le\normiii{\vec{\bx}-\mathcal P_h\vec{\bx}}+\normiii{\mathcal P_h\vec{\bx}-\vec{\bx}_h}\le\Big(1+{\gamma_\sharp^{-1}C_{\mathcal M}}\Big)\normiii{\vec{\bx}-\mathcal P_h\vec{\bx}},
\]
which is \eqref{error}. The rate follows from \Cref{iestimatelemma}, since $1+{\gamma_\sharp^{-1}C_{\mathcal M}}=1+{\gamma_\sharp^{-1}}(3+2C_\Gamma)\le(1+3{\gamma_\sharp^{-1}})(1+C_\Gamma)$ and $\gamma_\sharp$ does not depend on $\delta_c$, $\varepsilon_r$ or $h$.

\emph{Proof of (b).} We keep the notation $\vec{\boldsymbol\eta}:=\vec{\bx}-\mathcal P_h\vec{\bx}=(\bm\eta_{\bp},\bm\eta_{\bzeta},\eta_{\varphi_1},\eta_u,\eta_{\varphi_2})$, so that $\bm\eta_{\bp}=\bp-\bPi_h\bp$, $\bm\eta_{\bzeta}=\bzeta-\bPi_h\bzeta$, $\eta_{\varphi_1}=\varphi_1-P_h\varphi_1$, $\eta_u=u-P_hu$ and $\eta_{\varphi_2}=\varphi_2^0-\mathcal I_h\varphi_2^0$. Since $\tilde g_2\in W_h$ is reproduced by $\mathcal I_h$, we have $\eta_{\varphi_2}=\varphi_2-\mathcal I_h\varphi_2$ (the interpolation error of the full potential), and, as shown in the proof of (a), $\eta_{\varphi_2}\in W_0$. The difference with respect to (a) is that, instead of bounding $\mathcal M(\vec{\boldsymbol\eta};\vec{\by}_h)$ by \Cref{continuity}, we   observe that several of its terms vanish. Let $\vec{\by}_h=(\bq_h,\bxi_h,\psi_{1,h},v_h,\psi_{2,h})\in\vec{\bY}_h$. Since $\vdiv\bQ_h\subseteq V_h$ and $\eta_{\varphi_1}$ and $\eta_u$ are $L^2(\Omega_1)$-orthogonal to $V_h$, we have $b(\bq_h,\eta_{\varphi_1})=\int_{\Omega_1}\eta_{\varphi_1}\,\vdiv\bq_h=0$ and $b(\bxi_h,\eta_u)=\int_{\Omega_1}\eta_u\,\vdiv\bxi_h=0$. By the commuting property, $\vdiv\bm\eta_{\bp}=(I-P_h)\vdiv\bp$ and $\vdiv\bm\eta_{\bzeta}=(I-P_h)\vdiv\bzeta$, therefore  $b(\bm\eta_{\bp},\psi_{1,h})=b(\bm\eta_{\bzeta},v_h)=0$. Finally, $a(\eta_u,\psi_{1,h})=a(\eta_u,v_h)=a(\eta_{\varphi_1},v_h)=0$ because $\psi_{1,h},v_h\in V_h$. Then, the definition of $\mathcal M$ reduces to
\begin{equation}\label{eq:Meta}
\mathcal M(\vec{\boldsymbol\eta};\vec{\by}_h)=a_1(\bm\eta_{\bp},\bq_h)+a_2(\bm\eta_{\bzeta},\bxi_h)+a_3(\eta_{\varphi_2},\psi_{2,h})-b_\Gamma(\bq_h,\eta_{\varphi_2})+b_\Gamma(\bm\eta_{\bp},\psi_{2,h}).
\end{equation}
We bound the five terms on the right-hand side of \eqref{eq:Meta} separately. First, the Cauchy--Schwarz inequality in $\bL^2(\Omega_1)$ (twice), in the inner product inducing $\|\bullet\|_{W,\varepsilon_r}$, and finally in $\mathbb R^3$, yields
\begin{align*}
&|a_1(\bm\eta_{\bp},\bq_h)+a_2(\bm\eta_{\bzeta},\bxi_h)+a_3(\eta_{\varphi_2},\psi_{2,h})|\\
&\quad\le\|\bm\eta_{\bp}\|_{0,\Omega_1}\|\bq_h\|_{0,\Omega_1}+\delta_c^{-{2}}\|\bm\eta_{\bzeta}\|_{0,\Omega_1}\,
\|\bxi_h\|_{0,\Omega_1}+\|\eta_{\varphi_2}\|_{W,\varepsilon_r}\|\psi_{2,h}\|_{W,\varepsilon_r}\\
&\quad\le\big(\|\bm\eta_{\bp}\|_{0,\Omega_1}^2+\delta_c^{-2}\|\bm\eta_{\bzeta}\|_{0,\Omega_1}^2+\|\eta_{\varphi_2}\|_{W,\varepsilon_r}^2\big)^{1/2}\big(\|\bq_h\|_{0,\Omega_1}^2+\delta_c^{-2}\|\bxi_h\|_{0,\Omega_1}^2+\|\psi_{2,h}\|_{W,\varepsilon_r}^2\big)^{1/2}\\
&\quad\le\normiii{\vec{\boldsymbol\eta}}\,\normiii{\vec{\by}_h}.
\end{align*}
Second, since $\eta_{\varphi_2}\in W_0$, the bound \eqref{eq:normaltrace} and the definition of $C_{\rm tr}$ (cf. paragraph following \eqref{eq:normaltrace}) yield 
\[
|b_\Gamma(\bq_h,\eta_{\varphi_2})|=|\langle\bq_h\cdot\bm n_1,E_\Gamma\eta_{\varphi_2}\rangle_{\partial\Omega_1}|\le\|\bq_h\|_{\bQ}\|E_\Gamma\eta_{\varphi_2}\|_{1/2,\partial\Omega_1}\le C_{\rm tr}\|\varphi_2-\mathcal I_h\varphi_2\|_{1,\Omega_2}\normiii{\vec{\by}_h}.
\]
In contrast with the second inequality in \eqref{eq:trace}, the $H^1(\Omega_2)$-norm of the interpolation error is not bounded   by its  norm $\|\bullet\|_{W,\varepsilon_r}$, which is what avoids the factor $\max\{1,\varepsilon_r^{-1/2}\}$. Third, by the definition \eqref{eq:Thetadef}, $|b_\Gamma(\bm\eta_{\bp},\psi_{2,h})|\le\Theta_h(\bp)\|\psi_{2,h}\|_{W,\varepsilon_r}\le\Theta_h(\bp)\normiii{\vec{\by}_h}$. Collecting these three bounds in \eqref{eq:Meta} we arrive at
\[
\sup_{\vec{\by}_h\in\vec{\bY}_h\setminus\{\bm0\}}\frac{|\mathcal M(\vec{\boldsymbol\eta};\vec{\by}_h)|}{\normiii{\vec{\by}_h}}\le\normiii{\vec{\boldsymbol\eta}}+C_{\rm tr}\|\varphi_2-\mathcal I_h\varphi_2\|_{1,\Omega_2}+\Theta_h(\bp).
\]
As in the proof of (a), the discrete inf-sup condition \eqref{idinfsup} and the Galerkin orthogonality give
\begin{align*}
\gamma_\sharp\normiii{\mathcal P_h\vec{\bx}-\vec{\bx}_h}&\le\sup_{\vec{\by}_h\in\vec{\bY}_h\setminus\{\bm0\}}\frac{\mathcal M(\mathcal P_h\vec{\bx}-\vec{\bx};\vec{\by}_h)}{\normiii{\vec{\by}_h}} =\sup_{\vec{\by}_h\in\vec{\bY}_h\setminus\{\bm0\}}\frac{|\mathcal M(\vec{\boldsymbol\eta};\vec{\by}_h)|}{\normiii{\vec{\by}_h}},
\end{align*}
where we used that $\mathcal M(\mathcal P_h\vec{\bx}-\vec{\bx};\vec{\by}_h)=-\mathcal M(\vec{\boldsymbol\eta};\vec{\by}_h)$ and that $\vec{\bY}_h$ is invariant under $\vec{\by}_h\mapsto-\vec{\by}_h$. The last two displays and the triangle inequality $\normiii{\vec{\bx}-\vec{\bx}_h}\le\normiii{\vec{\boldsymbol\eta}}+\normiii{\mathcal P_h\vec{\bx}-\vec{\bx}_h}$ yield \eqref{error2}.
\end{myproof}

The quantity $\Theta_h(\bp)$ is the only term in \eqref{error2} that is not a standard interpolation error. It is estimated in the next lemma, which shows that, uniformly in $\varepsilon_r$, it loses at most half an order with respect to the interpolation error of the normal flux.

\begin{lemma}[Interface consistency]\label{lem:interface}
Let $\bp\in\bH^1(\Omega_1)$ be such that $(\bp\cdot\bm n_1)|_F\in H^{k_1+1}(F)$ for all $F\in\mathcal F_h^\Gamma$. Then there exists $C_I>0$, depending only on $d$, $k_1$, $k_2$ and the shape regularity of $\mathcal T_h$, such that
\begin{equation}\label{eq:Theta}
\Theta_h(\bp)\le C_I\,h_\Gamma^{k_1+1/2}\min\big\{1,h_\Gamma\,\varepsilon_r^{-1/2}\big\}\,|\bp\cdot\bm n_1|_{k_1+1,\Gamma_h}.
\end{equation}
\end{lemma}
\begin{proof}
Let $\psi_h\in W_h^0\setminus\{0\}$ and write $\bm\eta_{\bp}:=\bp-\bPi_h\bp$. Since $\bp\in\bH^1(\Omega_1)$, the trace theorem gives $\bp\cdot\bm n_1\in L^2(\partial\Omega_1)$, and the normal trace of $\bPi_h\bp$ is piecewise polynomial, so that $\bm\eta_{\bp}\cdot\bm n_1\in L^2(\partial\Omega_1)$. Then the discussion following the definition of $b_\Gamma$ in \Cref{sec:model} shows that
\[
b_\Gamma(\bm\eta_{\bp},\psi_h)=\langle\bm\eta_{\bp}\cdot\bm n_1,\psi_h\rangle_\Gamma=\sum_{F\in\mathcal F_h^\Gamma}\langle\bm\eta_{\bp}\cdot\bm n_1,\psi_h\rangle_F.
\]
By \eqref{eq:RTtrace}, $\bm\eta_{\bp}\cdot\bm n_1|_F=(I-Q_F)(\bp\cdot\bm n_1|_F)$ is $L^2(F)$-orthogonal to $\mathbb P_{k_1}(F)$. Hence, for each $F\in\mathcal F_h^\Gamma$ and any $\mu_F\in\mathbb P_{k_1}(F)$, the Cauchy--Schwarz inequality, \eqref{eq:facetapprox} and $h_F\le h_{T_F}$ imply
\begin{equation}\label{eq:facet}
\big|\langle\bm\eta_{\bp}\cdot\bm n_1,\psi_h\rangle_F\big|=\Big|\langle(I-Q_F)(\bp\cdot\bm n_1),(\psi_h-\mu_F)\rangle_F\Big|\le C_A\,h_{T_F}^{k_1+1}|\bp\cdot\bm n_1|_{k_1+1,F}\,\|\psi_h-\mu_F\|_{0,F}.
\end{equation}
We now bound the last factor in two different ways, {employing} the discrete trace inequality \cite[Lemma 12.8]{Ern1}: there exists $C_{\rm dt}>0$, depending only on $d$, $k_2$ and the shape regularity of $\mathcal T_h$, such that $\|w\|_{0,F}\le C_{\rm dt}h_T^{-1/2}\|w\|_{0,T}$ for every $T\in\mathcal T_h$, every facet $F$ of $T$, and every $w\in\mathbb P_{k_2}(T)$. Choosing $\mu_F=0$ and applying this inequality with $T=T_F$ and $w=\psi_h|_{T_F}$, we obtain
\[
\|\psi_h\|_{0,F}\le C_{\rm dt}h_{T_F}^{-1/2}\|\psi_h\|_{0,T_F}.
\]
Alternatively, choosing $\mu_F:=\overline{\psi}_{T_F}:=|T_F|^{-1}\int_{T_F}\psi_h$ (a constant, hence an element of $\mathbb P_{k_1}(F)$), applying the discrete trace inequality to $w=\psi_h|_{T_F}-\overline{\psi}_{T_F}$, and then the Poincar\'e inequality $\|w-\overline w_T\|_{0,T}\le\pi^{-1}h_T\|\nabla w\|_{0,T}$, valid on the convex set $T$ \cite{PayneWeinberger1960,Bebendorf2003}, we obtain
\[
\|\psi_h-\overline{\psi}_{T_F}\|_{0,F}\le C_{\rm dt}h_{T_F}^{-1/2}\|\psi_h-\overline{\psi}_{T_F}\|_{0,T_F}\le C_{\rm dt}\,\pi^{-1}h_{T_F}^{1/2}\|\nabla\psi_h\|_{0,T_F}.
\]
We insert each of these two bounds in \eqref{eq:facet}, sum over $F\in\mathcal F_h^\Gamma$, employ$h_{T_F}\le h_\Gamma$ and the Cauchy--Schwarz inequality in $\mathbb R^{\#\mathcal F_h^\Gamma}$, and take into account that a given element $T\subset\overline\Omega_2$ coincides with $T_F$ for at most $d+1$ facets $F\in\mathcal F_h^\Gamma$. This leads to
\begin{align*}
|b_\Gamma(\bm\eta_{\bp},\psi_h)|&\le C_AC_{\rm dt}(d+1)^{1/2}\,h_\Gamma^{k_1+1/2}\,|\bp\cdot\bm n_1|_{k_1+1,\Gamma_h}\|\psi_h\|_{0,\Omega_2},\\
|b_\Gamma(\bm\eta_{\bp},\psi_h)|&\le C_AC_{\rm dt}(d+1)^{1/2}\pi^{-1}\,h_\Gamma^{k_1+3/2}\,|\bp\cdot\bm n_1|_{k_1+1,\Gamma_h}\|\nabla\psi_h\|_{0,\Omega_2}.
\end{align*}
Finally, the definition of $\|\bullet\|_{W,\varepsilon_r}$ gives $\|\psi_h\|_{0,\Omega_2}\le\|\psi_h\|_{W,\varepsilon_r}$ and $\|\nabla\psi_h\|_{0,\Omega_2}\le\varepsilon_r^{-1/2}\|\psi_h\|_{W,\varepsilon_r}$. Taking the minimum of the two resulting bounds, dividing by $\|\psi_h\|_{W,\varepsilon_r}$ and taking the supremum over $\psi_h$, we obtain \eqref{eq:Theta} with $C_I:=C_AC_{\rm dt}(d+1)^{1/2}$ (recall that $\pi^{-1}<1$).
\end{proof}

Combining \Cref{iestimatelemma}, \Cref{errthm}(b), and \Cref{lem:interface}, we obtain the following rates of convergence, in the spirit of \Cref{errthm}(a) but with constants that do not depend on the parameters.

\begin{corollary}[Parameter-robust rates of convergence]\label{cor:rates}
In addition to the assumptions of \Cref{iestimatelemma}, suppose that $(\bp\cdot\bm n_1)|_F\in H^{k_1+1}(F)$ for all $F\in\mathcal F_h^\Gamma$ (this holds, for instance, if $\bp\in\bH^{k_1+2}(\Omega_1)$, since $\bm n_1$ is constant on each $F$, and then $|\bp\cdot\bm n_1|_{k_1+1,\Gamma_h}\lesssim\|\bp\|_{k_1+2,\Omega_1}$ by the trace theorem). Then there exists $C>0$, independent of $h$, $\delta_c$ and $\varepsilon_r$, such that
\begin{align}
\normiii{\vec{\bx}-\vec{\bx}_h}\lesssim 
\,&h^{k_1+1}\Big(|\bp|_{k_1+1,\Omega_1}^{}+|\vdiv\bp|_{k_1+1,\Omega_1}^{}+\delta_c^{-1}|\bzeta|_{k_1+1,\Omega_1}^{}\nonumber\\
&\qquad\qquad+|\vdiv\bzeta|_{k_1+1,\Omega_1}^{}+|\varphi_1|_{k_1+1,\Omega_1}^{}+|u|_{k_1+1,\Omega_1}^{}\Big)\nonumber\\
&+h^{k_2}\big(1+\varepsilon_r^{1/2}\big)|\varphi_2|_{k_2+1,\Omega_2}^{}+h^{k_1+1/2}\min\big\{1,h\,\varepsilon_r^{-1/2}\big\}|\bp\cdot\bm n_1|_{k_1+1,\Gamma_h}^{}.\label{eq:rates}
\end{align}
In particular, if $k_2=k_1+1$ and $\varepsilon_r\le1$, then $\normiii{\vec{\bx}-\vec{\bx}_h}=\mathcal O(h^{k_1+1/2})$ with a constant independent of $\delta_c$ and $\varepsilon_r$ (so that the rate is uniform provided the norms of the solution in \eqref{eq:rates} are bounded), and $\normiii{\vec{\bx}-\vec{\bx}_h}=\mathcal O(h^{k_1+1})$, again with a constant independent of $\delta_c$ and $\varepsilon_r$, whenever $h\le\varepsilon_r^{1/2}$.
\end{corollary}
\begin{proof}
We bound the three terms on the right-hand side of \eqref{error2}. The first {term} is estimated by \eqref{ierrestimate}, and $h^{k_2}(\varepsilon_r^{1/2}+h)\le(1+{\rm diam}(\Omega))\,h^{k_2}(1+\varepsilon_r^{1/2})$. For the second {term}, the local Lagrange interpolation \eqref{eq:H1approx} 
yields $\|\varphi_2-\mathcal I_h\varphi_2\|_{1,\Omega_2}\lesssim(1+h)h^{k_2}|\varphi_2|_{k_2+1,\Omega_2}$. The third {term} 
is estimated by \Cref{lem:interface}, using $h_\Gamma\le h$ and that $t\mapsto t^{k_1+1/2}\min\{1,t\,\varepsilon_r^{-1/2}\}$ is nondecreasing. Since $\gamma_\sharp$, $C_{\rm tr}$, $C_A$, and $C_I$ do not depend on $h$, $\delta_c$ or $\varepsilon_r$, this proves \eqref{eq:rates}. If $k_2=k_1+1$, the first two terms in \eqref{eq:rates} are $\mathcal O(h^{k_1+1})$, while the third one is bounded by $h^{k_1+1/2}$ and, for $h\le\varepsilon_r^{1/2}$, by $h^{k_1+3/2}\varepsilon_r^{-1/2}\le h^{k_1+1}$.
\end{proof}

\begin{remark}\label{rem:errrobust}
\Cref{errthm}(a) is robust in the sense that the constant in \eqref{error} does not depend on $\delta_c$, and depends on $\varepsilon_r$ only through $C_\Gamma$. Part (b) and \Cref{cor:rates} remove the latter dependence. The factor $\varepsilon_r^{-1/2}$ in $C_{\mathcal M}$ arises from bounding $b_\Gamma$ in the weighted norm of $W$ (cf. \eqref{eq:trace}), whereas in \eqref{eq:Meta} the two interface terms only involve interpolation errors, which can be estimated directly. The price to pay is a loss of at most half an order in the pre-asymptotic regime $h>\varepsilon_r^{1/2}$, in which the mesh does not resolve the length scale $\varepsilon_r^{1/2}$ of \eqref{IP}$_2$. Heuristically, the error in $\varphi_2$ is driven by the error in the normal flux $\bp\cdot\bm n_1$ on $\Gamma$ (cf. the last equation of \eqref{IPh}), and for $h>\varepsilon_r^{1/2}$ the discrete potential absorbs this flux error within the layer of elements of $\Omega_2$ adjacent to $\Gamma$, which costs the factor $h^{-1/2}$ of the discrete trace inequality. Since parts (a) and (b) of \Cref{errthm} hold simultaneously, the error is bounded by the minimum of the two corresponding bounds; for $k_2=k_1+1$ these are of order $(1+C_\Gamma)h^{k_1+1}$ and $h^{k_1+1}+h^{k_1+1/2}\min\{1,h\varepsilon_r^{-1/2}\}$, respectively. Note that for $\varepsilon_r\le1$ the bound in (b) is never worse than that in (a), up to constants, so that (b) implies (a) in that regime.  The half-order reduction is visible in Example 1 of \Cref{sec:numerical} for $k_1=1$ and $\varepsilon_r=10^{-6}$. However, the Sobolev norms of the exact solution in \eqref{ierrestimate} do depend on the parameters: for small $\delta_c$, $u$, and $\bzeta$ develop boundary layers of width $O(\delta_c)$ at $\partial\Omega_1$ (since the condition on $u$ is lost in the limit $\delta_c\to0$), and for small $\varepsilon_r$, $\varphi_2$ develops a layer of width $O(\varepsilon_r^{1/2})$ at $\Gamma$. Uniform convergence rates would require layer-adapted meshes.
\end{remark}

The Maxwell stress tensor, denoted by $\btau$, in the electrolyte occupying the domain $\Omega_{1}$ with dielectric constant $\varepsilon_{1}$, is defined in physical units as
\begin{equation}\label{stress}
    \btau^{\rm phys} = \varepsilon_{1}^{} \bp\bp^T - \frac{\varepsilon_1}{2}\bp^T\bp\mathbb{I} +\varepsilon_{1}^{} l_c^2\bigg[ \nabla(\vdiv\bp)^T\bp\mathbb{I} - \bp\nabla(\vdiv\bp)^T - \nabla(\vdiv\bp)\bp^T + \frac{1}{2}(\vdiv\bp)^2\mathbb{I}\bigg],
\end{equation}
where in this formula only $\bp=\nabla\varphi_{1}^{}$ represents the electric field in physical units (with dimensional potential and position), $l_c$ is the correlation length (a length), and $\mathbb{I}$ denotes the identity matrix. Problem \eqref{IP}, however, is dimensionless  (see \eqref{PP}): lengths are scaled by the Debye length $\lambda_D$ and potentials by the thermal voltage $k_{B}T/e$, so that $l_c=\delta_c\lambda_D$. Writing \eqref{stress} in terms of the dimensionless variables (still denoted by $\bp$, $\nabla$, $\vdiv$), the physical field is $\frac{k_{B}T}{e\lambda_D}\bp$, every derivative contributes a factor $\lambda_D^{-1}$, and, since $u=-\vdiv\bp$ and $\bzeta=\delta_c^2\nabla u$, the correlation terms transform as
\[
l_c^2\nabla(\vdiv\bp) \mapsto \delta_c^2\lambda_D^2\,\frac{k_{B}T}{e\lambda_D^{3}}\nabla(\vdiv\bp)=-\frac{k_{B}T}{e\lambda_D}\,\bzeta,\qquad
l_c^2(\vdiv\bp)^2\ \mapsto\ \Big(\frac{k_{B}T}{e\lambda_D}\Big)^2\delta_c^2u^2 .
\]
Therefore, $\btau^{\rm phys}=\varepsilon_1\big(\frac{k_{B}T}{e\lambda_D}\big)^2\btau=\varepsilon_1\big(\frac{k_{B}T}{e}\big)^2\frac{\delta_c^2}{l_c^2}\btau$, where the dimensionless  $\btau$ is given in terms of the mixed variables by
\begin{equation}\label{mstress}
    \btau = \bp\bp^T - \frac12\bp^T\bp\mathbb{I} -\bzeta^T\bp\mathbb{I} + \bp\bzeta^T + \bzeta\bp^T + \frac{\delta_c^2}{2}u^2\mathbb{I}.
\end{equation}
Since the factor $\varepsilon_1(k_{B}T/(e\lambda_D))^2$ is a constant, in the analysis and in all the numerical examples we work with the dimensionless tensor \eqref{mstress}; the physical parameters $\varepsilon_1$ and $l_c$ enter only through this factor and through $\delta_c=l_c/\lambda_D$. 

The Maxwell stress tensor $\btau$ (which belongs to $\mathbb{L}^{1}(\Omega_{1})$ in general, and to $\mathbb{L}^{2}(\Omega_{1})$ under the assumptions of \Cref{stressthm} below) can be approximated by 
\begin{equation}\label{discretestress}
        \btau_{h}^{}= \bp_{h}^{}\bp_{h}^{T} - \frac12\bp_{h}^{T}\bp_{h}^{}\mathbb{I} -\bzeta_{h}^{T}\bp_{h}^{}\mathbb{I} + \bp_{h}^{}\bzeta_{h}^{T} + \bzeta_{h}^{}\bp_{h}^{T} + \frac{\delta_c^2}{2}u_{h}^{2}\mathbb{I} \in [\mathbb{P}_{2k_1+2}(\mathcal{T}_h)]^{d\times d},
\end{equation}
where $(\bp_{h}^{},\bzeta_{h}^{},\varphi_{1,h}^{},u_{h}^{},\varphi_{2,h}^{})$ is the solution of \eqref{IPh}.

The following corollary demonstrates that the post-processed approximation of the Maxwell stress tensor converges  in the $\mathbb{L}^2$-norm, with optimal order except in the lowest-order case $k_1=0$ in 3D.

\begin{corollary}[Convergence of the Maxwell stress]\label{stressthm}
    Under the regularity assumptions of \Cref{iestimatelemma}, assume that $\bp,\bzeta\in\bL^\infty(\Omega_1)$, $u\in L^\infty(\Omega_1)$, and that the mesh family is quasi-uniform, and let $m:=\min\{k_1+1,k_2\}$. Then
    \begin{equation}\label{stresserror}
            \|\btau-\btau_{h}^{}\|_{0,\Omega_1} \lesssim h^{m}+h^{2m-d/2},
    \end{equation}
    where the hidden constant depends on $\delta_c$ (through the factor $1+\delta_c+\delta_c^2$), on $C_\Gamma$ and on the norms of the exact solution. In particular, \eqref{stresserror} is of optimal order $\mathcal O(h^m)$ if $m\ge d/2$, which, for $k_2=k_1+1$, holds in all cases except $k_1=0$ in 3D, where it gives $\mathcal O(h^{1/2})$. If $\varepsilon_r\le1$, using \Cref{errthm}(b) and \Cref{cor:rates} instead of \Cref{errthm}(a) removes the dependence on $C_\Gamma$, at the price of replacing $h^m$ by $h^{k_1+1/2}\min\{1,h\,\varepsilon_r^{-1/2}\}+h^{m}$ in \eqref{stresserror}. The same estimate holds for $\btau^{\rm phys}$, multiplied by $\varepsilon_1(k_{B}T/(e\lambda_D))^2$.
\end{corollary}
\begin{proof}
    The tensor \eqref{mstress} is a linear combination of products $fg$ with $f,g\in\{\bp,\bzeta,u,\mathbb I\}$, and \eqref{discretestress} is the same combination of the products $f_hg_h$ with $f_h,g_h\in\{\bp_h,\bzeta_h,u_h,\mathbb I\}$. Writing $fg-f_hg_h=f(g-g_h)+(f-f_h)g-(f-f_h)(g-g_h)$ and using the H\"older inequality with exponents $(\infty,2)$ in the first two terms and $(4,4)$ in the last one, we get
    \begin{equation}\label{eq:quad_bound}
        \|fg - f_hg_h\|_{0,\Omega_1} \le \|f\|_{0,\infty,\Omega_1} \|g-g_h\|_{0,\Omega_1} + \|g\|_{0,\infty,\Omega_1} \|f-f_h\|_{0,\Omega_1}+\|f-f_h\|_{0,4,\Omega_1}\|g-g_h\|_{0,4,\Omega_1},
    \end{equation}
    so that only the exact solution needs to be bounded in $L^\infty(\Omega_1)$. Let $\bm e_{\bp}:=\bp-\bp_h$, $\bm e_{\bzeta}:=\bzeta-\bzeta_h$, $e_u:=u-u_h$, and let $M_\infty$ be the sum of the $L^\infty(\Omega_1)$-norms of $\bp$, $\bzeta$ and $u$. Applying \eqref{eq:quad_bound} to each of the six terms of $\btau-\btau_h$, and using $|\bm a\bm b^T|=|\bm a||\bm b|$ and $|\mathbb{I}|=\sqrt d$ for the Frobenius norm, we obtain
    \begin{equation}\label{eq:stresssplit}
    \|\btau-\btau_h\|_{0,\Omega_1}\lesssim M_\infty\big(\|\bm e_{\bp}\|_{0,\Omega_1}+\|\bm e_{\bzeta}\|_{0,\Omega_1}+\delta_c^2\|e_u\|_{0,\Omega_1}\big)+\big(\|\bm e_{\bp}\|_{0,4,\Omega_1}+\|\bm e_{\bzeta}\|_{0,4,\Omega_1}\big)^2+\delta_c^2\|e_u\|_{0,4,\Omega_1}^2.
    \end{equation}
    For the first group of terms, $\|\bm e_{\bp}\|_{0,\Omega_1},\|e_u\|_{0,\Omega_1}\le\normiii{\vec{\bx}-\vec{\bx}_h}$ and $\|\bm e_{\bzeta}\|_{0,\Omega_1}\le\delta_c\normiii{\vec{\bx}-\vec{\bx}_h}$, and $\normiii{\vec{\bx}-\vec{\bx}_h}\lesssim h^m$ by \Cref{errthm}(a) and \Cref{iestimatelemma} (for $\varepsilon_r\lesssim1$). It remains to bound the $L^4(\Omega_1)$-norms of the errors in the last group of terms of \eqref{eq:stresssplit}; this requires no additional regularity, since only the energy-norm error estimates and an inverse inequality are used. We split, e.g., $\bm e_{\bp}=(\bp-\bPi_h\bp)+(\bPi_h\bp-\bp_h)$. For the continuous part, we employ \eqref{eq:L4approx}. 
    For the discrete part, the local inverse inequality $\|\bq_h\|_{0,4,T}\lesssim h_T^{-d/4}\|\bq_h\|_{0,T}$ for polynomials \cite[Lemma 12.1]{Ern1} and the quasi-uniformity of the mesh give $\|\bPi_h\bp-\bp_h\|_{0,4,\Omega_1}\lesssim h^{-d/4}\|\bPi_h\bp-\bp_h\|_{0,\Omega_1}\le h^{-d/4}\normiii{\mathcal P_h\vec{\bx}-\vec{\bx}_h}$. Summing the fourth powers of the local bounds and using $\sum_Ta_T^4\le(\sum_Ta_T^2)^2$, and recalling from the proof of \Cref{errthm} that $\normiii{\mathcal P_h\vec{\bx}-\vec{\bx}_h}\le\gamma_\sharp^{-1}C_{\mathcal M}\normiii{\vec{\bx}-\mathcal P_h\vec{\bx}}\lesssim h^m$, we arrive at $\|\bm e_{\bp}\|_{0,4,\Omega_1}\lesssim h^{m-d/4}$. In the same way (with $P_h$ in place of $\bPi_h$ for $u$, and with the factor $\delta_c$ for $\bzeta$, since $\|\bPi_h\bzeta-\bzeta_h\|_{0,\Omega_1}\le\delta_c\normiii{\mathcal P_h\vec{\bx}-\vec{\bx}_h}$), $\|\bm e_{\bzeta}\|_{0,4,\Omega_1}\lesssim\delta_c h^{m-d/4}$ and $\|e_u\|_{0,4,\Omega_1}\lesssim h^{m-d/4}$. Inserting these bounds in \eqref{eq:stresssplit} gives \eqref{stresserror}. The robust variant follows in the same way, bounding $\normiii{\vec{\bx}-\vec{\bx}_h}$ and $\normiii{\mathcal P_h\vec{\bx}-\vec{\bx}_h}$ through \eqref{error2} and \Cref{cor:rates}.
\end{proof}

The assumption $\bp,\bzeta,u\in L^\infty(\Omega_1)$ is only used for the two terms of \eqref{eq:quad_bound} that are linear in the error, and it follows from the regularity of \Cref{iestimatelemma} whenever $k_1+1>d/2$. Without it, using the H\"older inequality with exponents $(4,4)$ also in these two terms, together with $\bp,\bzeta,u\in H^1(\Omega_1)\subset L^4(\Omega_1)$ ($d\le3$) and the above $L^4$-bounds of the errors, the same argument gives $\|\btau-\btau_h\|_{0,\Omega_1}\lesssim h^{m-d/4}+h^{2m-d/2}$, that is, a convergent but suboptimal estimate.

\section{Numerical results}\label{sec:numerical}

In this section, we present a series of numerical experiments designed to validate the theoretical convergence estimates derived in \Cref{errthm} and \Cref{cor:rates} and to demonstrate the physical applicability of the proposed \gls{fem}. To ensure balanced, optimal convergence rates across all mixed and primal variables, we fix the polynomial degrees as $k_{1}=k$ and $k_{2}=k+1$ for $k \geq 0$ throughout our tests.

Our numerical investigations are fourfold. First, we validate the parameter robust convergence of the error with respect to $\delta_c$ and $\varepsilon_r$ in \Cref{subsec:example01}. Second, in \Cref{subsec:example02}-\Cref{subsec:example04}, we conduct convergence studies with several manufactured exact smooth solutions in different geometric configurations, validating the optimal convergence of the field variables as well as the Maxwell electrostatic stress tensor. Third, we utilise a closed-form solution for a single charged particle to confirm the convergence rates and their robustness with respect to the correlation length for realistic permittivity ratios, and to illustrate overscreening in \Cref{subsec:example05}. Fourth, in \Cref{subsec:example06}, we compute the interaction force between two charged particles to illustrate the charge reversal phenomenon in electrolytes with multivalent ions.
All numerical implementations are performed using the open-source finite element library \texttt{Gridap} \cite{Verdugo2022}.

\subsection{Example 1: Two-dimensional convergence tests with a straight interface and parameter robustness}\label{subsec:example01}

For our first example, we consider a square domain $\Omega = (-1,1)^2$ partitioned by a vertical interface at $x=0$, as illustrated in \Cref{fig:domainB}, that is, $\Omega_1=(-1,0)\times(-1,1)$, $\Omega_2=(0,1)\times(-1,1)$ and $\Gamma=\{0\}\times(-1,1)$. We prescribe the exact electrostatic potentials
$\varphi_1=\sin(x+y)$, $\varphi_2=\exp(x+y)$, 
which are discontinuous across $\Gamma$, with boundary conditions and volumetric source terms derived accordingly. In particular, the jump $j=\varphi_1-\varphi_2$ and the flux datum $g=\nabla\varphi_1\cdot\bm n_1+\varepsilon_r\nabla\varphi_2\cdot\bm n_2$ do not vanish on $\Gamma$, so that both transmission conditions are tested, and $g_3=u=-\Delta\varphi_1=2\sin(x+y)$ on $\partial\Omega_1$. Consequently, $\bp$, $u$, $\varphi_1$ and $\varphi_2$ do not depend on the parameters, whereas $\bzeta=\delta_c^2\nabla u$. The exact solution has no layers, and all the Sobolev norms in \eqref{eq:rates} are bounded uniformly for $\delta_c,\varepsilon_r\in(0,1]$ (note that $\delta_c^{-1}|\bzeta|_{k+1,\Omega_1}=\delta_c|\nabla u|_{k+1,\Omega_1}$). This allows us to observe the dependence of the error on the parameters through the stability and consistency constants only. We consider all $15$ combinations of $\delta_c\in\{10,2,1,10^{-3},10^{-6}\}$ and $\varepsilon_r\in\{1,10^{-3},10^{-6}\}$, the polynomial degrees $k\in\{0,1\}$, and a sequence of five unstructured meshes with $h$ ranging from $0.532$ to $0.033$.

For $\delta_c\ll1$, the term $a_2(\bzeta_h,\bxi_h)=\delta_c^{-2}\int_{\Omega_1}\bzeta_h\cdot\bxi_h$ produces matrix entries of order $\delta_c^{-2}$ next to entries of order one, while the unknown $\bzeta_h$ is of order $\delta_c^2$. With a sparse direct solver this leads to a severe loss of accuracy: in preliminary computations with $k\in\{0,1,2\}$, the relative residual of the linear system was of order $10^{-7}$--$10^{-6}$ for $\delta_c=10^{-3}$, which already spoils the convergence of $\bzeta_h$ for $k=2$. It reaches $0.5$ for $\delta_c=10^{-6}$, in which case the relative error of $\bzeta_h$ was larger than one for all $k$. We therefore solve for $\widehat{\bzeta}_h:=\delta_c^{-2}\bzeta_h$, that is, we replace $a_2(\bzeta_h,\bxi_h)$ by $\int_{\Omega_1}\widehat{\bzeta}_h\cdot\bxi_h$ and $b(\bzeta_h,v_h)$ by $\delta_c^2\,b(\widehat{\bzeta}_h,v_h)$ in \eqref{IPh}, the right-hand side being unchanged. This change of variable, which is the discrete counterpart of $\delta_c^{-2}\bzeta=\nabla u$, does not modify the discrete solution in exact arithmetic, and with it the relative residual is below $2\cdot10^{-11}$ in all the computations reported below.

\Cref{table:eg1} reports the total error $\mathsf e_h:=\normiii{\vec{\bx}-\vec{\bx}_h}$, measured in the parameter-weighted norm of \Cref{sec:wellposedness}, together with the experimental convergence rates. The errors decay with the optimal order $\mathcal O(h^{k+1})$ in all cases, and they are essentially independent of $\delta_c$: for $\delta_c=10^{-3}$ and $\delta_c=10^{-6}$ they coincide in the three digits shown (the larger values for $\delta_c=1$ stem from $\|\vdiv(\bzeta-\bzeta_h)\|_{0,\Omega_1}$, since $\bzeta=\nabla u$ is not small in that case). Decreasing $\varepsilon_r$ does not increase the errors either. The only noticeable effect is a slow decrease of the rate for $k=1$ and $\varepsilon_r=10^{-6}$, from $1.96$ to $1.81$. This is consistent with \Cref{cor:rates}: since $h>\varepsilon_r^{1/2}=10^{-3}$ on all the meshes, the interface consistency term of \Cref{lem:interface} is only bounded by $\mathcal O(h^{3/2})$, and it becomes visible as soon as it is no longer dominated by the remaining contributions, which are $\mathcal O(h^{2})$. {For $\delta_c>1$ the errors grow because $|\vdiv\bzeta|_{k+1,\Omega_{1}}^{}=\delta_c^2|\Delta u|_{k+1,\Omega_{1}}^{}$ enters \eqref{ierrestimate}.}
The individual components of the error, which we do not report in detail, support this interpretation: for $k=1$ and $\varepsilon_r=10^{-6}$, the experimental rates on the last three meshes are $1.81$, $1.53$ and $1.46$ for $\|\varphi_2-\varphi_{2,h}\|_{W,\varepsilon_r}$, and $1.89$, $1.82$ and $1.74$ for $\|\bp-\bp_h\|_{\bQ}$, whereas $\|\varphi_1-\varphi_{1,h}\|_{0,\Omega_1}$ and $\|u-u_h\|_{0,\Omega_1}$ converge with order $2.00$. For $\varepsilon_r=10^{-3}$, $\varepsilon_r^{1/2}\approx0.032$ is close to the size of the finest mesh, and accordingly the rate of the $\varphi_2$-component, after decreasing to $1.31$, recovers to $1.81$ on the last mesh, while the total error converges with order two. For $k=0$ no reduction is observed on these meshes. Finally, the relative errors of all the individual variables, including $\bzeta$, converge with the optimal order for all values of $\delta_c$, with the exceptions just described for $k=1$ and $\varepsilon_r\le10^{-3}$; for $\varepsilon_r=10^{-6}$ and $\delta_c\ll1$, the reduction also affects the (very small) component $\bzeta$, whose relative error converges with order $1.82$ ($\delta_c=10^{-3}$) and $1.79$ ($\delta_c=10^{-6}$) on the finest mesh.

\begin{table}[t!]
\centering
\resizebox{\textwidth}{!}{\begin{tabular}{ccccccccccccc}
\toprule
 & & & \multicolumn{2}{c}{$\delta_c=10$} & \multicolumn{2}{c}{$\delta_c=2$} & \multicolumn{2}{c}{$\delta_c=1$} & \multicolumn{2}{c}{$\delta_c=10^{-3}$} & \multicolumn{2}{c}{$\delta_c=10^{-6}$} \\
\cmidrule(lr){4-5}\cmidrule(lr){6-7}\cmidrule(lr){8-9}\cmidrule(lr){10-11}\cmidrule(lr){12-13}
$\varepsilon_r$ & order & $h$ & $\mathsf e_h$ & Rate & $\mathsf e_h$ & Rate & $\mathsf e_h$ & Rate & $\mathsf e_h$ & Rate & $\mathsf e_h$ & Rate \\
\midrule
 &  & 0.532 & 6.27e+01 & $\star$ &  2.77e+00 & $\star$ & 1.19e+00 & $\star$ & 9.52e-01 & $\star$ & 9.52e-01 & $\star$ \\
 &  & 0.266 & 3.14e+01 & 1.00    &  1.39e+00 & 1.00    & 5.96e-01 & 1.00 & 4.75e-01 & 1.00 & 4.75e-01 & 1.00 \\
$1$ & $k=0$ & 0.133 & 1.57e+01 & 1.00    &  6.93e-01 & 1.00    &  2.98e-01 & 1.00 & 2.38e-01 & 1.00 & 2.38e-01 & 1.00 \\
 &  & 0.067 & 7.85e+00 & 1.00    &  3.47e-01 & 1.00    & 1.49e-01 & 1.00 & 1.19e-01 & 1.00 & 1.19e-01 & 1.00 \\
 &  & 0.033 & 3.93e+00 & 1.00    &  1.73e-01 & 1.00    & 7.45e-02 & 1.00 & 5.94e-02 & 1.00 & 5.94e-02 & 1.00 \\
\cmidrule(lr){2-13}
 &  & 0.532 & 3.54e+00 & $\star$   &  1.65e-01 & $\star$ & 8.03e-02 & $\star$ & 6.72e-02 & $\star$ & 6.72e-02 & $\star$ \\
 &  & 0.266 & 8.96e-01 & 1.98      &  4.17e-02 & 1.98    & 2.03e-02 & 1.99 & 1.69e-02 & 1.99 & 1.69e-02 & 1.99 \\
 & $k=1$ & 0.133 & 2.25e-01 & 2.00      &  1.05e-02 & 2.00    & 5.08e-03 & 1.99 & 4.26e-03 & 1.99 & 4.26e-03 & 1.99 \\
 &  & 0.067 & 5.62e-02 & 2.00      &  2.62e-03 & 2.00    & 1.27e-03 & 2.00 & 1.07e-03 & 2.00 & 1.07e-03 & 2.00 \\
 &  & 0.033 & 1.41e-02 & 2.00      &  6.54e-04 & 2.00    & 3.19e-04 & 2.00 & 2.67e-04 & 2.00 & 2.67e-04 & 2.00 \\
\midrule
 &  & 0.532 & 6.27e+01 & $\star$   &  2.65e+00 & $\star$ & 8.78e-01 & $\star$ & 5.08e-01 & $\star$ & 5.08e-01 & $\star$ \\
 &  & 0.266 & 3.14e+01 & 1.00      &  1.33e+00 & 1.00    & 4.39e-01 & 1.00 & 2.52e-01 & 1.01 & 2.52e-01 & 1.01 \\
$10^{-3}$ & $k=0$ & 0.133 & 1.57e+01 & 1.00      &  6.63e-01 & 1.00    & 2.19e-01 & 1.00 & 1.26e-01 & 1.00 & 1.26e-01 & 1.00 \\
 &  & 0.067 & 7.85e+00 & 1.00      &  3.32e-01 & 1.00    & 1.10e-01 & 1.00 & 6.29e-02 & 1.00 & 6.28e-02 & 1.00 \\
 &  & 0.033 & 3.92e+00 & 1.00      &  1.66e-01 & 1.00    & 5.49e-02 & 1.00 & 3.14e-02 & 1.00 & 3.14e-02 & 1.00 \\
\cmidrule(lr){2-13}
 &  & 0.532 & 3.54e+00 & $\star$ & 1.54e-01 & $\star$ & 5.37e-02 & $\star$ & 3.08e-02 & $\star$ & 3.08e-02 & $\star$ \\
 &  & 0.266 & 8.96e-01 & 1.98 & 3.89e-02 & 1.98 & 1.36e-02 & 1.98 & 7.91e-03 & 1.96 & 7.91e-03 & 1.96 \\
 & $k=1$ & 0.133 & 2.25e-01 & 2.00 & 9.76e-03 & 1.99 & 3.45e-03 & 1.98 & 2.04e-03 & 1.95 & 2.04e-03 & 1.95 \\
 &  & 0.067 & 5.62e-02 & 2.00 & 2.45e-03 & 2.00 & 8.70e-04 & 1.99 & 5.22e-04 & 1.97 & 5.22e-04 & 1.97 \\
 &  & 0.033 & 1.41e-02 & 2.00 & 6.11e-04 & 2.00 & 2.16e-04 & 2.01 & 1.28e-04 & 2.02 & 1.28e-04 & 2.02 \\
\midrule
 &  & 0.532 & 6.27e+01 & $\star$ & 2.65e+00 & $\star$ & 8.78e-01 & $\star$ & 5.07e-01 & $\star$ & 5.07e-01 & $\star$ \\
 &  & 0.266 & 3.14e+01 & 1.00 & 1.33e+00 & 1.00 & 4.39e-01 & 1.00 & 2.51e-01 & 1.01 & 2.51e-01 & 1.01 \\
$10^{-6}$ & $k=0$ & 0.133 & 1.57e+01 & 1.00 & 6.63e-01 & 1.00 & 2.19e-01 & 1.00 & 1.26e-01 & 1.00 & 1.26e-01 & 1.00 \\
 &  & 0.067 & 7.85e+00 & 1.00 & 3.32e-01 & 1.00 & 1.10e-01 & 1.00 & 6.28e-02 & 1.00 & 6.28e-02 & 1.00 \\
 &  & 0.033 & 3.92e+00 & 1.00 & 1.66e-01 & 1.00 & 5.49e-02 & 1.00 & 3.14e-02 & 1.00 & 3.14e-02 & 1.00 \\
\cmidrule(lr){2-13}
 &  & 0.532 & 3.54e+00 & $\star$ & 1.54e-01 & $\star$ & 5.37e-02 & $\star$ & 3.08e-02 & $\star$ & 3.08e-02 & $\star$ \\
 &  & 0.266 & 8.96e-01 & 1.98 & 3.89e-02 & 1.98 & 1.36e-02 & 1.98 & 7.92e-03 & 1.96 & 7.92e-03 & 1.96 \\
 & $k=1$ & 0.133 & 2.25e-01 & 2.00 & 9.77e-03 & 1.99 & 3.47e-03 & 1.97 & 2.08e-03 & 1.93 & 2.08e-03 & 1.93 \\
 &  & 0.067 & 5.62e-02 & 2.00 & 2.46e-03 & 1.99 & 8.97e-04 & 1.95 & 5.65e-04 & 1.88 & 5.65e-04 & 1.88 \\
 &  & 0.033 & 1.41e-02 & 2.00 & 6.19e-04 & 1.99 & 2.38e-04 & 1.92 & 1.62e-04 & 1.81 & 1.62e-04 & 1.81 \\
\bottomrule
\end{tabular}}
\caption{Example 1: total error $\mathsf e_h=\normiii{\vec{\bx}-\vec{\bx}_h}$, in the parameter-weighted norm of \Cref{sec:wellposedness}, and experimental convergence rates, for $k\in\{0,1\}$ and all the combinations of $\delta_c\in\{10,2,1,10^{-3},10^{-6}\}$ and $\varepsilon_r\in\{1,10^{-3},10^{-6}\}$. Here and in the subsequent tables, the $\star$ denotes that no experimental convergence rate is computed at the first refinement level.}
\label{table:eg1}
\end{table}

\subsection{Example 2: Three-dimensional convergence tests with a planar interface}\label{subsec:example02}

\begin{figure}[t!]
    \centering
    \subfloat[Back-to-back domain   configuration.\label{fig:domain4}]{\begin{minipage}[b]{0.48\textwidth}
        \centering
        \raisebox{2mm}{\includegraphics{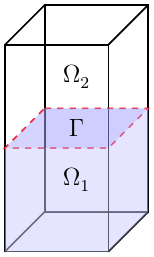}}
\end{minipage}}
\hfill
\subfloat[Sample of an unstructured coarse mesh.\label{fig:3dmesh}]{\begin{minipage}[b]{0.48\textwidth}
    \centering
    \includegraphics[width=0.47\linewidth]{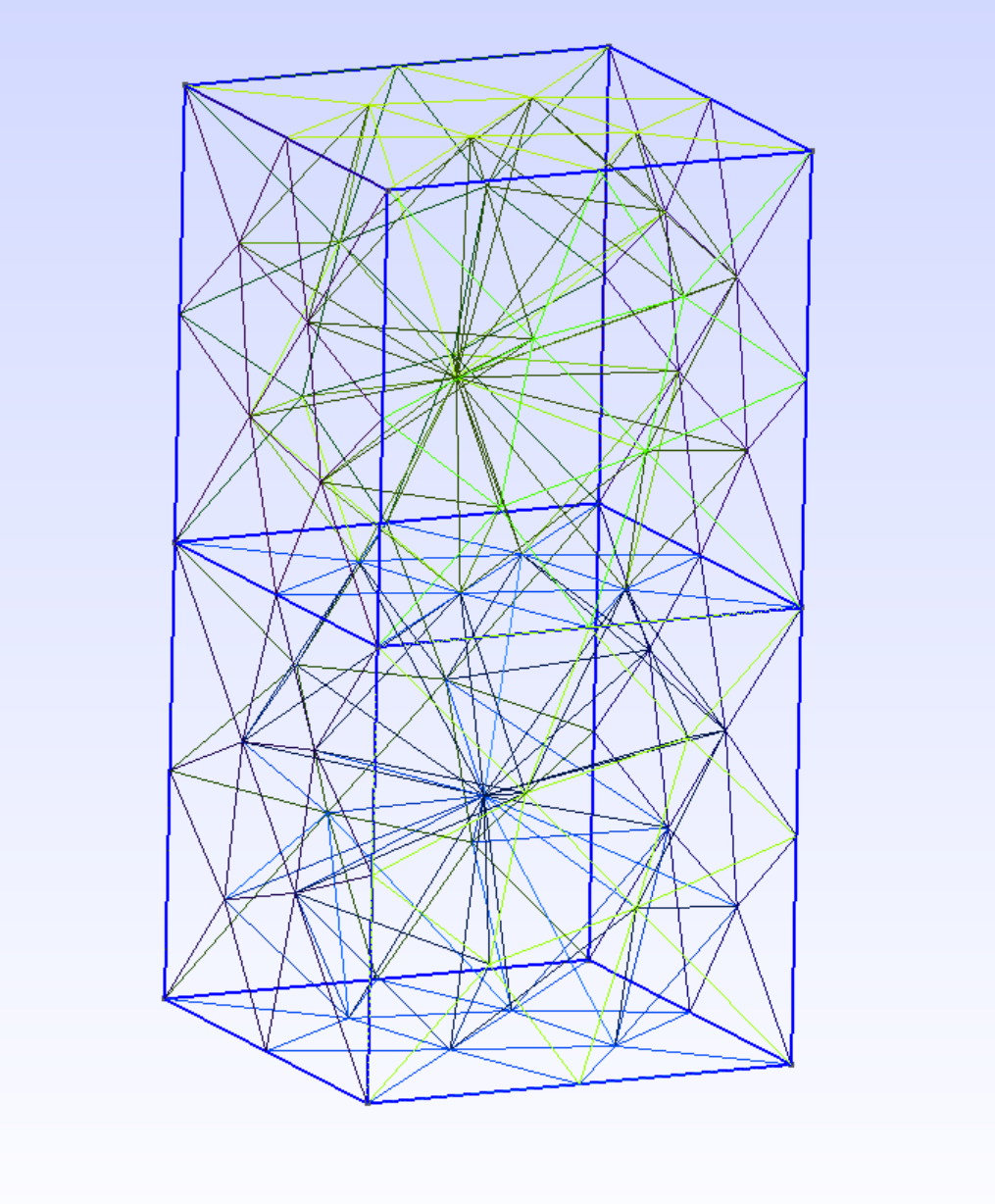}
\end{minipage}}
\caption{Three dimensional domain, interface configuration, and sample mesh used in Example 2.}
\end{figure}

To evaluate the performance of the proposed method in a three-dimensional setting, we consider the domain $\Omega = (0,1) \times (0,1) \times (0,2)$ partitioned by a horizontal interface at $z=1$, as depicted in \Cref{fig:domain4}. We prescribe the exact analytical solution as
$ \varphi_1=\sin(\pi x)\sin(\pi y)\sin(\pi z)$, $\varphi_2=\sin(\pi x)\sin(\pi y)\sin(\pi z)$, 
and select the model parameters $\varepsilon_r=2$ and $\delta_c=1$. A representative unstructured coarse mesh for this geometry is illustrated in \Cref{fig:3dmesh}. 

Consistent with the previous example, we discretise the domain using  $k \in \{0,1\}$. The computed error metrics and convergence rates are compiled in \Cref{table:eg2}. These numerical results verify that the mixed--primal \gls{fem} robustly maintains optimal $\mathcal{O}(h^{k+1})$ convergence rates across all field variables in three dimensions, confirming once more the theoretical results of \Cref{errthm}.

\begin{table}[t!]
\centering
\resizebox{\textwidth}{!}{%
\begin{tabular}{cccccccccccc}
\toprule
order & $h$ & $\|\bp-\bp_h\|_{\bQ}$ & Rate & $\|\bzeta-\bzeta_h\|_{\bQ,\delta_c}$ & Rate & $\|\varphi_1-\varphi_{1,h}\|_{0,\Omega_1}$ & Rate & $\|u-u_h\|_{0,\Omega_1}$ & Rate & $\|\varphi_2-\varphi_{2,h}\|_{W,\varepsilon_r}$ & Rate \\
\midrule
       & 0.492 & 5.05e+00 & $\star$ & 1.48e+02 & $\star$ & 1.71e-01 & $\star$ & 4.98e+00 & $\star$ & 1.71e+00 & $\star$ \\
$k=0$  & 0.246 & 2.67e+00 & 0.922 & 7.87e+01 & 0.916 & 8.96e-02 & 0.937 & 2.63e+00 & 0.922 & 9.91e-01 & 0.787 \\
       & 0.123 & 1.35e+00 & 0.986 & 3.98e+01 & 0.982 & 4.50e-02 & 0.993 & 1.33e+00 & 0.985 & 5.40e-01 & 0.876 \\
       & 0.062 & 6.75e-01 & 0.996 & 2.00e+01 & 0.996 & 2.25e-02 & 0.999 & 6.66e-01 & 0.996 & 2.80e-01 & 0.948 \\
\midrule
       & 0.492 & 1.10e+00 & $\star$ & 3.24e+01 & $\star$ & 3.74e-02 & $\star$ & 1.09e+00 & $\star$ & 4.22e-01 & $\star$ \\
$k=1$  & 0.246 & 3.11e-01 & 1.823 & 9.20e+00 & 1.815 & 1.04e-02 & 1.841 & 3.09e-01 & 1.824 & 1.48e-01 & 1.513 \\
       & 0.123 & 8.10e-02 & 1.943 & 2.40e+00 & 1.941 & 2.71e-03 & 1.945 & 8.03e-02 & 1.943 & 4.05e-02 & 1.867 \\
       & 0.062 & 2.04e-02 & 1.986 & 6.05e-01 & 1.985 & 6.84e-04 & 1.986 & 2.02e-02 & 1.986 & 1.06e-02 & 1.940 \\
\bottomrule
\end{tabular}}
\caption{Error history (individual errors and experimental convergence rates) for   Example 2.}
\label{table:eg2}
\end{table}

\subsection{Example 3: Two-dimensional convergence tests with disjoint immersed hexagonal inclusions}\label{subsec:example03}

\begin{figure}[t!]
    \centering
    \includegraphics{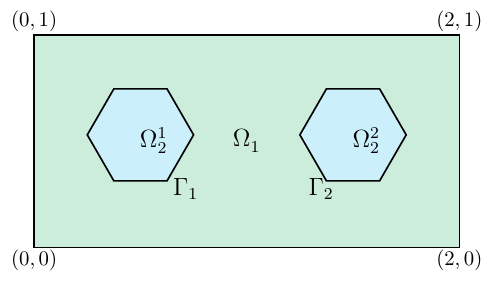}
    \caption{Domain used in Example 3 (immersed hexagons in a rectangular domain).}
    \label{fig:domain2}
\end{figure}

In this example, we prescribe the exact analytical solution as
$\varphi_1=\sin(\pi x)\sin(\pi y)$,  $\varphi_2=\sin(\pi x)\sin(\pi y)$, 
on a rectangular domain $\Omega = (0,2) \times (0,1)$. To test a more complex geometrical configuration, the solute domain $\Omega_{2}^{}$ comprises two disjoint regular hexagons, centred  at $(0.5,0.5)$ and $(1.5,0.5)$ with circumradii of $0.2$. These inclusions are   immersed within the solvent domain $\Omega_{1}^{}$, as illustrated in \Cref{fig:domain2}, so that $\partial\Omega_2\setminus\Gamma=\emptyset$ and no boundary condition is imposed on $\varphi_2$ (cf.\ \Cref{rem:domains}).  The model parameters are set to $\varepsilon_r=2$ and $\delta_c=1$ (the stress is the dimensionless tensor \eqref{mstress}, not involving  $\varepsilon_1$ or $l_c$ separately). 

\begin{figure}[t!]
    \centering
    \includegraphics[width=\linewidth]{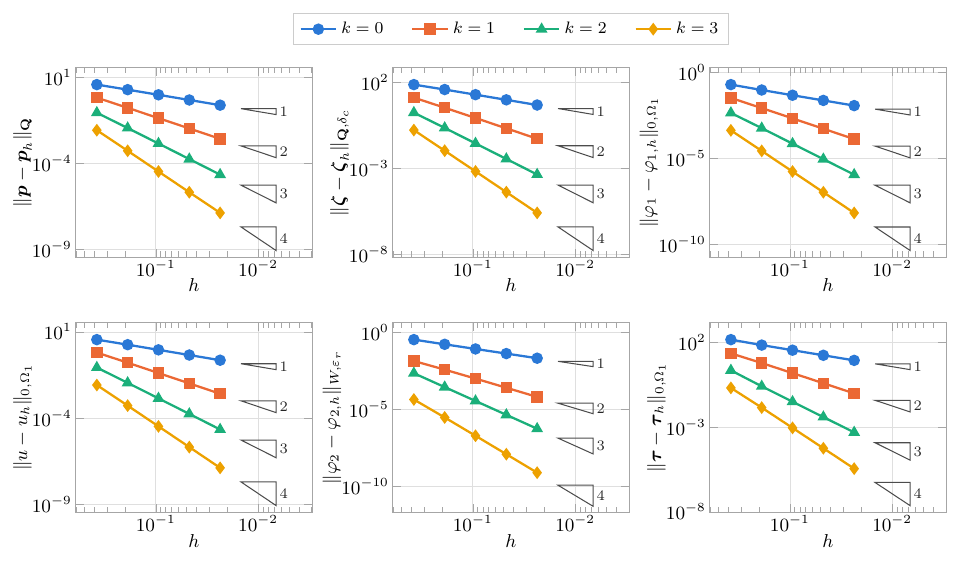}
    \vspace{-5mm}
    \caption{Example 3: errors of all field variables (in the norms of \Cref{sec:wellposedness}) and of the Maxwell stress tensor, for polynomial degrees $k\in\{0,1,2,3\}$. The triangles indicate the slope $k+1$.}
    \label{errfig:3}\label{fig:stresserr3}
\end{figure}

The primary objective of this configuration is to assess the  optimal convergence of the proposed scheme in the presence of disjoint, polygonal immersed interfaces. The error profiles for each variable, evaluated across different polynomial degrees $k \in \{0,1,2,3\}$, are presented in \Cref{errfig:3}. The log--log error curves align well with the $\mathcal{O}(h^{k+1})$ reference slopes, confirming that the numerical convergence rates remain optimal and conform strictly to the theoretical estimates of \Cref{errthm}, regardless of the topological complexity of the interface.

Furthermore, we explicitly evaluate the approximation of the Maxwell electrostatic stress tensor. 
As depicted in the last panel of \Cref{fig:stresserr3}, the convergence rate of the stress tensor in the $L^2$-norm is in  agreement with the theoretical bounds established in \Cref{stressthm}. This confirms that the mixed framework preserves high accuracy  approximations of derived physical quantities.

\subsection{Example 4: Two-dimensional convergence tests with immersed circular inclusions}\label{subsec:example04}

\begin{figure}[t!]
    \centering
    \includegraphics{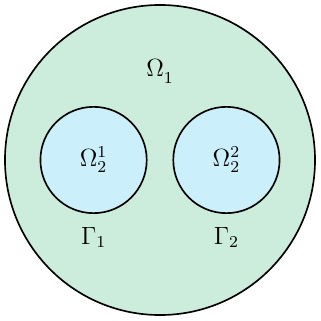}
    \caption{Domain in Example 4 (two immersed disks in a circular domain).}
    \label{fig:domain3}
  \end{figure}

We now prescribe the exact analytical solutions
$\varphi_1 = x^2(y^2-1)e^{x+y}$, $\varphi_2 = x(x-1)(y^2-1)e^{x+y}$, 
which explicitly exhibit a discontinuity (jump) across the interface. The global domain $\Omega$ is the unit disk centred at the origin, containing two disjoint circular solute inclusions of radius $0.2$, centred at $(-0.4,0)$ and $(0.4,0)$, respectively, as depicted in \Cref{fig:domain3}. The model parameters are set to $\varepsilon_r=2$ and $\delta_c=1$.

\begin{figure}[t!]
    \centering
    \begin{minipage}[b]{0.48\textwidth}
        \centering
        \begin{overpic}[width=0.8\textwidth]{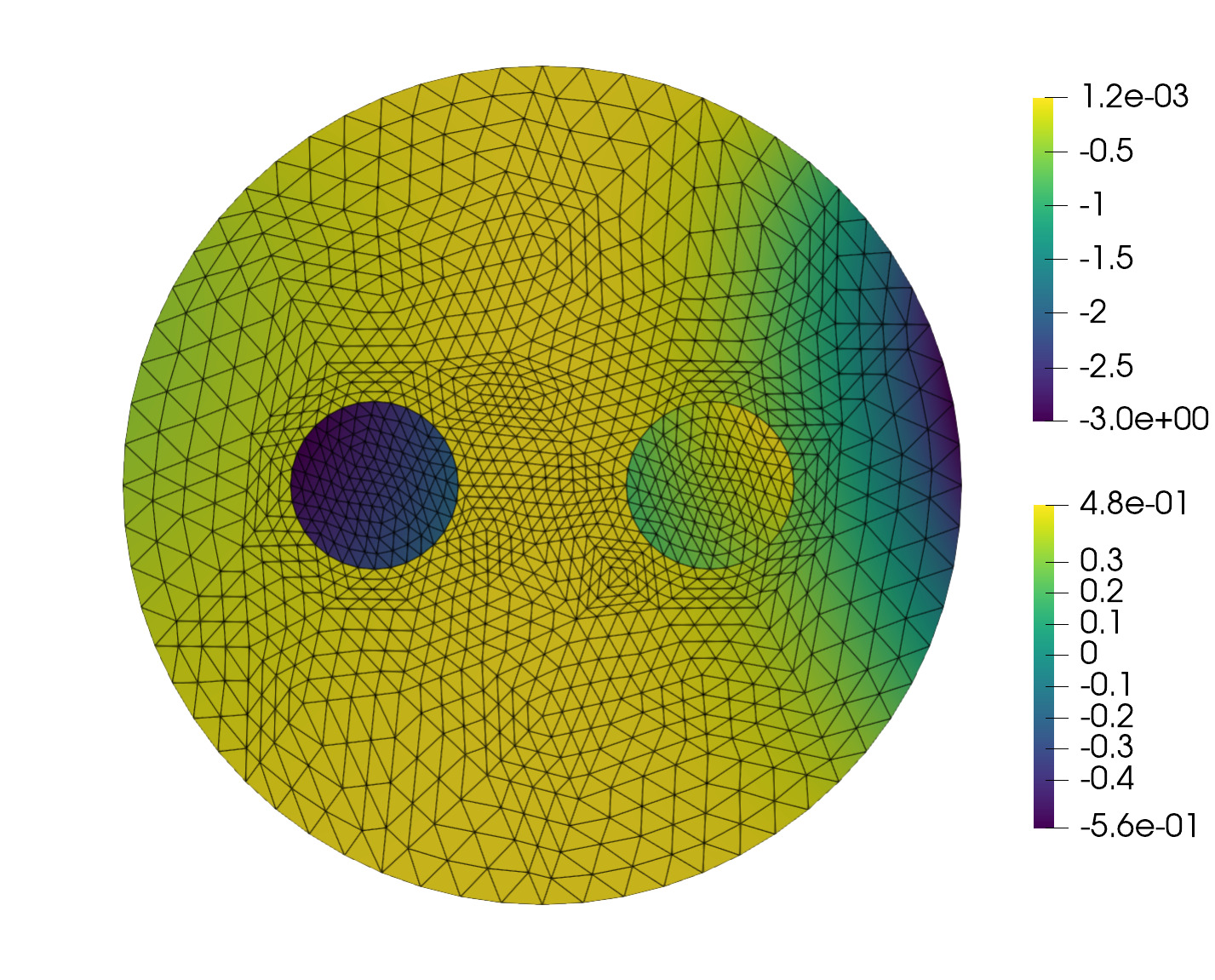}
            \put(95,60){$\varphi_{1,h}$}
            \put(95,25){$\varphi_{2,h}$}
        \end{overpic}
    \end{minipage}
    \hfill
    \begin{minipage}[b]{0.48\textwidth}
    \begin{overpic}[width=0.9\textwidth]{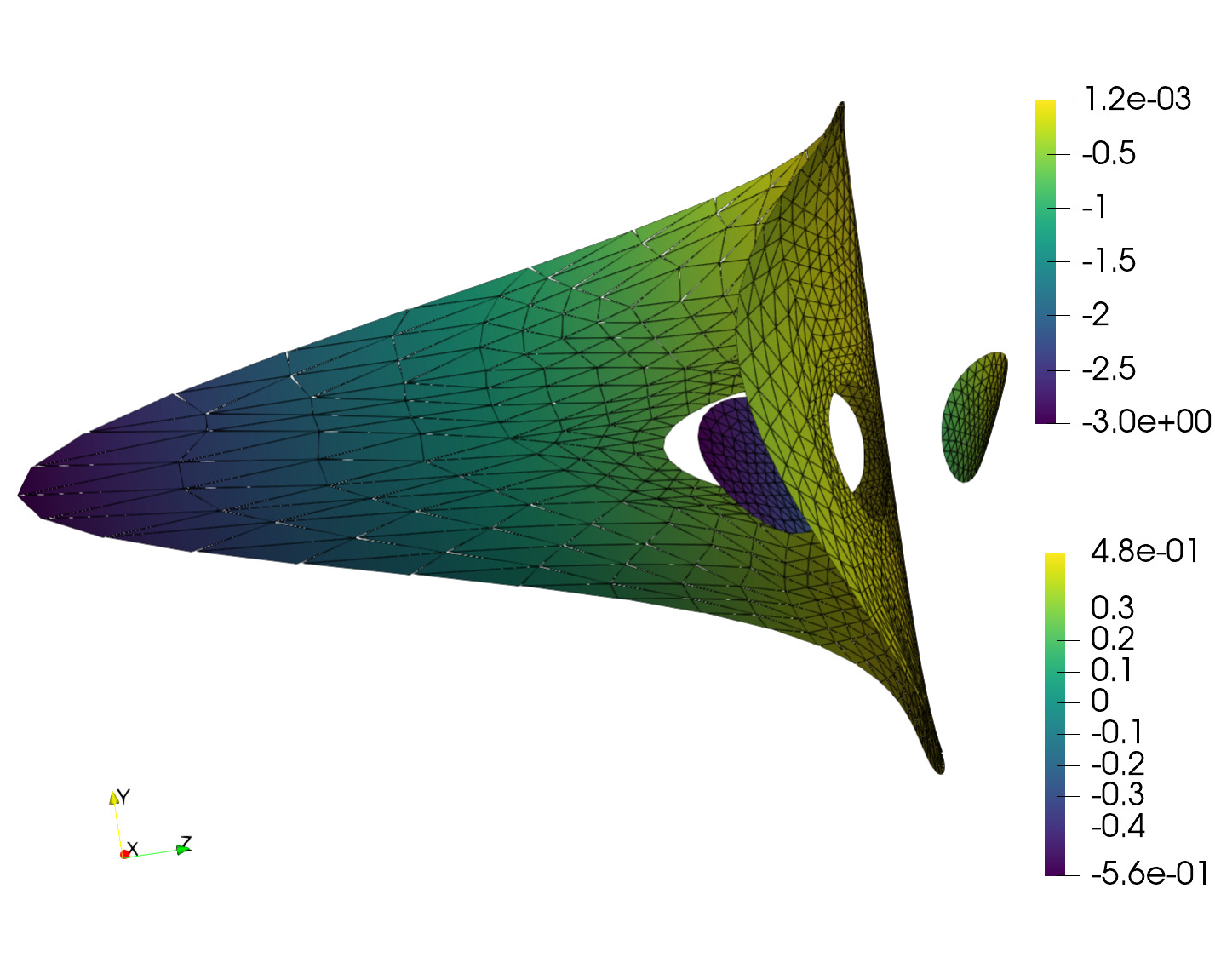}
            \put(95,60){$\varphi_{1,h}$}
            \put(95,25){$\varphi_{2,h}$}
        \end{overpic}
    \end{minipage}\\
    \begin{minipage}[b]{0.32\textwidth}
        \centering
        \begin{overpic}[width=0.99\textwidth]{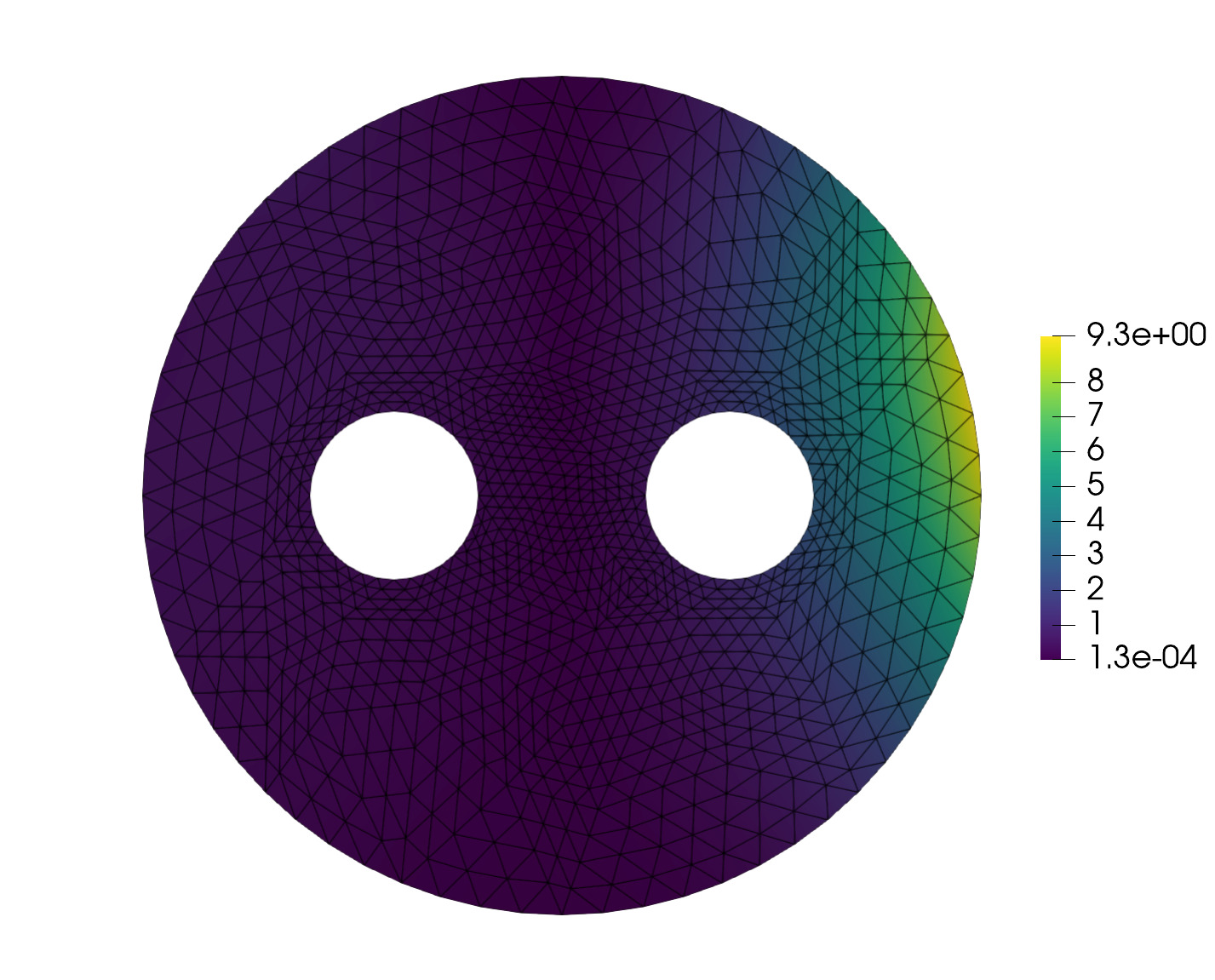}
            \put(95,40){$|\bp_{h}|$}
        \end{overpic}
    \end{minipage} 
    \begin{minipage}[b]{0.32\textwidth}
        \centering
        \begin{overpic}[width=0.99\textwidth]{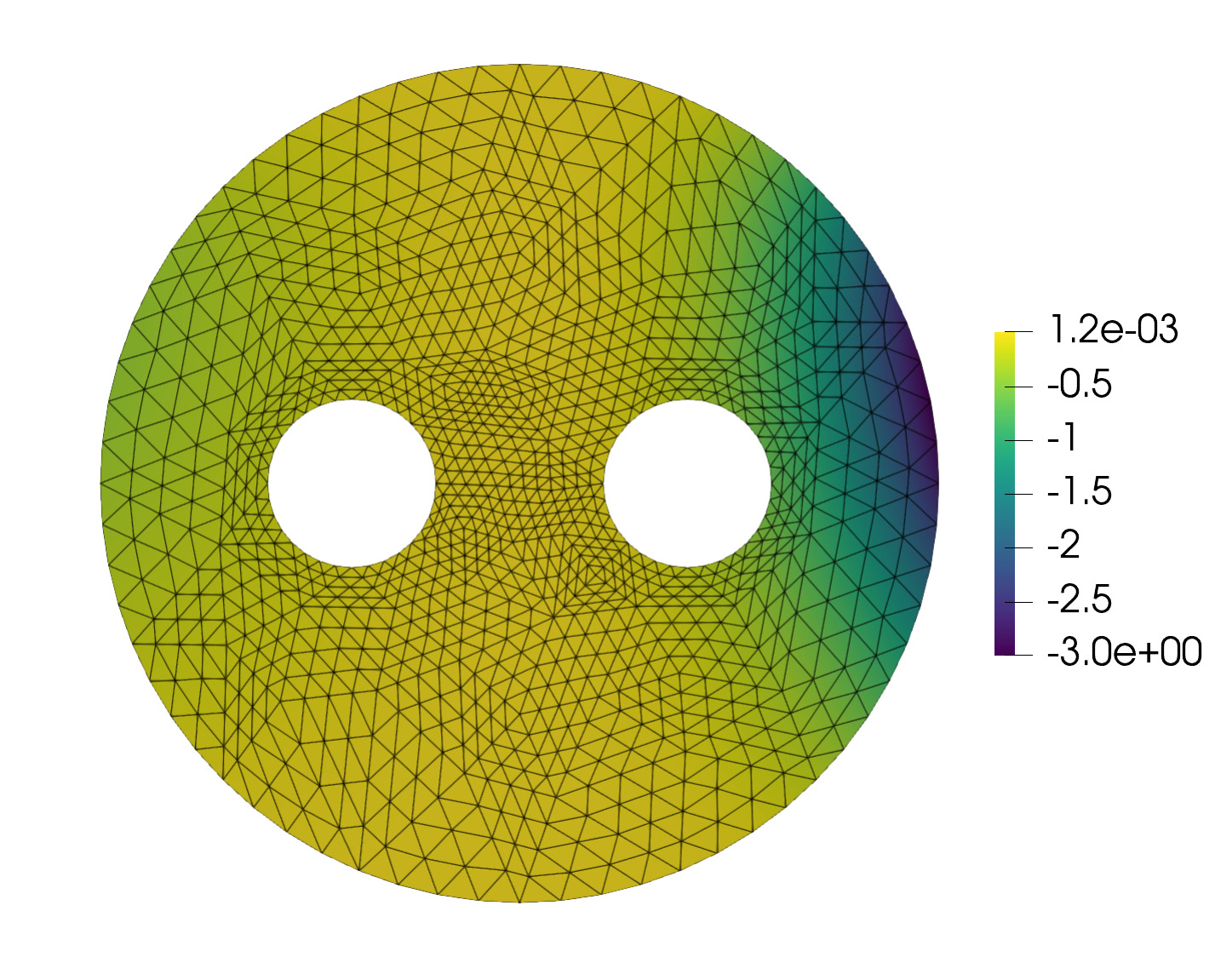}
            \put(95,40){$u_{h}^{}$}
        \end{overpic}
    \end{minipage} 
    \begin{minipage}[b]{0.32\textwidth}
        \centering
        \begin{overpic}[width=0.99\textwidth]{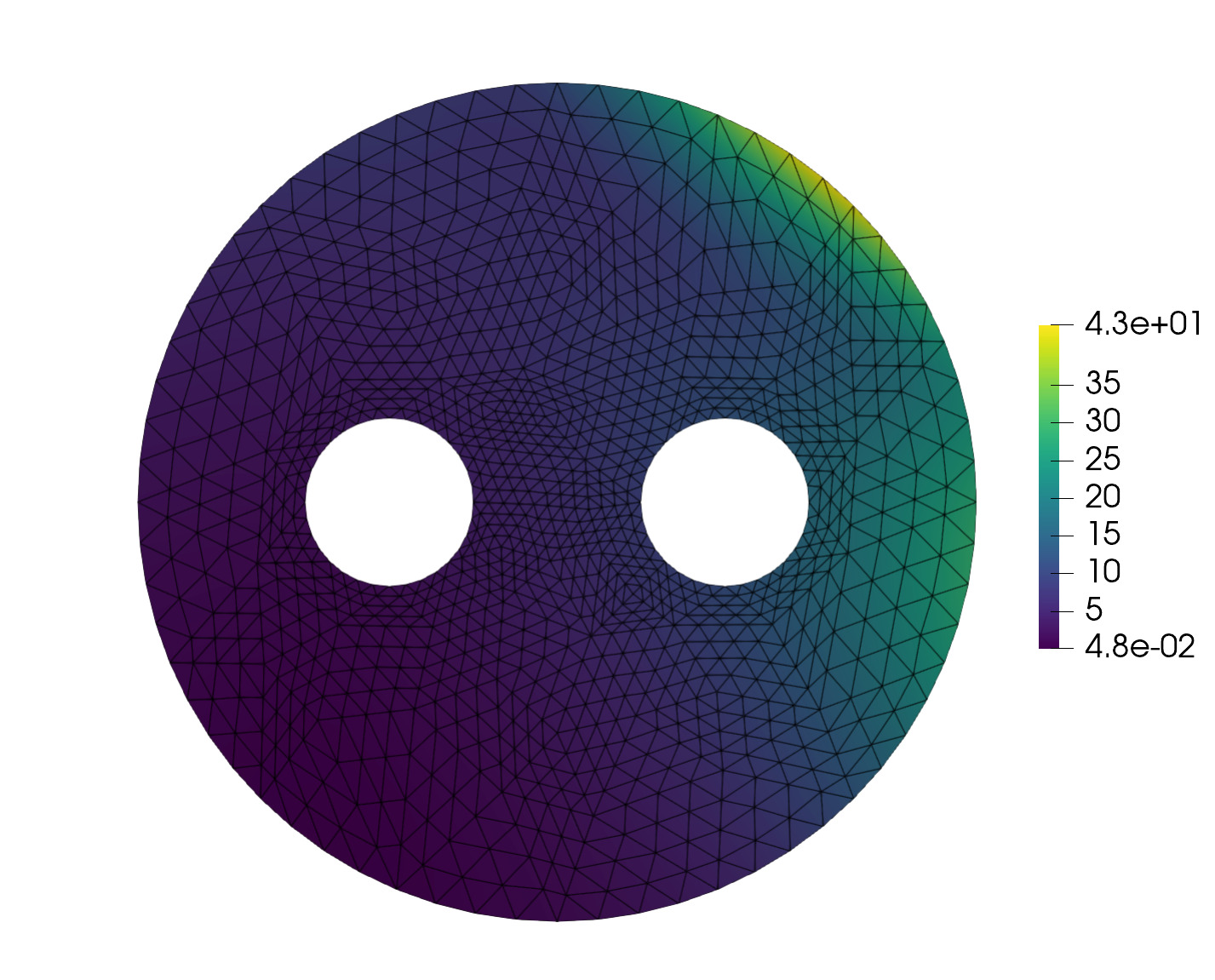}
            \put(95,40){$|\bzeta_{h}^{}|$}
        \end{overpic}
    \end{minipage}
    \caption{Plots of approximate solutions computed with polynomial degree $k=1$ at mesh refinement level 3 in Example 4: electrostatic potentials $\varphi_{1,h}^{},\varphi_{2,h}^{}$ (top left), elevation plot of the potentials (top right), magnitude of $\bp_{h}^{}$ (bottom left), $u_{h}^{}$ (bottom centre) and magnitude of $\bzeta_{h}^{}$ (bottom right).}
    \label{fig:soln4}
\end{figure}

The approximate solutions computed using a polynomial degree of $k=1$ are visualised in \Cref{fig:soln4}. To  highlight the interfacial jump in the electrostatic potential, an elevation plot (warp by scalar) of the approximate primal variables is also provided in the top-right panel. 

The   error metrics are compiled in \Cref{table:eg4},  confirming that the proposed scheme  maintains optimal convergence rates for all primary and mixed variables,   consistently with the theoretical bounds from \Cref{errthm}. The curved interfaces and the outer boundary are approximated by polygonal lines, and the exact solution and the interface and boundary data are evaluated on the resulting discrete domain $\Omega_h$. Hence, the reported errors measure the discretisation error only and do not include the geometric error. Furthermore, the numerical errors and corresponding convergence rates for the Maxwell electrostatic stress tensor, also reported in \Cref{table:eg41}, align with the theoretical predictions of \Cref{stressthm} also for curved boundaries and discontinuous potential fields.

\begin{table}[t!]
\centering
\resizebox{\textwidth}{!}{%
\begin{tabular}{cccccccccccc}
\toprule
order & $h$ & $\|\bp-\bp_h\|_{\bQ}$ & Rate & $\|\bzeta-\bzeta_h\|_{\bQ,\delta_c}$ & Rate & $\|\varphi_1-\varphi_{1,h}\|_{0,\Omega_1}$ & Rate & $\|u-u_h\|_{0,\Omega_1}$ & Rate & $\|\varphi_2-\varphi_{2,h}\|_{W,\varepsilon_r}$ & Rate \\
\midrule
       & 0.329 & 1.43e+00 & $\star$ & 1.59e+01 & $\star$ & 2.02e-01 & $\star$ & 1.17e+00 & $\star$ & 4.91e-02 & $\star$ \\
$k=0$  & 0.176 & 7.57e-01 & 1.016 & 8.54e+00 & 0.989 & 1.07e-01 & 1.007 & 6.11e-01 & 1.035 & 2.50e-02 & 1.073 \\
       & 0.088 & 3.84e-01 & 0.981 & 4.36e+00 & 0.973 & 5.46e-02 & 0.980 & 3.09e-01 & 0.988 & 1.26e-02 & 0.991 \\
       & 0.044 & 1.93e-01 & 0.996 & 2.19e+00 & 0.993 & 2.74e-02 & 0.995 & 1.55e-01 & 0.998 & 6.33e-03 & 0.998 \\
\midrule
       & 0.329 & 1.66e-01 & $\star$ & 1.81e+00 & $\star$ & 1.95e-02 & $\star$ & 1.44e-01 & $\star$ & 2.16e-03 & $\star$ \\
$k=1$  & 0.176 & 4.54e-02 & 2.064 & 4.88e-01 & 2.088 & 5.37e-03 & 2.053 & 3.92e-02 & 2.071 & 6.11e-04 & 2.006 \\
       & 0.088 & 1.17e-02 & 1.964 & 1.25e-01 & 1.971 & 1.38e-03 & 1.966 & 1.01e-02 & 1.967 & 1.60e-04 & 1.942 \\
       & 0.044 & 2.95e-03 & 1.990 & 3.15e-02 & 1.993 & 3.47e-04 & 1.992 & 2.54e-03 & 1.992 & 4.06e-05 & 1.980 \\
\midrule
       & 0.329 & 1.45e-02 & $\star$ & 1.26e-01 & $\star$ & 1.32e-03 & $\star$ & 1.32e-02 & $\star$ & 7.91e-05 & $\star$ \\
$k=2$  & 0.176 & 1.98e-03 & 3.168 & 1.70e-02 & 3.186 & 1.86e-04 & 3.117 & 1.80e-03 & 3.175 & 1.16e-05 & 3.053 \\
       & 0.088 & 2.55e-04 & 2.965 & 2.18e-03 & 2.971 & 2.41e-05 & 2.958 & 2.31e-04 & 2.967 & 1.53e-06 & 2.933 \\
       & 0.044 & 3.22e-05 & 2.991 & 2.74e-04 & 2.993 & 3.04e-06 & 2.990 & 2.91e-05 & 2.992 & 1.94e-07 & 2.982 \\
\midrule
       & 0.329 & 8.70e-04 & $\star$ & 6.14e-03 & $\star$ & 8.18e-05 & $\star$ & 8.10e-04 & $\star$ & 1.61e-06 & $\star$ \\
$k=3$  & 0.176 & 5.97e-05 & 4.265 & 4.18e-04 & 4.277 & 5.75e-06 & 4.224 & 5.55e-05 & 4.267 & 1.24e-07 & 4.087 \\
       & 0.088 & 3.85e-06 & 3.965 & 2.69e-05 & 3.970 & 3.74e-07 & 3.956 & 3.58e-06 & 3.966 & 8.23e-09 & 3.922 \\
       & 0.044 & 2.43e-07 & 3.992 & 1.69e-06 & 3.994 & 2.36e-08 & 3.990 & 2.26e-07 & 3.993 & 5.25e-10 & 3.978 \\
\bottomrule
\end{tabular}}
\caption{Errors and convergence rates in Example 4.}
\label{table:eg4}
\end{table}

\begin{table}[t!]
\centering 
{\begin{tabular}{ccccccccc}
\toprule
 & \multicolumn{2}{c}{$k=0$}
 & \multicolumn{2}{c}{$k=1$}
 & \multicolumn{2}{c}{$k=2$}
 & \multicolumn{2}{c}{$k=3$} \\
\cmidrule(lr){2-3} \cmidrule(lr){4-5} \cmidrule(lr){6-7} \cmidrule(lr){8-9}
$h$ & $\|\btau -\btau_{h}\|_{{\mathbb{L}^2(\Omega_{1}^{})}}$ & Rate & $\|\btau -\btau_{h}\|_{{\mathbb{L}^2(\Omega_{1}^{})}}$ & Rate & $\|\btau -\btau_{h}\|_{{\mathbb{L}^2(\Omega_{1}^{})}}$ & Rate &$\|\btau -\btau_{h}\|_{{\mathbb{L}^2(\Omega_{1}^{})}}$ & Rate \\
\midrule
0.329 & 4.4799e+01 & $\star$ & 4.8881e+00 & $\star$ & 3.3207e-01 & $\star$ & 1.6217e-02 & $\star$ \\
0.176 & 2.5250e+01 & 0.913 & 1.3813e+00 & 2.011 & 4.7899e-02 & 3.082 & 1.1760e-03 & 4.176 \\
0.088 & 1.3110e+01 & 0.948 & 3.6126e-01 & 1.940 & 6.2954e-03 & 2.936 & 7.7368e-05 & 3.937 \\
0.044 & 6.6291e+00 & 0.985 & 9.1805e-02 & 1.979 & 8.0048e-04 & 2.980 & 4.9210e-06 & 3.981 \\
\bottomrule
\end{tabular}}
\caption{Convergence results of Maxwell stress tensor in Example 4.}
\label{table:eg41}
\end{table}

\subsection{Example 5: A charged particle in an electrolyte, closed-form solution and robustness}\label{subsec:example05}

In this example and the next one we investigate the electrostatics of charged particles immersed in an electrolyte. While classical mean-field theories predict standard Coulomb's behaviour (like charges repel and opposite charges attract),  multivalent ions at sufficiently high concentrations can induce counterintuitive effects, such as the attraction of similarly charged particles. This   {charge reversal} (or {overcharging}), is driven by the strong accumulation of multivalent counter-ions near charged surfaces.

Within the solvent $\Omega_{1}$, the   potential $\varphi_{1}^{}$ is governed by the fourth-order modified \gls{pf} equation derived from  \cite{BSK2011}. In the solute regions $\Omega_{2}$ (the interior of the particles), the potential $\varphi_{2}^{}$ satisfies {the Poisson equation, since the particles are dielectric and contain no mobile ions (nor, here, fixed volume charges)}. The coupled dimensional system reads:
$$ \varepsilon_1(l_c^2 \Delta^2\varphi_1 - \Delta\varphi_1) = \rho(\varphi_1) \quad \text{in } \Omega_{1}^{}, \qquad
 -\varepsilon_2 \Delta\varphi_{2}^{} = 0 \quad \text{in } \Omega_{2}^{}. $$
Here, $\varepsilon_{1}$ and $\varepsilon_{2}$ denote the dielectric permittivities of the solvent and solute, respectively, and $l_c$ is the electrostatic correlation length. The charge density is given by the Boltzmann distribution 
 $\rho(\varphi_1) = z_{1} e C_{\infty}^{+}e^{-z_{1}\frac{e}{k_{B}T}\varphi_1} + z_{2} e C_{\infty}^{-}e^{-z_{2}\frac{e}{k_{B}T}\varphi_1}$,
where $z_{1}$ and $z_{2}$ are the valencies of the cations and anions, $C_{\infty}^{+}$ and $C_{\infty}^{-}$ are their respective far-field concentrations, $e$ is the elementary charge, $k_{B}$ is the Boltzmann constant, and $T$ is the absolute temperature. We assume a binary $z_{1}:z_{2}$ electrolyte satisfying the electroneutrality condition $z_{1}C_{\infty}^{+} = -z_{2}C_{\infty}^{-}$.
 This distribution neglects the steric (crowding) effects that are included in the Fermi-like distributions of \cite{BSK2011,Liu2013}.
 
To circumvent scaling issues and round-off errors associated with SI units at the nanometre scale, we non-dimensionalise the system using the Debye length $\lambda_D=\sqrt{\frac{\varepsilon_{1}k_{B}T}{e^2C_{\infty}^{+}z_{1}(z_{1}-z_{2})}}$. Introducing the dimensionless variables $\hat{\textbf{x}}=\frac{\textbf{x}}{\lambda_D}$, $\delta_{c}=\frac{l_c}{\lambda_D}$, and $\hat\varphi_{i}=\frac{e}{k_{B}T}\varphi_{i}$, we obtain the following dimensionless nonlinear interface problem \cite{Liu2025MixedPFInterface, Liu2025InterfacePF}:
\begin{align}\label{PP}
	\delta_{c}^{2}\Delta^{2}\hat\varphi_{1}^{}-\Delta\hat\varphi_{1}^{} &= \frac{e^{-z_{1}\hat\varphi_{1}^{}} - e^{-z_{2}\hat\varphi_{1}^{}}}{z_{1}-z_{2}^{}} && \text{ in }\Omega_{1},\nonumber\\
    -\varepsilon_{r}\Delta\hat\varphi_{2}^{}&=0&& \text{ in }\Omega_{2},\nonumber\\
    \hat\varphi_{1}^{}&=\hat\varphi_{2}^{} && \text{ on }\Gamma,\\
    \nabla\hat\varphi_{1}^{}\cdot \bm{n}_1+\varepsilon_{r}\nabla\hat\varphi_2\cdot\bm{n}_2 &=\hat\sigma:=\frac{e\lambda_{D}}{k_{B}T\varepsilon_{1}^{}}\sigma && \text{ on }\Gamma,\nonumber\\
    \hat\varphi_{1}^{}&=0 &&\text{ on }\partial\Omega_{1}^{}\setminus\Gamma ,\nonumber\\
    -\Delta\hat\varphi_{1}^{}&=\varphi_{0} && \text{ on }\partial\Omega_{1},\nonumber
\end{align}
where $\varepsilon_{r}=\varepsilon_2/\varepsilon_1$ is the relative permittivity and $\sigma$ is the surface charge density on the solute. For $i=1,2$, $\varepsilon_{i}^{}$ are defined by $\varepsilon_{i}^{}=\epsilon_{0}^{}\epsilon_{i}^{}$ . The cation concentration $C_{\infty}^{+}$ is related to the molar concentration $c_{0}$ (in mM) via $c_{0} = \frac{C_{\infty}^{+}}{N_A}$, where $N_{A}$ is Avogadro's number. {Note that $\hat u=-\Delta\hat\varphi_1$ satisfies $-\delta_c^2\Delta\hat u+\hat u=\hat\rho(\hat\varphi_1)$ in $\Omega_1$, by the first equation of \eqref{PP}; that is, $\hat u$ is the ionic charge density smoothed over the correlation length (and $\hat u=\hat\rho$ if $\delta_c=0$), in units of the bulk charge scale $ez_1(z_1-z_2)C_\infty^+$. The condition $-\Delta\hat\varphi_1=\varphi_0$ on $\Gamma$ thus prescribes this smoothed charge density at contact with the particle. Since no value for $\varphi_0$ is reported in \cite{Liu2025MixedPFInterface,Liu2025InterfacePF}, we take $\varphi_0=-\operatorname{sign}(\sigma)$, that is, $\varphi_0=1$ for the negative surface charge of \Cref{table:eg5}: a contact charge density of the order of the bulk charge scale, with the sign of the counter-ions (for $\sigma<0$, the value that the local linearised relation $\hat\rho=-\hat\varphi_1$ would give for a contact potential of $-k_{B}T/e$). The datum $\varphi_0$ acts as an additional interfacial source, and the forces in Example 6 depend on it.} All   parameters employed in Examples 5-6 are detailed in \Cref{table:eg5}. For the three concentrations considered below, $c_0={7.04}$, $15.85$ and ${63.4}$\,mM, one has $\lambda_D={1.5}$, $1$ and ${0.5}$\,nm, hence $\delta_c={2/3}$, $1$ and ${2}$, and $\hat\sigma={-4.09}$, ${-}2.73$ and ${-1.36}$, respectively. {We restrict the concentrations to $c_0\le63.4$\,mM, for which the Debye length $\lambda_D\ge0.5$\,nm is not smaller than the radius of a hydrated trivalent cation (about $0.45$\,nm for La$^{3+}$ \cite{Nightingale1959}), and $l_c=1$\,nm is of the order of its diameter. For much larger concentrations, such as $c_0=1585$\,mM ($\lambda_D=0.1$\,nm), the Debye length would be smaller than the ions themselves, and a continuum description of the electrolyte is no longer justified. The lowest concentration is chosen so that $\delta_c>1/2$, which is required for overscreening (see below).}

\begin{table}[t!]
\centering
\begin{tabular}{clcc}
\toprule
\textbf{Notation} & \textbf{Physical Definition} & \textbf{Value} & \textbf{Unit} \\
\midrule
$R$ & Radius of the electrolyte domain & 10 & nm \\
$a$ & Radius of charged particles & 1 & nm \\
$e$ & Elementary charge & $1.6 \times 10^{-19}$ & C \\
$k_{B} T$ & Boltzmann energy & $4.14 \times 10^{-21}$ & J (or N\,m) \\
$\epsilon_{0}$ & Permittivity of vacuum & $8.85 \times 10^{-12}$ & C$^2$/(N\,m$^2$) \\
$\epsilon_{1}$ & Electrolyte dielectric constant & 80 & -- \\
$\epsilon_{2}$ & Particle dielectric constant & 2 & -- \\
$z_1$ & Valency of cation & 3 & -- \\
$z_2$ & Valency of anion & $-1$ & -- \\
$l_c$ & Electrostatic correlation length & 1 & nm \\
$\sigma$ & Surface charge & {$-0.05$} & C/m$^2$ \\
$N_{A}$ & Avogadro constant & $6.02 \times 10^{23}$ & mol$^{-1}$ \\
$c_0$ & Far-field cation concentration & {7.04}, 15.85, {63.4} & mM \\
\bottomrule
\end{tabular}
\caption{Model parameters used in Examples 5 and 6.}
\label{table:eg5}
\end{table}

Focusing on a $3:-1$ electrolyte (i.e., $z_1=3, z_2=-1$), the dimensionless charge density reduces to $\hat\rho(\hat\varphi_1)=\frac{1}{4}(e^{-3\hat\varphi_{1}^{}} - e^{\hat\varphi_{1}^{}})$. Linearising this expression around $\hat{\varphi}_{1}=0$ yields the approximation $-\hat{\varphi}_{1}$. Substituting this linearisation into \eqref{PP} produces a reduced linear model that reads exactly as \eqref{PP} except for the first equation,
\begin{align}\label{LP}
	\delta_{c}^{2}\Delta^{2}\hat\varphi_{1}^{}-\Delta\hat\varphi_{1}^{} +\hat\varphi_{1}^{} &= 0\quad \text{ in }\Omega_{1},\\
\nonumber (\text{and the other equations are as in \eqref{PP}}). &\end{align}
This system corresponds to \eqref{IP} {without the zeroth-order term in (\ref{IP})$_2$ ($a_3(\varphi_2,\psi_2)=\varepsilon_r\int_{\Omega_2}\nabla\varphi_2\cdot\nabla\psi_2$ in \eqref{IPw}). Since $\Omega_2$ is immersed, this form vanishes for functions that are constant in each particle; the constants are determined through the coupling with $\Omega_1$, and the discrete problems are uniquely solvable.} {For the negative surface charge of \Cref{table:eg5}, the counter-ions are the trivalent cations, which is the configuration in which charge reversal is expected \cite{Besteman2005}. The nonlinear model distinguishes counter-ions from co-ions, whereas the linearised one   is invariant under a simultaneous change of sign of $\hat\sigma$ and $\varphi_0$, so that, with $\varphi_0=-\operatorname{sign}(\sigma)$, its solution simply changes sign with $\sigma$.}

\begin{figure}[t!]
    \centering
    \includegraphics[width=0.9\linewidth]{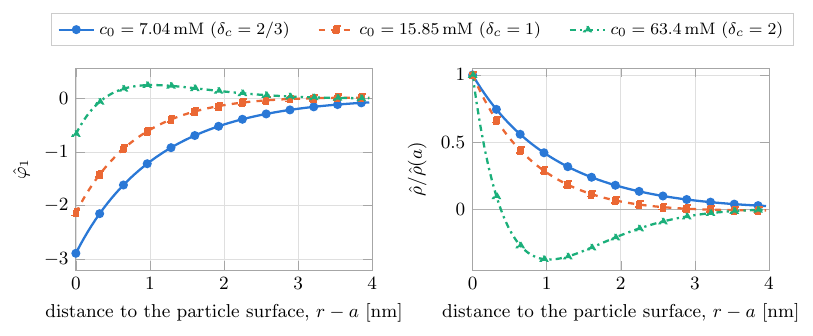}
    \vspace{-3mm}
    \caption{Example 5: closed-form radial profiles of the potential $\hat\varphi_1$ (left) and of the normalised ionic charge density {$\hat\rho/\hat\rho(a)$} (right) around a single particle, for the three concentrations{ and $\sigma=-0.05$\,C/m$^2$}. The sign change of $\hat\rho$ represents  overscreening.}
    \label{fig:eg5profiles}
\end{figure}

Consider a single particle $\Omega_2=\{|\hat\bx|<a\}$ inside the electrolyte $\Omega_1=\{a<|\hat\bx|<R\}$. Since $\delta_c^2\Delta^2-\Delta+1=\delta_c^2(\Delta-\lambda_+)(\Delta-\lambda_-)$ with $\lambda_\pm=\frac12\delta_c^{-2}[ 1\pm\sqrt{1-4\delta_c^2}]$, 
the radial solution of \eqref{LP} is
\begin{equation*}
    \hat\varphi_1(r)=\sum_{s=\pm}\big(A_sK_0(\sqrt{\lambda_s}\,r)+B_sI_0(\sqrt{\lambda_s}\,r)\big),\qquad \hat\varphi_2(r)=C,
\end{equation*}
where $I_0,K_0$ are the modified Bessel functions {(the only radial harmonic functions that are regular at $r=0$ are the constants, so that the flux condition on $\Gamma$ reduces to $-\hat\varphi_1'(a)=\hat\sigma$)}, and the five constants are determined by the three interface conditions and the two outer boundary conditions in \eqref{PP}. For $\delta_c>1/2$, which is the case for all the concentrations considered here, $\lambda_\pm$ are complex conjugate, and the potential and the ionic charge density $\hat\rho=-\hat\varphi_1$ decay in an oscillatory manner. This is the overscreening mechanism of the \gls{bsk} theory behind charge reversal \cite{BSK2011}. \Cref{fig:eg5profiles} shows these radial profiles for the three concentrations{. For the three concentrations the potential is negative near the particle, and the ionic charge density changes sign at a distance of about $8.4$, $3.2$ and $0.4$\,nm from the particle surface for $c_0=7.04$, $15.85$ and $63.4$\,mM, respectively (with a very small amplitude for the two lower concentrations, which is not visible in the figure for $c_0=7.04$\,mM). In all cases $\hat\rho$ decays in an oscillatory manner, and the oscillations become more pronounced as $\delta_c$ increases.}

\begin{figure}[t!]
    \centering
    \includegraphics[width=\linewidth]{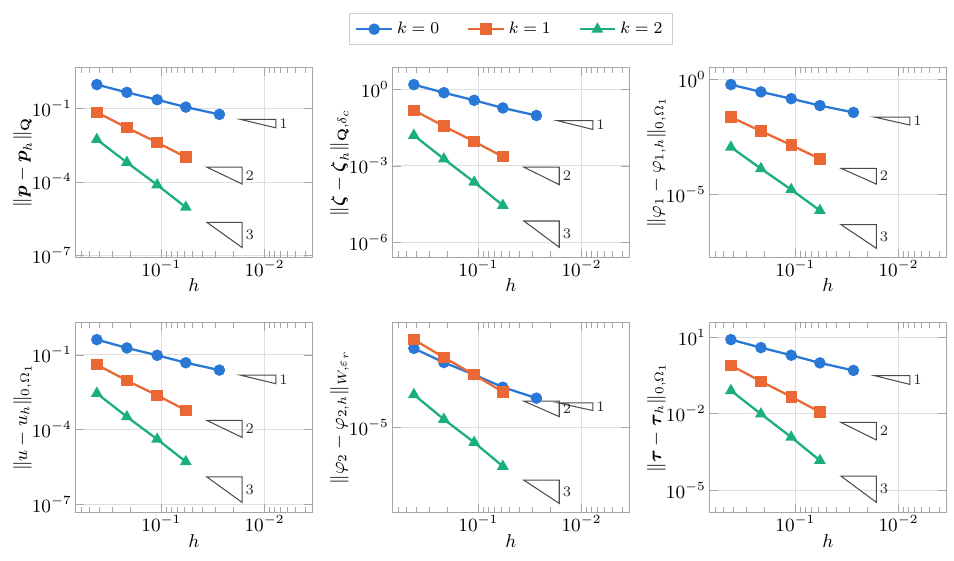}
    \vspace{-5mm}
    \caption{Example 5: errors with respect to the closed-form solution for $\delta_c=1$, $\varepsilon_r=2/80$ and $k\in\{0,1,2\}$. The triangles indicate the slope $k+1$.}
    \label{fig:eg5conv}
\end{figure}

We use this closed-form solution as a benchmark, with $a=1$, $R=4$, $\varepsilon_r=2/80$, $\hat\sigma=2.73$, and $\delta_c\in\{1,2,10\}${ (the value $\delta_c=10$ lies outside the physical range discussed above, and is included only to test the robustness of the method)}, on a sequence of unstructured meshes of the polygonal approximation of the annulus and the disk, evaluating the exact solution and the interface and boundary data on the discrete domain as in Example 4. This tests the transmission conditions with a realistic, small permittivity ratio. \Cref{fig:eg5conv} shows the errors for $\delta_c=1$ and $k\in\{0,1,2\}$: all variables and the Maxwell stress converge with the optimal order $\mathcal O(h^{k+1})$ {(the error of $\varphi_2$ even converges faster, since the exact $\hat\varphi_2$ is constant and its error is only due to the coupling with $\Omega_1$; for $k=1$ it is larger than for $k=0$ on the coarser meshes, owing to a discrete interface component of $\varphi_{2,h}$ that is not seen by the solvent, scales as $\varepsilon_r^{-1}$ and decays as $\mathcal O(h^{k+2})$)}. \Cref{table:eg5robust} compares the errors for the three values of $\delta_c$ on the same mesh: the errors in the weighted norms of \Cref{sec:wellposedness} are of the same size for $\delta_c=1$, $2$ and $10$, in agreement with the parameter-robust estimates of \Cref{errthm} (the stress error grows with $\delta_c$ because the stress itself grows with $\delta_c$). 

\begin{table}[t!]
\centering
\begin{tabular}{ccccccc}
\toprule
$\delta_c$ & $\|\bp-\bp_h\|_{\bQ}$ & $\|\bzeta-\bzeta_h\|_{\bQ,\delta_c}$ & $\|\varphi_1-\varphi_{1,h}\|_{0,\Omega_1}$ & $\|u-u_h\|_{0,\Omega_1}$ & $\|\varphi_2-\varphi_{2,h}\|_{W,\varepsilon_r}$ & $\|\btau-\btau_h\|_{0,\Omega_1}$ \\
\midrule
{$1$} & {1.01e-03 (2.16)} & {2.24e-03 (2.15)} & {3.56e-04 (2.17)} & {5.90e-04 (2.16)} & {3.59e-04 (2.72)} & {1.16e-02 (2.15)} \\
{$2$} & {8.48e-04 (2.15)} & {1.95e-03 (2.15)} & {3.39e-04 (2.18)} & {2.92e-04 (2.16)} & {3.59e-04 (2.72)} & {2.08e-02 (2.15)} \\
{$10$} & {7.89e-04 (2.15)} & {2.93e-03 (2.15)} & {3.35e-04 (2.18)} & {9.79e-05 (2.17)} & {3.59e-04 (2.72)} & {1.70e-01 (2.15)} \\
\bottomrule
\end{tabular}
\caption{Example 5: errors (and experimental rates with respect to the previous mesh) for $k=1$, $\varepsilon_r=2/80$, $h=0.058$, and three values of the correlation length.}
\label{table:eg5robust}
\end{table}

\subsection{Example 6: Interaction force between two charged particles}\label{subsec:example06}

We now place two  particles of radius $r=1$\,nm, centred at $(\pm(r+dL/2),0)$, so that their surface-to-surface distance is $dL$, inside the disk of radius $R=10$\,nm (see \Cref{fig:domain3}), and solve both the linearised model \eqref{LP} and the nonlinear model \eqref{PP} (the latter with Newton's method, initialised with the linear solution){. As in Example 5, there are no ions inside the particles, so that the solute equation is $-\varepsilon_r\Delta\hat\varphi_2=0$. The particles carry opposite surface charges, $\sigma=0.05$\,C/m$^2$ on the left particle and $\sigma=-0.05$\,C/m$^2$ on the right one, and, as in Example 5, $\varphi_0=-\operatorname{sign}(\sigma)$ on each interface, that is, $\varphi_0=-1$ on $\Gamma_1$ and $\varphi_0=1$ on $\Gamma_2$, so that the prescribed contact charge density has the sign of the counter-ions of each particle}. The force exerted by the electrolyte on the left particle is computed from the total stress (the Maxwell stress \eqref{mstress} minus the osmotic pressure):
\[
\hat{\bm F}_1=\int_{\Gamma_1}\big(\hat\btau-\hat\Pi\,\mathbb{I}\big)\bm n_{2}\ds,
\]
where $\bm n_2=-\bm n_1$ is the outward normal of the particle, $\hat\btau$ is the dimensionless stress \eqref{mstress} evaluated at the solution $(\hat\bp,\hat\bzeta,\hat u)$ of the dimensionless problem, and $\hat\Pi=\hat\varphi_1^2/2$ for the linearised model and $\hat\Pi=\frac{1}{12}\big(e^{-3\hat\varphi_1}-1+3(e^{\hat\varphi_1}-1)\big)$ for the nonlinear one. 
{The osmotic term is needed because, by \eqref{PP}, $\vdiv\hat\btau=-\hat\rho\,\nabla\hat\varphi_1$ in $\Omega_1$: at equilibrium, this electric body force on the ions is balanced by the solvent pressure, which also acts on the particle. Since $\nabla\hat\Pi=-\hat\rho\,\nabla\hat\varphi_1$, the total stress satisfies $\vdiv(\hat\btau-\hat\Pi\,\mathbb{I})=\bm 0$, so the force does not depend on the integration contour; for $\delta_c=0$ it reduces to the classical Poisson--Boltzmann force (Maxwell plus osmotic stress) \cite{Neu1999}, and its extension to the \gls{bsk} model follows from the momentum balance in \cite{BSK2012}.
}
Since the computations are two-dimensional, the physical force is a force per unit length of two parallel cylinders, $\bm F_1=\varepsilon_1\frac{(k_{B}T/e)^2}{\lambda_D}\hat{\bm F}_1$ (the factor $\varepsilon_1(k_{B}T/(e\lambda_D))^2$ relating $\btau^{\rm phys}$ and $\hat\btau$ in \Cref{sec:error}, times $\lambda_D$ from the line element; the osmotic pressure is scaled in the same way), which we report in pN/nm.

We verified the computations as follows: (i) an isolated particle feels a force of order $10^{-5}$; (ii) $\hat{\bm F}_1+\hat{\bm F}_2\approx\bm 0${ (for the nonlinear model, $|\hat F_{1,x}+\hat F_{2,x}|\le0.08$)}; (iii) the force computed on an auxiliary circle of radius $r+dL/4$ around the left particle agrees with the one computed on $\Gamma_1$; (iv) {doubling the domain radius to $R=20$\,nm does not change the forces appreciably; and (v) the forces are stable under refinement of the mesh near the particles. For $c_0=7.04$\,mM the surface potentials are the largest, and the nonlinear solution has a steep layer at the negatively charged particle, where the cation density $e^{-3\hat\varphi_1}$ is large; the reported values for this concentration are therefore computed on the mesh refined once more near the particles.}

\begin{figure}[t!]
    \centering
    \includegraphics[width=\linewidth]{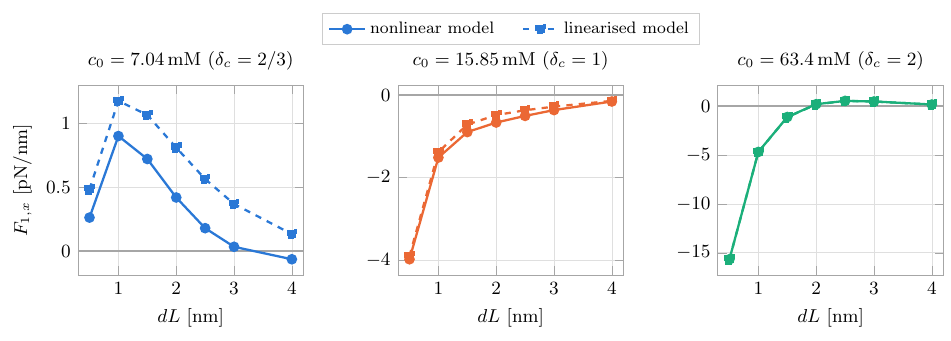}
    \vspace{-5mm}
    \caption{Example 6: horizontal force per unit length on the left particle (positive: attraction) versus the distance between the particles, for the nonlinear (solid) and linearised (dashed) models.{ The surface charges are $\sigma=0.05$\,C/m$^2$ (left particle) and $\sigma=-0.05$\,C/m$^2$ (right particle), and $\varphi_0=-\operatorname{sign}(\sigma)$.}}
    \label{fig:eg6}
\end{figure}

The forces are reported in \Cref{tab:eg6forces} and \Cref{fig:eg6}, where positive values mean attraction. {For $c_0=7.04$\,mM the particles attract each other, as for oppositely charged particles in the classical \gls{pb} theory, except for a very weak repulsion predicted by the nonlinear model at $dL=4$\,nm; the attraction is maximal at $dL\approx1$\,nm and decreases at shorter distances. For $c_0=15.85$\,mM the force is repulsive at all the distances considered, and it decays monotonically with $dL$. For $c_0=63.4$\,mM the interaction changes sign with the distance: the particles repel each other for $dL\le1.5$\,nm and attract each other weakly beyond (with a maximum of about $0.53$\,pN/nm at $dL=2.5$\,nm). A repulsion between oppositely charged particles, and a force that changes sign with the distance, cannot be captured by the classical \gls{pb} theory, for which oppositely charged particles attract each other. Here, the datum $\varphi_0$ imposes at each interface a layer of counter-ions with a charge density of the bulk scale, which remains in the gap as the particles approach, and the correlations make the charge density around each particle oscillate (see \Cref{fig:eg5profiles}); the sign change of the force for $c_0=63.4$\,mM occurs at a distance comparable with that of the first sign change of $\hat\rho$ in Example 5. The force therefore depends on $\varphi_0$: with $\varphi_0=1$ on both interfaces, for instance, its sign is reversed at most distances for $c_0=63.4$\,mM. The linearised and nonlinear models are close for $c_0=15.85$ and $63.4$\,mM, where the surface potentials are moderate ($\hat\varphi_1\approx1.6$--$2.3$ and $0.7$--$0.8$ on $\Gamma_1$, respectively), and differ for $c_0=7.04$\,mM, where they are the largest ($\hat\varphi_1\approx2$--$3$ on $\Gamma_1$), the nonlinear attraction being weaker.} 

\begin{table}[t!]
\centering
\begin{tabular}{ccccccc}
\toprule
 & \multicolumn{2}{c}{$c_0={7.04}$\,mM} & \multicolumn{2}{c}{$c_0={15.85}$\,mM} & \multicolumn{2}{c}{$c_0={63.4}$\,mM} \\
\cmidrule(lr){2-3} \cmidrule(lr){4-5} \cmidrule(lr){6-7}
$dL$ [nm] & linear & nonlinear & linear & nonlinear & linear & nonlinear \\
\midrule
0.5 & {0.4841} & {0.2622} & {-3.907} & {-3.986} & {-15.67} & {-15.66} \\
1.0 & {1.176} & {0.901} & {-1.384} & {-1.517} & {-4.692} & {-4.689} \\
1.5 & {1.066} & {0.7212} & {-0.7171} & {-0.8993} & {-1.137} & {-1.121} \\
2.0 & {0.8101} & {0.4197} & {-0.4863} & {-0.667} & {0.1794} & {0.1945} \\
2.5 & {0.5641} & {0.1794} & {-0.3691} & {-0.5067} & {0.5244} & {0.533} \\
3.0 & {0.3688} & {0.03348} & {-0.281} & {-0.3653} & {0.4675} & {0.4722} \\
4.0 & {0.1331} & {-0.06418} & {-0.1418} & {-0.157} & {0.1597} & {0.1611} \\
\bottomrule
\end{tabular}
\caption{Example 6: horizontal force per unit length (pN/nm) on the left particle{, for $\sigma=0.05$\,C/m$^2$ and $\varphi_0=-1$ on the left particle and $\sigma=-0.05$\,C/m$^2$ and $\varphi_0=1$ on the right one}; positive values mean attraction.}
\label{tab:eg6forces}
\end{table}

\section{Conclusion and future work}\label{sec:conclusion}
In this paper, we introduce and analyse a linearised fourth-order/second-order \gls{pf} interface problem arising from the \gls{bsk} theory \cite{BSK2011}. This model is essential for rigorously studying electrostatic interactions while accounting for electrostatic correlations. To solve the governing equations, we develop a mixed--primal \gls{fem}. While our approach employs the same four-variable splitting for the fourth-order solvent problem as in \cite{Liu2025MixedPFInterface}, we   circumvent the need for artificial stabilisation by utilising stable Raviart--Thomas finite element approximations for the vector fields $\bp$ and $\bzeta=-\delta_c^2\nabla(\vdiv\bp)$. We establish the well-posedness and parameter-robust stability of this formulation, derive optimal a priori error estimates in parameter-weighted norms, and validate the theoretical convergence rates through numerical experiments. We also show that the error estimates are robust with respect to the correlation length and to the permittivity ratio, up to a loss of at most half an order in the pre-asymptotic regime $h>\varepsilon_r^{1/2}$, and we confirm this behaviour numerically for $\delta_c\in[10^{-6},10]$ and $\varepsilon_r\in[10^{-6},1]$; for small $\delta_c$, this requires solving for the rescaled flux $\delta_c^{-2}\bzeta_h$ in order to avoid round-off errors in the linear solver. A closed-form solution for a single charged particle confirms the optimal convergence and its robustness with respect to the correlation length for a realistic, small permittivity ratio. Finally, we compute the interaction force between two charged particles with the full (Maxwell and osmotic) stress, which shows {that, as the concentration (and with it the correlation length relative to the Debye length) increases, the attraction between oppositely charged particles turns into a repulsion at short distances, and the interaction changes sign with the distance}. 

{We have identified a number of directions for} future research. First, addressing the second-order solute problem in primal form compromises local mass conservation and flux accuracy (particularly across interfaces with highly contrasting dielectric constants) while requiring at least   $H^1$ regularity. We will also recast the second-order problem in mixed form, but that will   necessitate Lagrange multipliers or Nitsche penalty terms to enforce the transmission conditions. 

We also intend to extend our   framework to the fully nonlinear fourth-order/second-order \gls{pf} interface problem. Depending on the physical scales of interest, this nonlinear framework will be coupled with either an underlying fluid flow (e.g., via the Stokes or Navier--Stokes equations \cite{HeSun2021}) or a poromechanical medium to model complex electrokinetic transport phenomena in deformable tissues \cite{gatica25}. Furthermore, we aim to explore more mathematically challenging interfacial mechanics, such as low-regularity interface conditions and alternative transmission mechanisms such as ion channels with selective filters \cite{keller2025model}, or similar processes that frequently arise in highly heterogeneous biological environments. Another challenge is the treatment of very irregular geometrical interfaces as those encountered in electrodiffusive interaction between neurons and surrounding tissue \cite{herlyng25}. 

Finally, investigating alternative advanced spatial discretisations for the fourth-order problem represents a direction of marked interest. $C^0$-interior penalty  schemes \cite{chen2025immersedc0interiorpenalty}, virtual element methods \cite{ma2025nitsche}, and hybrid high-order  formulations \cite{dong2021hybrid} could provide enhanced flexibility in handling general polygonal or polyhedral meshes and complex interface geometries without relying on sub-triangulations. {In addition, all the computations in this paper use a sparse direct solver. Since the inf-sup constant of \Cref{th:dwellposed} does not depend on $h$, $\delta_c$ or $\varepsilon_r$, and the continuity constant of \Cref{continuity} depends only on $\varepsilon_r$ (through $C_\Gamma$), the Riesz map associated with the parameter-weighted inner product inducing $\normiii{\bullet}$ (block diagonal, with $\bH(\mathrm{div})$-type blocks for $\bp$ and $\bzeta$, mass matrices for $\varphi_1$ and $u$, and a weighted $H^1$-block for $\varphi_2$) is a natural candidate for a block-diagonal preconditioner that is robust with respect to $h$ and $\delta_c$, following the operator preconditioning framework of \cite{MardalWinther2011}; robustness with respect to $\varepsilon_r$ would require a modification of the interface part of the norm.}




\bibliographystyle{siam}
\bibliography{ref}

\end{document}